
\input amstex
\documentstyle{amsppt}

\loadbold
\loadeurm
\loadeurb
\loadeusm

\magnification=\magstep1
\hsize=6.7truein
\vsize=9.5truein
\hcorrection{-0.1truein} 
\vcorrection{-0.2truein} 

\font\sc=cmcsc10
  
\font \smallrm=cmr10 at 11truept
 at 7truept
\font \smallbf=cmbx10 at 11truept 
 at 11truept
\font \smallsl=cmsl10 at 11truept

\baselineskip=14pt

\def \loongrightarrow {\relbar\joinrel\relbar\joinrel\rightarrow}
\def \llongrightarrow
 {\relbar\joinrel\relbar\joinrel\relbar\joinrel\rightarrow}
\def \llongtwoheadrightarrow
 {\relbar\joinrel\relbar\joinrel\relbar\joinrel\twoheadrightarrow}
\def \longtwoheadrightarrow
 {\relbar\joinrel\relbar\joinrel\twoheadrightarrow}
\def \longhookrightarrow {\lhook\joinrel\relbar\joinrel\rightarrow}
\def \llonghookrightarrow
 {\lhook\joinrel\relbar\joinrel\relbar\joinrel\relbar\joinrel\rightarrow}

\def \gerg {\frak{g}}

\def \gerb {\frak{b}}
\def \gerU {\frak{U}}
\def \gerF {\frak{F}}
\def \calU {\Cal{U}}
\def \calF {\Cal{F}}
\def \gersl {\frak{sl}}
\def \gergl {\frak{gl}}
\def \H {\hbox{\bf H}}
\def \gerH {\frak{H}}
\def \calH {\Cal{H}}

\def \uqg {U_q(\gerg)}
\def \geruqg {\gerU_q(\gerg)}
\def \caluqg {\calU_q(\gerg)}
\def \gerueg {\gerU_{\,\varepsilon}(\gerg)}
\def \calueg {\calU_\varepsilon(\gerg)}
\def \uqgl {U_q(\gergl_n)}
\def \geruqgl {\gerU_q(\gergl_n)}
\def \caluqgl {\calU_q(\gergl_n)}
\def \geruegl {\gerU_{\,\varepsilon}(\gergl_n)}
\def \caluegl {\calU_\varepsilon(\gergl_n)}

\def \caluunogl {\calU_1\!(\gergl_n)}
\def \uqsl {U_q(\gersl_n)}
\def \geruqsl {\gerU_q(\gersl_n)}
\def \caluqsl {\calU_q(\gersl_n)}
\def \geruesl {\gerU_{\,\varepsilon}(\gersl_n)}

\def \fqm {F_q[M_n]}
\def \gerfqm {\gerF_q[M_n]}
\def \calfqm {\calF_q[M_n]}
\def \gerfunom {\gerF_1[M_n]}

\def \gerfem {\gerF_\varepsilon[M_n]}
\def \calfem {\calF_\varepsilon[M_n]}
\def \fqg {F_q[G]}
\def \gerfqg {\gerF_q[G]}
\def \calfqg {\calF_q[G]}
\def \gerfeg {\gerF_\varepsilon[G]}
\def \calfeg {\calF_\varepsilon[G]}

\def \fqgl {F_q[{GL}_n]}
\def \gerfqgl {\gerF_q[{GL}_n]}
\def \calfqgl {\calF_q[{GL}_n]}
\def \gerfegl {\gerF_\varepsilon[{GL}_n]}
\def \calfegl {\calF_\varepsilon[{GL}_n]}
\def \gerfunogl {\gerF_1[{GL}_n]}
\def \calfunogl {\calF_1[{GL}_n]}
\def \fqsl {F_q[{SL}_n]}
\def \gerfqsl {\gerF_q[{SL}_n]}
\def \calfqsl {\calF_q[{SL}_n]}
\def \gerfesl {\gerF_\varepsilon[{SL}_n]}
\def \calfesl {\calF_\varepsilon[{SL}_n]}
\def \gerfunosl {\gerF_1[{SL}_n]}
\def \calfunosl {\calF_1[{SL}_n]}

\def \uqgs {U_q(\gerg^*)}
\def \geruqgs {\gerU_q(\gerg^*)}
\def \caluqgs {\calU_q(\gerg^*)}

\def \caluegs {\calU_\varepsilon(\gerg^*)}
\def \uqgls {U_q\big({\gergl_n}^{\!\!*}\big)}
\def \geruqgls {\gerU_q\big({\gergl_n}^{\!\!*}\big)}
\def \caluqgls {\calU_q\big({\gergl_n}^{\!\!*}\big)}
\def \geruegls {\gerU_{\,\varepsilon}\big({\gergl_n}^{\!\!*}\big)}

\def \geruunogls {\gerU_1\!\big({\gergl_n}^{\!\!*}\big)}

\def \uqsls {U_q\big({\gersl_n}^{\!\!*}\big)}
\def \geruqsls {\gerU_q\big({\gersl_n}^{\!\!*}\big)}
\def \caluqsls {\calU_q\big({\gersl_n}^{\!\!*}\big)}
\def \geruesls {\gerU_{\,\varepsilon}\big({\gersl_n}^{\!\!*}\big)}

\def \geruunosls {\gerU_1\!\big({\gersl_n}^{\!\!*}\big)}

\def \gerFr {\frak{Fr}}
\def \calFr {\Cal{F}r}

\def \m {\text{\rm m}}
\def \e {\text{\rm e}}
\def \f {\text{\rm f}}
\def \g {\text{\rm g}}

\def \l {\text{\rm l}}
\def \t {\hbox{\bf t}}
\def \a {\hbox{\bf a}}
\def \b {\hbox{\bf b}}
\def \c {\hbox{\bf c}}
\def \d {\hbox{\bf d}}
\def \ebar {\overline{E}}
\def \fbar {\overline{F}}

\def \N {\Bbb N}
\def \Z {\Bbb Z}
\def \C {\Bbb C}
\def \Q {\Bbb Q}

\def \Qq {\Q (q)}
\def \Zqqm {\Z \!\left[ q, q^{-1}\right]}
\def \Qqqm {\Q \!\left[ q, q^{-1}\right]}
\def \Zeps {\Z_\varepsilon}
\def \Qeps {\Q_\varepsilon}

\document


\topmatter

{\ } 

\vskip-33pt  

 \hfill   {\smallrm {\smallsl Journal of Algebra\/}  {\smallbf 315}  (2007), 761--800}   
\hskip19pt   {\ }  

\vskip51pt  

\title
  $ \boldkey{F}_\boldkey{q} \boldkey{[}
{\boldkey{M}}_{\hskip0,3pt\boldkey{n}}
\boldkey{]} \, $,  $ \boldkey{F}_\boldkey{q} \boldkey{[}
{\boldkey{G}\boldkey{L}}_{\hskip0,7pt\boldkey{n}} \boldkey{]} $
AND  $ \boldkey{F}_\boldkey{q} \boldkey{[}
{\boldkey{S}\boldkey{L}}_{\hskip0,7pt\boldkey{n}}
\boldkey{]} $  \\
  AS QUANTIZED HYPERALGEBRAS
\endtitle

\author
       Fabio Gavarini${}^\dagger$,  \ Zoran Raki\'{c}$^{\,\ddagger}$
\endauthor

\leftheadtext{ Fabio Gavarini{\,}, \ \  Zoran Raki\'{c} }
\rightheadtext{ $ F_q[M_{\,n}] \, $,  $ F_q[{GL\,}_n] $ 
and  $ F_q[{SL\,}_n] $  as quantized hyperalgebras \quad }

\affil
  $ {}^\dagger $Universit\`a di Roma ``Tor Vergata'' ---
Dipartimento di Matematica  \\
  Via della Ricerca Scientifica 1, I-00133 Roma --- ITALY  \\
                             \\
  \hbox{ $ {}^\ddagger $Univerzitet u Beogradu ---
Matemati\v{c}ki fakultet }  \\
  \hbox{ Studentski trg 16, 11000 Beograd --- SERBIA }  \\
\endaffil

\address\hskip-\parindent
  Fabio Gavarini  \newline
  \indent   Universit\`a degli Studi di Roma ``Tor Vergata''  
---   Dipartimento di Matematica  \newline
  \indent   Via della Ricerca Scientifica 1, I-00133 Roma, ITALY  
---   gavarini\@{}mat.uniroma2.it  \newline
  {}  \newline
  Zoran Raki\'{c}  \newline
  \indent   Univerzitet u Beogradu  
---   Matemati\v{c}ki fakultet  \newline
  \indent   Studentski trg 16, 11000 Beograd, SERBIA  
---   zrakic\@{}matf.bg.ac.yu
\endaddress

\abstract
   Within the quantum function algebra  $ F_q[{GL}_n] $,  we study
the subset  $ \calfqgl $   --- introduced in [Ga1] ---   of all
elements of  $ \fqgl $  which are  $ \Zqqm $--valued  when paired
with  $ \caluqgl \, $,  the unrestricted  $ \Zqqm $--integral 
form of  $ \uqgl $  introduced by De Concini, Kac and Procesi.  In
particular we obtain a presentation of it by generators and relations,
and a PBW-like theorem.  Moreover, we give a direct proof that 
$ \calfqgl $  is a Hopf subalgebra of  $ \fqgl $,  and that  $ \,
\calfqgl\Big|_{q=1} \!\!\! \cong U_\Z({\gergl_n}^{\!*}) \, $. 
We describe explicitly its specializations at roots of 1, say 
$ \varepsilon $,  and the associated quantum Frobenius (epi)morphism
from  $ \calfegl \, $  to  $ \, \calfunogl \cong U_\Z({\gergl_n}^{\!*})
\, $,  also introduced in [Ga1].  The same analysis is done for 
$ \, \calfqsl \, $  and (as key step) for  $ \, \calfqm \, $.
\endabstract

\endtopmatter

\footnote""{Keywords: \ {\sl Hopf Algebras, Quantum Groups}.}

\footnote""{2000 {\it Mathematics Subject Classification:} \
Primary 16W30, 17B37; Secondary 81R50.}

\footnote""{Partially supported by a cooperation agreement between
the Department of Mathematics of the University of Rome ``Tor
Vergata'' and the Faculty of Mathematics of the University of
Belgrade (2000--2003).  The first author was also supported by
the Italian Ministry of Teaching, University and Research (MIUR),
within the project COFIN 2003 "Group actions: algebraic and geometric
aspects".  The second author was also supported by the Serbian Ministry
of Science and Environmental Protection, project No.~144032D.}  

\vskip9pt

\centerline {\bf Introduction }

\vskip13pt

   Let  $ G $  be a semisimple, connected, simply connected affine
algebraic group over  $ \Q \, $, and  $ \gerg $  its tangent Lie
algebra.  Let  $ \uqg $  be the Drinfeld-Jimbo quantum group over 
$ \gerg \, $  defined (after Jimbo) as a Hopf algebra over the field 
$ \Q(q) \, $,  \, where  $ q $  is an indeterminate.  After Lusztig,
one has an integral form over  $ \Zqqm $,  say  $ \geruqg \, $,  which
for  $ \, q \rightarrow 1 \, $  specializes to  $ U_\Z(\gerg) \, $, 
the integral  $ \Z $--form  of  $ U(\gerg) $  defined by Kostant
(see [CP], \S 9.3, and references therein, and [DL], \S\S 2--3). 
As  $ U_\Z(\gerg) $  is usually called ``hyperalgebra'', we call 
$ \geruqg $  {\sl ``quantum''  ({\rm or}  quantized) hyperalgebra''.} 
In particular, as  $ U_\Z(\gerg) $  is generated by divided powers
(in the simple root vectors) and binomial coefficients (in the simple
coroots) so  $ \geruqg $  is generated by quantum analogues of divided
powers and of binomial coefficients.  Moreover, if  $ \varepsilon $  is
a root of 1 with odd order  $ \ell \, $,  if  $ \Zeps $  is the formal
extension of  $ \Z $  by  $ \epsilon $  (see \S 1.4) and  $ \gerueg $ 
is the corresponding specialization of  $ \geruqg \, $,  there is a
Hopf algebra epimorphism  $ \; {\frak{Fr}}_\gerg^{\,\Z} \, \colon \,
\gerueg \longtwoheadrightarrow \Zeps \otimes_\Z U_\Z(\gerg) \, $, 
\, described as an  ``$ \ell $-th  root operation'' on generators. 
This is a quantum analogue of the Fro\-benius morphism in positive
characteristic, so is called  
   \hbox{{\sl quantum Frobenius morphism for\/} 
$ \gerg \, $.}   
                                                   \par
   In a Hopf-dual setting, one constructs ([DL], \S\S 4--6) a Hopf
algebra  $ \fqg $  of matrix coefficients of  $ \uqg \, $,  and a 
$ \Zqqm $--form  $ \gerfqg $  of it which specializes to  $ F_\Z[G] $ 
(the algebra of regular functions over  $ G_\Z \, $,  the algebraic
scheme of  $ \Z $--points  of  $ G \, $)  as a Hopf algebra, for  $ \,
q \rightarrow 1 \, $.  In particular,  $ \gerfqg $  is nothing but the
set of ``functions'' in  $ \fqg $  which take values in  $ \Zqqm $  when
``evaluated'' on  $ \geruqg \, $:  in a word, the  $ \Zqqm $--valued 
functions on  $ \geruqg \, $.  When specializing at roots of 1 (with
notation as above) there is a Hopf algebra monomorphism  $ \;
{\gerFr}_G^{\,\Z} \, \colon \,  F_\Z[G] \, \llonghookrightarrow
\, \gerfeg \; $  dual to the above epimorphism and described, roughly,
as an  ``$ \ell $-th  power operation'' on generators.  This also is
a quantum analogue of the classical Frobenius morphism, which is
therefore called the  {\it quantum Frobenius morphism for}  $ G \, $.   
                                                 \par   
   The quantization  $ \geruqg $  of  $ U_\Z(\gerg) $  endows the latter
with a co-Poisson (Hopf) algebra structure which makes  $ \gerg $  into
a Lie bialgebra; similarly,  $ \gerfqg $  endows  $ F_\Z[G] $  with a
Poisson (Hopf) algebra structure which makes  $ G $  into a Poisson group. 
The Lie bialgebra structure on  $ \gerg $  is exactly the one induced by
the Poisson structure on  $ G \, $.  Then one can consider the dual Lie
bialgebra  $ \gerg^* $,  and dual Poisson groups  $ G^\star $  having 
$ \gerg^* $  as tangent Lie bialgebra.   
                                                 \par   
   Lusztig's  $ \Zqqm $--integral  forms  $ \geruqg $  and  $ \gerfqg $ 
are said to be  {\sl restricted}.  On the other hand, another 
$ \Zqqm $--integral  form of  $ \uqg \, $,  \, say  $ \caluqg $, 
has been introduced by De Concini and Procesi (cf.~[CP], \S 9.2,
and [DP], \S 12.1 --- the original construction is over  $ \Bbb{C}
\big[q,q^{-1}\big] $,  but it works the same over  $ \Zqqm $  too),
called  {\sl unrestricted}.  It is generated by suitably rescaled
quantum root vectors and by toral quantum analogues of simple
root vectors, and for  $ \, q \rightarrow 1 \, $  specializes to 
$ F_\Z[G^*] $  (notation as before).  Moreover, at roots of 1 there
is a Hopf algebra monomorphism  $ \; \calFr_\gerg^{\,\Z} \, \colon
\, F_\Z[G^*] = \calU_1(\gerg) \, \longhookrightarrow \, \calueg \, $,  \; defined on
generators as an  ``$ \ell $-th  power operation''.  This is a quantum
analogue of the classical Frobenius morphisms (for  $ G^* $),  strictly
parallel to  $ \calFr_\gerg^{\,\Z} $  above, so is called the  {\it
quantum Frobenius morphism for\/}  $ G^* \, $.  In the Hopf-dual
setting, one can consider   --- [Ga1], \S\S 4, 7 ---   as ``dual''
of  $ \caluqg $  the subset  $ \calfqg $  of ``functions'' in  $ \fqg $ 
which take values in  $ \Zqqm $  when ``evaluated'' on  $ \caluqg \, $; 
\, this subset is a Hopf subalgebra, which for  $ \, q \rightarrow 1
\, $  specializes to  $ U_\Z(\gerg^*) \, $.  When specializing at
roots of 1 there is a Hopf epimorphism  $ \; \calFr_G^{\,\Zeps} \,
\colon \, \calfeg \, \llongtwoheadrightarrow \, \Zeps \otimes_\Z
U_\Z(\gerg^*) \, $,  \; dual to the previous monomorphism; this
again is a quantum analogue of the classical Frobenius morphism
(for  $ \gerg^* \, $),  hence it is called the  {\it quantum
Frobenius morphism for\/}  $ \gerg^* \, $.   
                                                   \par
   In this paper we provide an explicit description of  $ \, \calfqg
\, $,  its specializations at roots of 1 and its quantum Frobenius
morphisms  $ \calFr_G^{\,\Zeps} $  for  $ \, G = SL_n \, $.  The
whole construction makes sense for  $ \, G = GL_n \, $  and  $ \,
G = M_n \, $  (the Poisson algebraic monoid of square matrices of
size  $ n \, $)  too, and we find similar results for them.  In fact,
we first approach the case of  $ \calfqm \, $,  for which the strongest
results are found; then from these we get those for  $ \calfqgl $ 
and  $ \calfqsl \, $.  
                                                   \par
   Our starting point is the well-known description of  $ \gerfqm $ 
by generators and relations, as a  $ \Zqqm $--algebra  generated
by the entries of a quantum  $ (n \times n) $--matrix  (see [APW],
Appendix).  In particular, this is an algebra of skew-commutative
polynomials, much like  $ \caluqg $  is just an algebra of
skew-commutative polynomials (which are Laurent in some variables). 
Dually, this leads us to expect that, like  $ \geruqgl \, $,  \,
also  $ \calfqm $  be generated by quantum divided powers and
quantum binomial coefficients: also, we expect that a suitable
PBW-like theorem holds for  $ \calfqm \, $,  like for  $ \caluqgl \, $.  
 Similarly, as  $ \; {\gerFr}_{M_n}^{\,\Z} \,
\colon \, F_\Z[M_n] \longhookrightarrow \gerfem \; $  is
defined on generators as an  ``$ \ell $-th  power operation'',
dually we expect   
   \hbox{that  $ \; \calFr_{M_n}^{\,\Zeps} \, \colon
\, \calfem \longtwoheadrightarrow \Zeps \otimes_\Z \!
U_\Z\big(\gergl_{\,n}^{\,*}\big) \; $}   
be given by an  ``$ \ell $--th  root operation'', much like 
$ \; \gerFr_{\gergl_n}^{\,\Zeps} \, \colon \, \geruegl \,
\longtwoheadrightarrow \, \Zeps \otimes_\Z U_\Z(\gergl_{\,n}) \, $.
%
                                         \par   
   In fact, all these conjectural
 expectations turn out to be true. 
   From this we get similar (yet slightly weaker) results for 
$ \calfqgl $  and  $ \calfqsl \, $.  On the way, we also improve
the (already known) above mentioned results about specializations
and quantum
Frobenius epimorphisms.
                                                   \par
   The intermediate step is the quantum group  $ \uqgs \, $,  analogue
for  $ \gerg^* $  of what  $ \uqg $  is for  $ \gerg $  (see [Ga1], \S
6).  In particular, there are integral  $ \Zqqm $--forms  $ \geruqgs $ 
and  $ \caluqgs $  of  $ \uqgs $  for which PBW theorems and
presentations hold.  Moreover, a Hopf algebra embedding  $ \, \calfqm
\longhookrightarrow \geruqgls \, $  exists, via which we ``pull
back'' a PBW-like basis and a presentation from  $ \geruqgls $  to 
$ \calfqm \, $.  These arguments work,  {\it mutatis mutandis},  for 
$ GL_n $  and  $ SL_n $  as well.  As aside results, we provide (in
\S 3)  explicit descriptions of these embeddings, and related results
which turn useful in studying specializations at roots of 1.   
                                                   \par  
   The present work bases upon the analysis of the case  $ \, n = 2
\, $,  \, which is treated in [GR].   
 \vskip3pt
  {\it  $ \underline{\hbox{\it Warning}} \, $:}  {\sl at the end of
the paper, we report a short list of the main symbols we use.}   

\vskip13pt

   \centerline{\sc acknowledgments}   
  The authors thank the referee for his/her careful analysis of the
work, and the number of valuable suggestions he/she offered to improve
the paper.   

\vskip11pt

   \centerline{\sc dedicatory}   
  {\smallsl This work grew out of a cooperation supported by an
official agreement between the Department of Mathematics of the
University of Rome ``Tor Vergata'' and the Faculty of Mathematics
of the University of Belgrade in the period 2000--2003.  This
agreement was the outcome of a common wish of peaceful, fruitful
partnership, as an answer to the military aggression of NATO
countries to the Federal Republic of Yugoslavia, which started
in spring of 1999.   
                                      \par   
   This paper is dedicated to the memory of all victims of that war. }   

\vskip1,3truecm

\centerline {\bf \S\; 1 \ Geometrical background and  $ q $--numbers }

\vskip13pt

   {\bf 1.1 Poisson structures on linear groups.}  Let  $ \, \gerg
:= \gergl_n({\Bbb Q}) \, $,  \, with its natural structure of Lie
algebra; it has basis given by the elementary matrices  $ \, \m_{i,j}
:= {\big( \delta_{\ell,i} \, \delta_{j,k} \big)}_{\ell, = 1, \dots,
n}^{k = 1, \dots, n} \, $  ($ \, \forall \; i, j = 1, \dots, n \, $). 
Define  $ \, e_i := \m_{i,i+1} \, $,  $ \, g_j := \m_{j,j} \, $, 
$ \, f_i := \m_{i+1,i} \, $,  \, and also  $ \, h_i := g_i -
g_{i+1} \, $  ($ \, i = 1, \dots, n-1 \, $,  $ \, j = 1, \dots,
n \, $):  then  $ \, \big\{\, e_i \, , \, g_k \, , \, f_i \,\big|\,
i = 1, \dots, n-1, \, k = 1, \dots, n \,\big\} \, $  is a set of Lie
algebra generators of  $ \gerg \, $.  Moreover, a Lie cobracket is
defined on  $ \gerg $  by  $ \, \delta(e_i) = h_i \otimes e_i - e_i
\otimes h_i \, $,  $ \, \delta(g_k) = 0 \, $,  $ \, \delta(f_i) = h_i
\otimes f_i - f_i \otimes h_i \, $  (for all  $ i $  and  $ k \, $), 
which makes  $ \gerg $  itself into a  {\sl Lie bialgebra\/}  ([CP],
\S 1.3.8).  Then  $ \, U(\gerg) \, $  is naturally a co-Poisson Hopf
algebra, whose co-Poisson bracket is the extension of the Lie cobracket
of  $ \gerg $  (the Hopf structure being standard).  Finally, Kostant's 
$ \Z $--integral  form of  $ U(\gerg) $   --- also called  {\sl
hyperalgebra}  ---   is the unital  $ \Z $--subalgebra  $ U_\Z(\gerg) $ 
of  $ U(\gerg) $  generated by the divided powers  $ \, f_i^{(m)} \, $, 
$ \, e_i^{(m)} \, $  and the binomial coefficients  $ \, \Big(\! {g_k
\atop m} \!\Big) \, $  (for all  $ i $,  $ k $,  and all  $ \, m \in
\N \, $),  where hereafter we use notation  $ \, x^{(m)} := x^m \big/
m! \, $  and  $ \, \Big(\! {t \atop m} \!\Big) := {{\, t (t-1) \cdots
(t-m+1) \,} \over {\, m! \,}} \, $.  This again is a co-Poisson Hopf
algebra (over  $ \Z $);  it is free as a  $ \Z $--module,  with PBW-like 
$ \Z $--basis  the set of ordered monomials  $ \; \Big\{\! \prod_{i<j}
e_{i,j}^{(\eta_{i,j})} \prod_{k=1}^n \Big(\! {g_k \atop \gamma_k} \!\Big)
\prod_{i>j} f_{i,j}^{(\varphi_{i,j})} \,\Big|\, \eta_{i,j}, \gamma_k,
\varphi_{i,j} \in \N \,\Big\} \, $  (w.r.t.~any total order of the pairs 
$ (i,j) $  with  $ \, i \not= j \, $), where  $ \, e_{i,j} := \m_{i,j}
\, $  and  $ \, f_{j,i} := \m_{j,i} \, $  for all  $ \, i < j \, $; 
see e.g.~[Hu], Ch.~VII.  
                                            \par
   As  $ \gergl_n(\Q) $  is a Lie bialgebra, by general theory  $ \,
G := {GL}_n(\Q) \, $  is then a  {\sl Poisson\/}  group, its Poisson
structure being the unique one which induces on  $ \, \gerg := \gergl_n(\Q)
\, $  the Lie bialgebra structure mentioned above.  Explicitly,  the
algebra  $ F[G] $  of regular functions on  $ G $  is the unital associative
commutative  $ \Q $--algebra  with generators  $ \, \bar{t}_{i,j} \, $  ($ \,
i, j = 1, \dots, n \, $)  and  $ D^{-1} $,  where  $ \, D := \text{\sl det}
\,\big(\bar{t}_{i,j}\big)_{i,j=1,\dots,n} \, $  is the determinant.  The
group structure on  $ G $  yields on  $ F[G] $  the natural Hopf structure
given by matrix product; in particular  $ D^{-1} $  is group-like.  The
Poisson structure   --- usually referred to as ``standard'' ---   is given
by (see, e.g., [BGY], \S\S 1.4--5)
 \vskip-9pt   
  $$  \matrix  
   \big\{\bar{t}_{i,j}, \bar{t}_{i,k}\big\} = \, \bar{t}_{i,j} \,
\bar{t}_{i,k}  \hskip13pt  \forall \; j < k \, ,  &  \big\{\bar{t}_{i,j},
\bar{t}_{\ell,j}\big\} = \, \bar{t}_{i,j} \, \bar{t}_{\ell,j}  \hskip13pt 
\forall \; i < \ell  \cr   
   \hskip5pt \phantom{\bigg|}  \big\{\bar{t}_{i,j}, \bar{t}_{\ell,k}\big\}
= \, 0  \hskip13pt  \forall \; i < \ell, k < j \, ,  &  \hskip31pt 
\big\{\bar{t}_{i,j}, \bar{t}_{\ell,k}\big\} = \, 2 \; \bar{t}_{i,k}
\, \bar{t}_{\ell,j}  \hskip13pt  \forall \; i < \ell, j< k   
   \endmatrix  $$   
 \vskip-11pt   
  $$  \big\{ D^{-1}, \bar{t}_{i,j} \big\} \, = \, - \, \big\{ D ,
\bar{t}_{i,j} \big\}  \hskip13pt  \forall \; i, j= 1, \dots, n \, .  $$   
   \indent   We shall also consider the Poisson group-scheme  $ G_\Z $ 
associated to  $ {GL}_n $,  \, for which a like analysis applies: in
particular, its function algebra  $ F[G_\Z] $  is a Poisson Hopf 
$ \Z $--algebra  with the same presentation as  $ F[G] $  but over
the ring  $ \Z \, $.   
                                           \par   
   Similar constructions hold with  $ \gersl_n $  replacing  $ \gergl_n $ 
and with  $ {SL}_n $  instead of   $ {GL}_n \, $:  \, one simply replaces
the  $ g_k $'s  with the  $ h_i $'s  ($ \, i=1, \dots, n-1 \, $).  In
particular,  $ \gersl_n(\Z) $  is a Lie  $ \Z $--sub-bialgebra  of 
$ \gergl_n(\Z) \, $,  and  $ U_\Z\big(\gersl_n\big) $  is a co-Poisson
Hopf subalgebra (over  $ \Z $)  of  $ \, U_\Z\big(\gergl_n\big) \, $, 
with a PBW-like  $ \Z $--basis  as above but with the  $ h_i $'s 
instead of the  $ g_k $'s;  moreover,  $ F \big[ ({SL}_n)_\Z \big] $ 
is the quotient Poisson-Hopf algebra of  $ F \big[ ({GL}_n)_\Z \big] $ 
modulo the principal ideal  $ \, \big( D - 1 \big) \, $.   
                                                \par
   Finally, the subalgebra of  $ F\big[({GL}_n)_{\Z}\big] $  generated
by the  $ \bar{t}_{i,j} $'s  alone clearly is a Poisson subbialgebra of 
$ F\big[({GL}_n)_{\Z}\big] \, $:  indeed, it is the algebra  $ F\big[
(M_n)_{\Z}\big] $  of regular functions of the  $ \Z $-scheme 
associated to the Poisson algebraic monoid  $ \, M_n \, $  of
all  $ (n \times n) $--matrices.   

\vskip9pt

   {\bf 1.2 Dual Lie bialgebras and dual Poisson groups.}  Let  $ \,
\gerg := \gergl_n\!(\Q) \, $  be a Lie bialgebra structure as in \S
1.1; then its dual  $ \, \gerg^* \, $  is a Lie algebra too.  Let 
$ \, \big\{ e_{i,j}^{\,*} \, , g_k^{\,*} \, , f_{j,i}^{\,*} \;\big|\;
\forall \, k \, , \, \forall \, i \! < \! j \, \big\} \, $  be the dual
basis to the basis of elementary matrices for  $ \gerg \, $.  Set  $ \,
\e_{\,i,j} := e_{i,j}^{\,*} \big/ 2 \, $,  $ \, \g_k := g_k^{\,*} \, $, 
$ \, \f_{j,i} := f_{j,i}^{\,*} \big/ 2 \, $  (for all  $ k $  and all 
$ \, i < j \, $),  and  $ \, \e_{\,i} := \e_{\,i,i+1} \, $,  $ \,
\f_{\,i} := \f_{\,i,i+1} \, $  for all  $ \, 1 \leq i < n \, $. 
Then  $ \gerg^* $  is the Lie algebra with generators  $ \,
\e_{\,1} \, , \dots , \e_{\,n-1} \, , \g_1 \, , \dots , \g_n
\, , \f_{\,1} \, , \dots $,  $ \f_{\,n-1} \, $  and relations   
 \vskip-14pt   
  $$  \displaylines{ 
   \big[ \g_i \, , \, \f_j \big] \, = \, \big( \delta_{i,j} - 
\delta_{i-1,j} \big) \, \f_j \; ,  \hskip9pt 
   \big[ \f_{\,k} \, , \, \e_\ell \big] \, = \, 0 \; ,  \hskip9pt 
   \big[ \g_i \, , \, \e_j \big] \, = \, \big( \delta_{i,j} -
\delta_{i-1,j} \big) \, \e_j  \hskip17pt \forall\; k, \ell , i, j  \cr 
   \big[ \g_k \, , \, \g_\ell \big] \, = \, 0 \; ,  \hskip11pt 
   \big[ \f_{\,i} \, , \, \f_j \big] \, = \, 0 \; ,  \hskip11pt 
   \big[ \e_{\,i} \, , \, \e_j \big] \, = \, 0  \hskip25pt 
\forall\; k, \ell \, ,  \hskip11pt  \forall\; i, j : 
\vert i-j \vert > 1  \cr   
   \big[ \f_{\,i} \, , \big[ \f_{\,i} \, , \, \f_j \big] \big]
= \, 0 \; ,  \hskip23pt  
   \big[ \e_{\,i} \, , \big[ \e_{\,i} \, , \, \e_j \big] \big]
= \, 0  \hskip31pt  \forall\; i, j : \vert i-j \vert = 1  \cr }  $$   
%
%
%
 \vskip-1pt   
\noindent   
 (just compute!  See also [Ga2], \S 1).  The Lie cobracket is given
(for  $ \, x \wedge y := x \otimes y - y \otimes x \, $)  by  
 \vskip-14pt   
  $$  \displaylines{ 
   \delta\big(\f_{\,i}\big) \, = \, \f_{\,i} \wedge
\big(\g_i - \g_{i+1}\big) + 2 \cdot \big( {\textstyle \sum_{j=1}^{i-1}}
\, \f_{\,i+1,j} \wedge \e_{j,i} + {\textstyle \sum_{j=i+2}^{n+1}}
\, \e_{\,i+1,j} \wedge \f_{j,i} \big) \; ,  
      \quad   1 \leq i \leq n\!-\!1 \, \phantom{.}  \cr   
   \qquad \qquad   \delta\big(\g_k\big) \, = \; 4 \cdot {\textstyle
\sum_{\ell=1}^{k-1}} \, \f_{\,k,\ell} \wedge \e_{\,\ell,k}
\, + \, 4 \cdot {\textstyle \sum_{\ell = k+1}^{n+1}} \,
\e_{\,k,\ell} \wedge \f_{\,\ell,k} \;\; , 
      \qquad \qquad \qquad   1 \leq k \leq n  \cr   
   \delta\big(\e_{\,i}\big) \, = \, \big( \g_i - \g_{i+1} \big)
\wedge \e_{\,i} + 2 \cdot \big( {\textstyle \sum_{\ell=1}^{i-1}}
\, \f_{\,i,\ell} \wedge \e_{\,\ell,i+1} \, + \, {\textstyle
\sum_{\ell=i+2}^{n+1}} \, \e_{\,i,\ell} \wedge \f_{\,\ell,i+1} \big) \; , 
      \quad   1 \leq i \leq n\!-\!1 \, .  \cr }  $$   
 \vskip-3pt   
%
%
   \indent   All these formul{\ae}  also provide a presentation of 
$ U(\gerg^*) $  as a co-Poisson Hopf algebra.  Finally, we define the
Kostant's  $ \Z $--integral  form, or hyperalgebra,  $ U_\Z(\gerg^*) $ 
of  $ U(\gerg^*) $  as the unital  $ \Z $--subalgebra  generated by
the divided powers  $ \, \f_i^{\,(m)} $,  $ \, \e_i^{(m)} $  and
binomial coefficients  $ \Big(\! {\g_k \atop m} \!\Big) $  (for all 
$ \, m \in \N \, $  and all  $ i $,  $ k \, $).  This again is a
co-Poisson Hopf  $ \Z $--algebra,  
%
%
 with
PBW-like  $ \Z $--basis  the set   
     \hbox{$ \; \Big\{\, \prod_{i<j} \e_{i,
j}^{(\eta_{i,j})} \prod_{k=1}^n \Big(\! {\g_k \atop \gamma_k} \!\Big)
\prod_{i>j} \f_{i,j}^{\,(\varphi_{i,j})} \,\Big|\, \eta_{i,j},
\gamma_k, \varphi_{i,j} \in \N \Big\} \; $}
  of ordered monomials
(w.r.t.~any total order as before).  Alternatively, one can take
the elements  $ \, \l_k := \g_1 \! + \cdots + \g_k \, $  ($ \, k
= 1, \dots, n \, $)  instead of the  $ \g_k $'s  in the definition
of  $ U_\Z(\gerg^*) $;  the presentation changes accordingly,
and  $ \, \Big\{ \prod_{i<j} \e_{i,j}^{(\eta_{i,j})} \prod_{k=1}^n
\!\! \Big(\! {\l_k \atop \lambda_k} \!\Big) \prod_{i>j} \f_{i,j}^{\,
(\varphi_{i,j})} \,\Big| \eta_{i,j}, \lambda_k, \varphi_{i,j} \in
\N \Big\} \; $   
   \hbox{is another PBW-like basis.}   
 \eject   
%
%
   A like description holds for  $ \, {\gersl_n(\Q)}^* \, $.  In
fact, 
%
%
  $ \, {\gersl_n
(\Q)}^* \! = {\gergl_n(\Q)}^* \!\Big/
 \Q \cdot \l_n \, $,  \,  
%
%
dually to the embedding  $ \, \gersl_n(\Q) \longhookrightarrow
\gergl_n(\Q) \, $;  \, thus one simply has to add to the presentation
of  $ {\gergl_n(\Q)}^* $  the additional relation  $ \, \g_1 \! + \cdots
+ \g_n = 0 \, $,  or  $ \, \l_n = 0 \, $.  Then  $ U_\Z\big(\gersl_n^{\,*}
\big) $  is the  $ \Z $--subalgebra  of  $ U \big(\gersl_n(\Q)^*\big) $ 
generated by divided powers and binomial coefficients as above, but
taking care of the additional relation); also,  $ \, U_\Z \big(
\gersl_n^{\,*} \big) \, $  is a quotient co-Poisson Hopf  $ \Z $--algebra 
of  $ \, U_\Z\big(\gergl_n^{\,*} \big) \, $.   
%
%
 Again,   
%
%
 admits 
the set  $ \, \Big\{ \prod_{i<j} \e_{i,
j}^{(\eta_{i,j})} \prod_{k=1}^{n-1} \Big(\! {\l_h \atop \lambda_h}
\!\Big) \prod_{i>j} \f_{i,j}^{\,(\varphi_{i,j})} \,\Big|\, \eta_{i,j},
\lambda_k, \varphi_{i,j} \in \N \,\Big\} \, $
 as  $ \Z $--basis.  
 Finally, note that  $ \, \Big\{\, \e_{i,j} := e_{i,j}^{\,*} \big/ 2 \, ,
\, \l_k \, , \, \f_{j,i} := f_{j,i}^{\,*} \big/ 2 \;\,\Big|\;\, 1 \leq k
\leq n-1 , \, 1 \leq i < j \leq n \,\Big\} \, $  is the basis of 
$ \gersl_n^{\,*} $  dual to the basis  $ \, \big\{\, e_{i,j} \, ,
h_k \, , f_{j,i} \;\big|\; 1 \leq k \leq n-1 , \, 1 \leq i < j
\leq n \,\big\} \, $  of  $ \gersl_n \, $.  
                                            \par   
   If  $ \, \gerg = \gergl_n \, $,  \, a simply connected algebraic
Poisson group with tangent Lie bialgebra  $ \gerg^* $  is the subgroup 
$ {}_s{G}^* $  of  $ \, G \times G \, $  made of all pairs  $ \, \big(
L, U \big) \in G \times G \, $  such that  $ L $  is lower triangular, 
$ U $  is upper triangular, and their diagonals are inverse to each other
(i.e.~their product is the identity matrix).  This is a Poisson subgroup
for the natural Poisson structure of  $ G \times G \, $.  Its centre
is  $ \, Z := \big\{ \big(z I , z^{-1} I \big) \,\big|\, z \in \Q
\setminus \{0\} \big\} \, $,  \, hence the associated adjoint group
is  $ \, {}_a{G}^* := {}_s{G}^* \Big/ Z \, $.  The same construction
defines Poisson group-schemes  $ {}_s{G}^*_\Z \, $  and 
$ {}_a{G}^*_\Z \; $.   
                                            \par   
   If  $ \, \gerg = \gersl_n \, $  the construction of dual Poisson
group-schemes  $ {}_s{G}^*_\Z \, $  and  $ {}_a{G}^*_\Z \; $  is
entirely similar, just taking  $ {SL}_n $  instead of  $ {GL}_n $ 
in the previous recipe.   

\vskip11pt

   {\bf 1.3  $ q $--numbers,  $ q $--divided  powers and  $ q $--binomial 
coefficients.}  Let  $ q $  be an indeterminate.  For all  $ \, n, s, k_1,
\dots, k_r \in \N \, $,  let  $ \; {(n)}_q := {{q^n - 1} \over {q - 1}}
\, $,  $ \; {(n)}_q! := \prod_{r=1}^n {(r)}_q \, $,  $ \; {\Big( {n \atop
s} \Big)}_{\!q} := {{{(n)}_q!} \over {{(s)}_q! {(n-s)}_q!} } \, $,  $ \;
{\Big( {n \atop {k_1, \dots, k_r}} \Big)}_{\!q} := {{{(n)}_q!} \over
{{(k_1)}_q! \cdots {(k_r)}_q!} } \;\; (\in \Z[q]) \, $,  \; then  $ \;
{[n]}_q := {{q^n - q^{-n}} \over {q - q^{-1}}} \, $,  $ \; {[n]}_q! :=
\prod_{r=1}^n {[r]}_q \, $,  $ \; {\Big[ {n \atop s} \Big]}_q \! :=
{{{[n]}_q!} \over {{[s]}_q! {[n-s]}_q!} } \;\; (\in \Zqqm) \, $, 
\, and also  $ {\Big( {-n \atop s} \Big)}_{\! q} \! := {(-1)}^s
q^{- n s - {s \choose 2}} {{{(n \! - \! 1 \! + \! s)}_q!} \over
{(s)}_q! } \;\; (\in \Z[q]) \, $.   
                                       \par   
   If  $ A $  is any  $ \Qq $--algebra  we define
$ q $--divided  powers and  $ q $--binomial  coefficients as
follows: for every  $ \, X \in A \, $,  \, each  $ \, n \in \N \, $
and each  $ \, c \in \Z \, $  we set  $ \; X^{(n)} := X^n \Big/ {[n]}_q!
\; $  and  $ \; \Big( {{X \, ; \, c} \atop n} \Big) := \prod_{s=1}^n {{\,
q^{c+1-s} X - 1 \,} \over {\, q^s - 1 \,}} \; $.  Furthermore, if 
$ \, Z \in A \, $  is  {\sl invertible\/}  we define also  $ \; \left[
{{X \, ; \, c} \atop n} \right] := \prod_{s=1}^n {{\, q^{c+1-s} Z^{+1}
- q^{s-1-c} Z^{-1} \,} \over {\, q^{+s} - q^{-s} \,}} \; $  for every 
$ \, n \in \N \, $  and  $ \, c \in \Z \, $.  We shall also consider the
elements  $ \; \left\{ {{X \, ; \, c} \atop {n \, , \, r}} \right\} \,
:= \, \sum_{s=0}^r q^{{s+1} \choose 2} {\Big( {r \atop s} \Big)}_{\!q}
\cdot \Big( {{X \, ; \, c+s} \atop {n-r}} \Big) \; $  for every  $ \,
X \in A \, $,  $ \, n, r \in \N \, $,  $ \, c \in \Z \, $.

                                             \par
   For later use, we remark that the  $ q $--divided  powers in 
$ \, X \in A \, $  satisfy the (obvious) relations   
  $$  X^{(r)} X^{(s)} \, = \left[ {{r+s} \atop s} \right]_q X^{(r+s)}
\; ,  \qquad  X^{(0)} = 1   \eqno (1.1)  $$
   Similarly, the  $ q $--binomial  coefficients in  $ \, X \in A \, $ 
satisfy the relations
  $$  \hbox{ $ \eqalign{  
   {{X \, ; c} \choose t} {{X \, ; c-t} \choose s} =  &
\; {{t+s} \choose t}_{\!\! q} {{X \, ; c} \choose {t+s}} \; ,
\qquad  {{X \, ; c+1} \choose t} - q^t {{X \, ; c} \choose t}
= {{X \, ; c} \choose {t-1}}  \cr
   {{X ; c} \choose m} {{X ; s} \choose n} =  &  \; {{X ; s} \choose n}
{{X ; c} \choose m} \; ,  \quad  {{X \, ; c} \choose t} =
{\textstyle \sum\limits_{p \geq 0}^{p \leq c,t}} \, q^{(c-p)(t-p)}
{c \choose p}_{\!\! q} {{X \, ; 0} \choose {t-p}}  \cr
   {{X \, ; c} \choose 0} = 1 \; ,  \quad  &  {{X \, ; -c} \choose t}
= {\textstyle \sum\limits_{p=0}^t} \, {(-1)}^p q^{- t(c+p) + p(p+1)/2}
{{p+c-1} \choose p}_{\!\! q} {{X \, ; 0} \choose {t-p}}  \cr
   {{X \, ; c+1} \choose t} -  &  {{X \, ; c} \choose t} = \,
q^{c-t+1} \left( 1 + (q-1) {{X \, ; 0} \choose 1} \! \right)
{{X \, ; c} \choose {t-1}}  \cr } $ }   \hskip17pt \hfill (1.2)  $$
%
%
 \eject   

   {\bf 1.4  Roots of 1 and specializations.} \, Let  $ \, \ell \in \N_+
\, $  be odd, set  $ \, \Zeps := Z[q] \Big/ \big(\phi_\ell(q)\big) \, $ 
where  $ \, \phi_\ell(q) \, $  is the  $ \ell $-th  cyclotomic polynomial
in  $ q \, $,  \, and let  $ \, \varepsilon := \overline{q} \, $,  \, a
(formal) primitive  $ \ell $-th  root of 1 in  $ \Zeps \, $.  Similarly,
let  $ \, \Qeps := \Q[q] \Big/ \big(\phi_\ell(q)\big) \, $,  \, the field
of quotients of  $ \Zeps \, $.  If  $ M $  is a module over  $ \Zqqm $ 
or  $ \Qqqm $  we set  $ \; M_\varepsilon := M \Big/ \big( \phi_\ell(q)
\big) \, M \, $,  \, which is a module over  $ \Zeps $  or over  $ \Qeps
\, $,  \, called the  {\sl specialization of  $ M $  at  $ \, q =
\varepsilon \,$}.  In fact, we have  $ \, M_\varepsilon \cong \Zeps
\otimes_{\Zqqm} M \, $  or  $ \, M_\varepsilon \cong \Zeps \otimes_{\Zqqm}
M \, $  as modules over  $ \Zeps $  or  $ \Qeps $  (in either case).  

\vskip11pt

   For later use, we recall now some results about  $ q $-numbers  and
their specializations:  

\vskip11pt

 \proclaim{Lemma 1.5}  Let  $ A $  be any\/ 
$ \Q\big[p,p^{-1}\big] $--algebra,  $ p $ 
being any indeterminate.  
                                        \par   
%
%
   (a) \, If  $ \, x_1, \dots, x_s \in A \, $  are such that  $ \;
x_i \, x_j = p \, x_j \, x_i \, $  for all  $ \, i < j \, $,  then 
for all  $ \, m \in \N \, $   
 \vskip-4pt
  $$  {\big(\, {\textstyle \sum_{k=1}^s} \, x_k \big)}^m \, = \,
{\textstyle \sum\limits_{k_1 + \cdots + k_s = m} {\Big(\! {m \atop
{k_1, \dots, k_s}} \!\Big)}_{\! p^{-1}}} x_1^{k_1} \cdots x_s^{k_s}  $$   
 \vskip-4pt
   (b) \, If  $ \, y_1, \dots, y_s \in A \, $  are such that  $ \;
y_i \, y_j = p^2 \, y_j \, y_i \, $  for all  $ \, i < j \, $,  then 
for all  $ \, m \in \N \, $   
 \vskip-3pt
  $$  {\big(\, {\textstyle \sum_{k=1}^s} \, y_k \big)}^{(m)} \, =
\, {\textstyle \sum\limits_{k_1 + \cdots + k_s = m}} p^{{{k_1 + 1}
\choose 2} + \cdots + {{k_s + 1} \choose 2} - {{m + 1} \choose 2}}
\cdot y_1^{(k_1)} \cdots y_s^{(k_s)}  $$   
\endproclaim   
 \vskip-23pt
\demo{Proof}  {\it (a)\/}  is the well-known, generalized 
$ q $--Leibniz'  identity, and  {\it (b)\/}  follows
trivially.   \hskip15pt \hfill \qed\break   
\enddemo   

\vskip1pt

 \proclaim{Lemma 1.6}  Let  $ A $  be any\/  $ \Qq $--algebra,  and let 
$ \, x, y, z, w \in A \, $  be such that  $ \; x \, w = q^2 \, w \, x
\, $,  $ \; x \, y = q \, y \, x \, $,  $ \; x \, z = q \, z \, x \; $ 
and  $ \; y \, z = z \, y \, $.  Then for all  $ \, m \in \N \, $  and 
$ \, t \in \Z \, $  we have    
  $$  \leqalignno{ 
   \left( {{x + \big( q - q^{-1} \big)^2 \, w \; ; \, t} \atop m} \right)
\, =  &  \, \sum_{r=0}^m \, q^{r(t-m)} \big( q - q^{-1} \big)^r \cdot
w^{(r)} \left\{ {x \, ; \, t} \atop {m \, , \, r} \right\} \, =   &  
\hbox{\it (a-1)}  \cr   
   =  &  \, \sum_{r=0}^m \, q^{r(t-m)} \big( q - q^{-1} \big)^r \cdot
\left\{ {x \, ; \, t - 2 \, r} \atop {m \, , \, r} \right\} \, w^{(r)}  
&   \hbox{\it (a-2)}  \cr   
   \left( {{x + \big( q - q^{-1} \big)^2 \, y \, z \; ; \, t} \atop m}
\right) \, =  &  \, \sum_{r=0}^m \, q^{r(t-m)} \big( q - q^{-1} \big)^r
\, {[r]}_q! \cdot y^{(r)} \left\{ {x \, ; \, t - r} \atop {m \, , \, r}
\right\} \, z^{(r)}   &  \hbox{\it (b)}  \cr }  $$   
                                                \par    
\noindent   
  If instead  $ \; w \, x = q^2 \, x \, w \, $,  $ \; y \, x = q \, x \,
y \, $,  $ \; z \, x = q \, x \, z \; $  and  $ \; z \, y = y \, z \, $, 
\, then  ($ \, \forall \; m \in \N \, , \, t \in \Z \, $)   
  $$  \leqalignno{ 
   \left( {{x + \big( q - q^{-1} \big)^2 \, w \; ; \, t} \atop m} \right)
\, =  &  \, \sum_{r=0}^m \, q^{r(t-m)} \big( q - q^{-1} \big)^r \cdot
\left\{ {x \, ; \, t} \atop {m \, , \, r} \right\} \, w^{(r)} \, =  
&   \hbox{\it (c-1)}  \cr   
   =  &  \, \sum_{r=0}^m \, q^{r(t-m)} \big( q - q^{-1} \big)^r \cdot
w^{(r)} \left\{ {x \, ; \, t - 2 \, r} \atop {m \, , \, r} \right\}   &  
\hbox{\it (c-2)}  \cr   
   \left( {{x + \big( q - q^{-1} \big)^2 \, y \, z \; ; \, t} \atop m}
\right) \, =  &  \, \sum_{r=0}^m \, q^{r(t-m)} \big( q - q^{-1} \big)^r
\, {[r]}_q! \cdot z^{(r)} \left\{ {x \, ; \, t - r} \atop {m \, , \, r}
\right\} \, y^{(r)}   &  \hbox{\it (d)}  \cr }  $$   
\endproclaim   

\demo{Proof}  Claims  {\it (a-1/2)\/}  and  {\it (b)\/}  are proved
in Lemma 6.2 of [GR].  Instead, claims  {\it (c-1/2)\/}  and  {\it
(d)\/}  follow directly from claims  {\it (a-1/2)\/}  and  {\it (b)\/} 
when one applies them to the algebra  $ A $  endowed with the  {\sl
opposite\/}  product.   \hskip225pt \hfill \qed\break    
\enddemo   

\vskip11pt

 \proclaim{Lemma 1.7}  ([GR], Lemma 6.3) Let  $ \varOmega $  be any 
$ \Zqqm $--algebra,  \, let  $ \, \varepsilon \, $  be a (formal)
primitive  $ \ell $-th  root of 1, with  $ \, \ell \in \N_+ \, $, 
\, and let  $ \, \varOmega_\varepsilon := \varOmega \Big/ (q -
\varepsilon) \, \varOmega \, $  be the specialization of 
$ \varOmega $  at  $ \, q = \varepsilon \, $  (a 
$ \Zeps $--algebra).  Then for each  $ \, x, y
\in \varOmega \, $  we have   
 \vskip-15pt
  $$  \varOmega \, \ni {{x \, ; \, 0} \choose \ell} \;\; \Longrightarrow
\;\, \big(x\big|_{q=\varepsilon}\big)^\ell = \, 1 \text{\ \ in \ } \!
\varOmega_\varepsilon \, ,  \qquad \!  \varOmega \ni y^{(\ell)} \;\;
\Longrightarrow \;\, \big( y\big|_{q=\varepsilon} \big)^\ell =
\, 0 \text{\ \ in \ } \! \varOmega_\varepsilon \; .  
\eqno \square  $$   
\endproclaim

 \vskip1,1truecm

\centerline {\bf \S\; 2 \ Quantum groups }

\vskip13pt

  {\bf 2.1 The quantum groups  $ \uqgl $  and  $ \uqsl \, $.} \, Let 
$ \uqgl $  be the Drinfeld-Jimbo quantization of  $ U(\gergl_n) \, $: 
\, following [Ji], \S 2 (with normalizations of [No], \S 1.2), we
define it as the unital  $ \Qq $--algebra  with generators  $ \, F_i
\, $,  $ G_j^{\pm 1} $,  $ E_i \, $  ($ 1 \leq i \leq n\!-\!1 \, ;
\, 1 \leq j \leq j \, $)  and relations   
  $$  \displaylines { 
   \hfill   G_i \, G_i^{-1} \, = \, 1 \, = \, G_i^{-1} \, G_i \, , 
\qquad  G_i^{\pm 1} \, G_j^{\pm 1} \, = \, G_j^{\pm 1} \, G_i^{\pm 1} 
\qquad   \hfill  \forall \; i, j \qquad  \cr
   \hfill   G_i \, F_j \, G_i^{-1} \, = \, q^{\delta_{i,j+1} - \delta_{i,j}}
\, F_j \, ,  \qquad  G_i \, E_j \, G_i^{-1} \, = \, q^{\delta_{i,j} -
\delta_{i,j+1}} \, E_j  \qquad   \hfill  \forall \, i, j \qquad  \cr
   \hfill   E_i \, F_j \, - \, F_j \, E_i \, = \, \delta_{i,j} \, {{\; G_i
\, G_{i+1}^{-1} \, - \, G_i^{-1} \, G_{i+1} \;} \over {\; q - q^{-1} \;}} 
\qquad \qquad   \hfill  \forall \; i, j  \qquad  \cr
   \hfill   \quad  E_i \, E_j \, = \, E_j \, E_i \, ,  \qquad  F_i \, F_j
\, = \, F_j \, F_i   \hfill  \forall \;\; i, j : \vert i - j \vert > 1 \,
\phantom{.} \;  \cr  
   \hfill   E_i^2 \, E_j \, - \, \left( q + q^{-1} \right) \, E_i \, E_j \,
E_i \, + \, E_j \, E_i^2 \, = \, 0   \hfill  \forall \;\; i, j : \vert
i - j \vert = 1 \, \phantom{.} \;  \cr  
   \hfill   F_i^2 \, F_j \, - \, \left( q + q^{-1} \right) \, F_i \, F_j \,
F_i \, + \, F_j \, F_i^2 \, = \, 0   \hfill  \forall \;\; i, j : \vert i -
j \vert = 1 \, . \;  \cr }  $$   
   \indent   Moreover,  $ \uqgl $  has a Hopf algebra structure, given by
  $$  \displaylines { 
   \hfill   \Delta \left( F_i \right) \, = \, F_i \otimes G_i^{-1} \,
G_{i+1} \, + \, 1 \otimes F_i \, ,  \qquad  S \left( F_i \right) \,
= \, - F_i \, G_i \, G_{i+1}^{-1} \, ,  \qquad  \epsilon \left( F_i
\right) \, = \, 0   \hfill \qquad  \forall \;\; i \,  \phantom{.} \;  \cr 
   \hfill   \qquad  \Delta \left( G_i^{\pm 1} \right) \, = \, G_i^{\pm 1}
\otimes G_j^{\pm 1} \, ,  \qquad \qquad \qquad  S \left( G_j^{\pm 1} \right)
\, = \, G_j^{\mp 1} \, ,  \quad \qquad  \epsilon \left( G_j^{\pm 1} \right)
\, = \, 1   \hfill \qquad  \forall \;\; j \, \phantom{.} \;  \cr 
   \hfill   \Delta \left( E_i \right) \, = \, E_i \otimes 1 \, + \, G_i
G_{i+1}^{-1} \otimes E_i \, ,  \qquad  S \left( E_i \right) \, = \, -
G_i^{-1} \, G_{i+1} \, E_i \, ,  \qquad  \epsilon \left( E_i \right)
\, = \, 0   \hfill \qquad  \forall \;\; i \, . \;  \cr }  $$   
   \indent   Consider also the elements  $ \; K_i^{\pm 1} := G_i^{\pm 1}
G_{i+1}^{\mp 1} \, $  (for all  $ \, i = 1, \dots, n-1 \, $).  The
unital  $ \Qq $--subalgebra  of  $ \uqgl $  generated by  $ \; \big\{
F_i \, , K_i^{\pm 1}, E_i \big\}_{i=1,\dots,n-1} \; $  is the well known
Drinfeld-Jimbo's quantum algebra  $ \uqsl \, $.  From the presentation
of  $ \uqgl $  one deduces one of  $ \uqsl $  too, and also sees that
the latter is even a Hopf subalgebra of the former.   

\vskip11pt

  {\bf 2.2 Quantum root vectors and quantum PBW theorem.} \, Now we
introduce quantum analogues of root vectors for both  $ \uqgl $  and 
$ \uqsl $,  following an idea due to Jimbo (as in [Ji], \S 2, up to
changing  $ \, q^{\pm 1} \longleftrightarrow q^{\mp 1} \, $).  Set
 \vskip-11pt
  $$  \matrix   
   E_{i,i+1} := E_i \; ,  &  \quad  E_{i,j} :=
\big[\, E_{i,k} \, , E_{k,j} \,\big]_{q^{-1}} \equiv \,
E_{i,k} \, E_{k,j} - q^{-1} \, E_{k,j} \, E_{i,k}   & 
\qquad  \forall \; i < k < j \; \phantom{.}  \\
   F_{i+1,i} := F_i \; ,  &  \quad  F_{j,i} :=
\big[\, F_{j,k} \, , F_{k,i} \,\big]_{q^{+1}} \equiv \,
F_{j,k} \, F_{k,i} - q^{+1} \, F_{k,i} \, F_{j,k}   & 
\qquad  \forall \; j > k > i \; .  \endmatrix  $$   
 \vskip-3pt
Then all these are quantum root vectors, in  $ \uqgl $  or  $ \uqsl $, 
and the very definition gives a complete set of commutation relations
among them.  For later use, we point out that they can be obtained
also via a general method given by Lusztig. 
                                        \par    
  Namely, let  $ \, \{\alpha_1,\dots,\alpha_{n-1}\} $  be a basis
of simple roots of  $ \gersl_n $  (or of  $ \gergl_n $),  and let 
$ s_1 $,  $ \dots $,  $ s_{n-1} $  be generators of its Weyl group 
$ W \, $.  Let  $ \, w_0 = s_{i_1} s_{i_2} s_{i_3} \cdots s_{i_{N-1}}
s_{i_N} \, $  be a reduced expression for the longest element  $ w_0 $ 
of  $ W \, $,  \, where  $ \, N := {n \choose 2} \, $.  Then  $ \,
\big\{\, \alpha^k := s_{i_1} s_{i_2} \cdots s_{i_{k-1}}(\alpha_{i_k})
\;\big|\; k=1, \dots, N \,\big\} \, $  is the set of all positive roots
of  $ \gersl_n \, $.  Now, Lusztig defines an action of the braid group
associated to  $ W $,  with generators  $ T_1 $,  $ \dots $,  $ T_{n-1}
\, $:  \, we consider it normalized as in [Lu], \S 3.  Lusztig's result
(see [CP], \S 8.1, and references therein) is that the elements  $ \,
E_{\alpha^k} := T_{i_1} T_{i_2} \cdots T_{i_{k-1}}(E_{i_k}) \, $  and 
$ \, F_{\alpha^k} := T_{i_1} T_{i_2} \cdots T_{i_{k-1}}(F_{i_k}) \, $ 
($ k = 1, \dots, N \, $)  are quantum analogues of root vectors,
respectively of weight  $ +\alpha^k $  and  $ -\alpha^k \, $.  Now
consider the sequence  $ \, (i_1,\dots,i_N) := (1,2,\dots,n-1,n,1,2,
\dots,n-2,n-1,1,2,\dots,3,4,1,2,3,1,2,1) \, $.  Then  $ \, w_0 = s_{i_1}
s_{i_2} s_{i_3} \cdots s_{i_{N-1}} s_{i_N} \, $  is a reduced expression
for  $ w_0 $,  and the previous recipe applies.  It is easy to check
([Ga2], \S 5.3) that, with this choice of the sequence  $ (i_1,\dots,
i_N), $  one has  $ \; E_{\alpha^k} = E_{i,j} \, $,  $ \, F_{\alpha^k}
= F_{j,i} \, $  for  $ \, \alpha^k = \sum_{t=i}^{j-1} \alpha_t \, $, 
\, for all  $ k \, $  (a key fact is that  $ \, E_{\alpha^k} = E_h \, $ 
and  $ F_{\alpha^k} = F_h \, $  if  $ \, \alpha^k = \alpha_h \, $).  So
our choice of quantum root vectors is a special case of Lusztig's.   
                                           \par   
   The quantum version of the PBW theorem claims that the
set of ordered monomials  $ \; B^g := \Big\{\! \prod_{i<j}
E_{i,j}^{\eta_{i,j}} \prod_{k=1}^n G_k^{\gamma_k} \prod_{i<j}
F_{j,i}^{\varphi_{j,i}} \,\Big|\, \eta_{i,j}, \gamma_k, \varphi_{j,i}
\in \N \,\Big\} \, $  (w.r.t.~any total order of the pairs  $ (i,j) $ 
with  $ \, i \not= j \, $)  is a  $ \Qq $--basis  of  $ \uqgl \, $, 
and the set  $ \; B^s := \Big\{\! \prod_{i<j} E_{i,j}^{\eta_{i,j}}
\prod_{k=1}^{n-1} K_k^{\kappa_k} \prod_{i<j} F_{j,i}^{\varphi_{j,i}}
\,\Big|\,   
 \allowbreak   
\eta_{i,j}, \kappa_k, \varphi_{j,i} \in \N \,\Big\} \, $ 
is a  $ \Qq $--basis of  $ \uqsl \, $:  see [CP], \S 8.1,
and references therein.   

\vskip11pt

  {\bf 2.3 Integral forms of  $ \uqgl $  and  $ \uqsl \, $.} \, After
Lusztig's idea for semisimple Lie algebras, adapted to the case of 
$ \gergl_n \, $,  one can define (see [CP], \S 9.3, and references
therein) a suitable  {\sl restricted\/}  integral form of  $ \uqgl $ 
over  $ \Zqqm $,  say  $ \geruqgl \, $:  \, following [DL], \S 3, it
is the unital  $ \Zqqm $--subalgebra  of  $ \uqgl $  generated by 
$ \; F_1^{(m)} $,  $ \dots $,  $ F_{n-1}^{(m)} $,  $ G_1^{\pm 1} $, 
$ \Big(\! {{G_1 \, ; \, c} \atop m} \!\Big) $,  $ \dots $, 
$ G_n^{\pm 1} $,  $ \Big(\! {{G_n \, ; \, c} \atop m} \!\Big) $, 
$ E_1^{(m)} $,  $ \dots $,  $ E_{n-1}^{(m)} \, $,  \, for all  $ \,
m \in \N \, $  and  $ \, c \in \Z \, $.  This  $ \geruqgl $  is a 
$ \Zqqm $--integral  form of  $ \uqgl $  as a Hopf algebra, and
specializes to  $ U_\Z(\gergl_n) $  for  $ \, q \mapsto 1 \, $,  \,
that is  $ \; \geruqgl \Big/ (q-1) \, \geruqgl \cong U_\Z(\gergl_n)
\; $  as co-Poisson Hopf algebras; therefore we call  $ \, \geruqgl
\, $  a  {\it quantum  (\/{\rm or}  quantized) hyperalgebra}. 
Moreover, a suitable PBW-like theorem holds for  $ \geruqgl \, $. 
Finally, for every root of 1, say  $ \varepsilon $,  of odd order 
$ \ell \, $,  a well defined  {\it quantum Frobenius morphism\/} 
$ \; \gerFr_{\gergl_n}^{\,\Zeps} \, \colon \, \geruegl \,
\relbar\joinrel\relbar\joinrel\twoheadrightarrow \, \Zeps \otimes_\Z
U_\Z(\gergl_n) \; $  exists (a Hopf algebra epimorphism): here  $ \,
\geruegl := {\big(\geruqgl\big)}_\varepsilon \, $   --- as in \S 1.4
---   and  $ \gerFr_{\gergl_n}^{\,\Zeps} $  is defined on generators
by dividing out the order of each quantum divided power and each
quantum binomial coefficient by  $ \ell \, $,  if this makes sense,
and mapping to zero otherwise.  The obvious parallel construction  
--- just replacing the  $ G_k^{\pm 1} $'s  by the  $ K_i^{\pm 1} $'s 
everywhere ---   provides  $ \geruqsl \, $,  the  {\sl restricted\/} 
$ \Zqqm $--integral  form of  $ \uqsl \, $:  \, it specializes to 
$ U_\Z(\gersl_n) $  for  $ \, q \mapsto 1 \, $,  \, and  {\it quantum
Frobenius morphisms\/}  $ \; \gerFr_{\gersl_n}^{\,\Zeps} \, \colon \,
\geruesl \, \relbar\joinrel\relbar\joinrel\twoheadrightarrow \, \Zeps
\otimes_\Z U_\Z(\gersl_n) \; $  exist for every  $ \epsilon $  as
above.  See [DL], \S\S 3, 6, or [CP], \S 9.3.   

\vskip3pt

   However, we are mainly interested now in the  {\sl unrestricted\/} 
integral forms and their specializations.  First set, as notation, 
$ \, \overline{X} := \big(q - q^{-1}\big) X \, $.  Following [DP],
\S 12.1, or [Ga1], \S 3.3 --- whose constructions apply to semisimple
Lie algebras, but also work for  $ \gergl_n $  {\sl and\/}  with  $ \Z $ 
instead of  $ \C $  ---   we define  $ \caluqgl $  to be the unital 
$ \Zqqm $--subalgebra  of  $ \uqgl $  generated by  $ \, \big\{\,
\fbar_{j,i} \, , \, G_k^{\pm 1} , \, \ebar_{i,j} \,\big\}_{1 \leq
i < j \leq n}^{k=1,\dots,n} \; $.  This  $ \caluqgl $  is a Hopf 
$ \Zqqm $--subalgebra  of  $ \uqgl $,  and is a free  $ \Zqqm $--module 
with PBW-like basis  $ \; \Cal{B}^g := \Big\{\! \prod_{i<j}
\ebar_{i,j}^{\,\eta_{i,j}} \prod_{k=1}^n G_k^{\gamma_k}
\prod_{i<j} \fbar_{j,i}^{\,\varphi_{j,i}} \,\Big|\,
 \allowbreak   
\eta_{i,j}, \gamma_k, \varphi_{j,i} \in \N \,\Big\} \, $  (monomials
being ordered w.r.t.~any total order as before).  Thus  $ \caluqgl $ 
is another  $ \Zqqm $--integral  form of  $ \uqgl $   --- after 
$ \geruqgl $  ---   as a Hopf algebra, in that it is a Hopf 
$ \Zqqm $--subalgebra  such that  $ \, \Qq \otimes_{\Z[q,q^{-1}]}
\caluqgl \cong \uqgl \, $.   
                                        \par   
   From [DP], \S 14.1, or [Ga1], Theorem 7.4, one has that  $ \caluqgl $ 
is a  {\sl quantization\/}  of  $ \, F\big[{({{}_s{GL}_n}^{\!*})}_\Z\big]
\, $,  \, i.e.~$ \; \caluunogl := {\big(\caluqgl\big)}_1 \cong F
\big[{({}_s{{GL}_n}^{\!*})}_\Z\big] \; $  as Poisson Hopf algebras
(notation of \S\S 1.1/4), where on left-hand side we consider the standard
Poisson structure inherited from  $ \caluqgl \, $.  Finally, let 
$ \ell \, $,  $ \varepsilon $  and  $ \Zeps $  be as in \S 1.4, and 
$ \; \caluegl := {\big( \caluqgl \big)}_\varepsilon \, $:  \, then there
is a Hopf  $ \Z $--algebra  embedding  $ \; \calFr_{\gergl_n}^{\,\Z}
\colon \, F\big[{({{}_s{GL}_n}^{\!*})}_\Z\big] \cong \, \caluunogl
\longhookrightarrow \caluegl \, $,  \; 
%
%
 called
 the  {\it quantum Frobenius morphism\/} 
for  $ \, {{}_s {GL}_n}^{\!*} \, $,  
 \, given by  $ \ell $--th  power on generators 
(see [Ga3], \S\S 3.6--7, for details).  It also
extends to a  $ \Zeps $--linear  map  $ \calFr_{\gergl_n}^{\,\Zeps} $ 
defined on  $ \, F \big[ {({{}_s{GL}_n}^{\!*})}_{\Zeps} \big] \cong
\, \Zeps \otimes_\Z \caluunogl \, $.   
                                        \par   
   Once again, the same constructions can be done with  $ \gersl_n $ 
instead of  $ \gergl_n \, $.  In this way one construct  $ \caluqsl $ 
simply following the recipe above but replacing the  $ G_k^{\pm 1} $'s 
with the  $ K_i^{\pm 1} $'s.  Then  $ \caluqsl $  is another 
$ \Zqqm $--integral  form of  $ \uqsl $  (after  $ \geruqsl $)  as a
Hopf algebra, and it is free as a  $ \Zqqm $--module  with (PBW-like) 
$ \Zqqm $--basis  the set of ordered monomials  $ \; \Cal{B}^s :=
\Big\{\! \prod_{i<j} \ebar_{i,j}^{\,\eta_{i,j}} \prod_{k=1}^{n-1}
K_k^{\kappa_k} \prod_{i<j} \fbar_{j,i}^{\,\varphi_{j,i}} \,\Big|\,
\eta_{i,j}, \kappa_k, \varphi_{j,i} \in \N \,\Big\} \, $.  Entirely
similar results to those for  $ \caluqgl $  then hold for  $ \caluqsl $ 
too, also regarding specializations at roots of 1 (including 1 itself)
and quantum Frobenius    
   \hbox{morphisms  $ \calFr_{\gersl_n}^{\,\Z} $  and 
         $ \calFr_{\gersl_n}^{\,\Zeps} $  (see [Ga3], \S 4).}      
%
%
                                        \par   
   The embedding  $ \; \uqsl \longhookrightarrow \uqgl \; $  restricts
to Hopf embeddings  $ \; \geruqsl \longhookrightarrow \geruqgl \, $  and 
$ \; \caluqsl \longhookrightarrow \caluqgl \, $.  The specializations of
the latter ones at  $ \, q = \varepsilon \, $  and at  $ \, q = 1 \, $ 
are compatible (in the obvious sense) with the quantum Frobenius
morphisms.   

\vskip11pt

  {\bf 2.4  Quantum function algebras, their integral forms and quantum
Frobenius morphisms.} \, Let  $ \fqm $  be the quantum function algebra
over  $ M_n $  introduced by Manin.  Namely (cf.~[No], \S 1.1,
and references therein),  $ \fqm $ is the unital associative 
$ \Qq $--algebra  with generators  $ \; t_{i,j} \; $  ($ \,
i, j = 1, \dots, n \, $)  and relations   
 \vskip-16pt   
  $$  \hbox{ $ \eqalign {
   {} \hskip11pt   t_{i,j} \, t_{i,k} = q \, t_{i,k} \, t_{i,j} \; ,
\quad \quad  t_{i,k} \, t_{h,k}  &  = q \, t_{h,k} \, t_{i,k}
\hskip56pt  \forall  \quad  j<k \, , \, i<h \, ,  \cr
   {} \hskip11pt   t_{i,l} \, t_{j,k} = t_{j,k} \, t_{i,l} \; ,
\qquad  t_{i,k} \, t_{j,l} - \, t_{j,l} \, t_{i,k}  &
= \left( q - q^{-1} \right) \, t_{i,l} \, t_{j,k}   \hskip20pt
\forall  \quad  i<j \, , \, k<l \, .  \cr } $ }   \eqno  $$
The  $ (n \times n) $--matrix  $ \, T := {\big( t_{i,j} \big)}_{i,j=1,
\dots, n} \, $  will be called a  {\sl  $ q $--matrix}.  We call  {\it
``quantum determinant''\/}  the element  $ \; D_q \equiv \text{\sl
det}_q \left(\! {\big(t_{k,\ell}\big)}_{k,\ell=1,\dots,n} \right) :=
{\textstyle \sum_{\sigma \in \Cal{S}_n}} {(-q)}^{l(\sigma)} t_{1{\,}\sigma(1)}
\, t_{2,\sigma(2)} \cdots t_{n,\sigma(n)} \; $,  \; belonging to the
center of  $ \fqm \, $.  This  $ \fqm $  is a bialgebra, with 
$ \; \Delta (t_{i,j}) = \sum\limits_{k=1}^n t_{i,k} \otimes t_{k,j} \; $  and  $ \; \epsilon(t_{i,j}) = \delta_{ij} \; $  (for all  $ i \, $, 
$ j \, $);  \, then  $ D_q $  is group-like, i.e.~$ \; \Delta(D_q) = D_q \otimes D_q \; $  and  $ \; \epsilon(D_q) = 1 \; $.   
%
%
                                        \par
   Finally,  $ \fqm $  admits the set 
 $ \, B_{M_n} := \Big\{\! \prod_{i<j} t_{i,j}^{\tau_{i,j}}
\prod_{k=1}^n t_{k,k}^{\tau_{k,k}} \prod_{i>j} t_{i,j}^{\tau_{i,j}}
\;\Big|\; \tau_{i,j} \in \N \;\, \forall \; i, j \,\Big\} \, $ 
as  $ \Qq $--basis,  where again factors in monomials are ordered
w.r.t.~any fixed total order of the pairs  $ (i,j) $  with  $ \,
i \not= j \, $  (see [Ga4],  Theorem 2.1{\it (a)},  and references
therein).   
                                        \par
   We define  $ \gerfqm $  to be the unital  $ \Zqqm $--subalgebra  of 
$ \fqm $  generated by the  $ t_{i,j} $'s  (for all  $ i $,  $ j \, $); 
\, this is also a  {\sl sub-bialgebra},  and admits  {\sl the same
presentation\/}  as  $ \fqm $  but over  $ \Zqqm \, $.  Therefore  $ \,
D_q \in \gerfqm \, $,  the latter is a free  $ \Zqqm $--module  with 
$ \Zqqm $--basis  $ B_{M_n} $  ([Ga4],  Theorem 2.1{\it (a)\/}),  and 
$ \gerfqm $  is a  {\sl  $ \Zqqm $--integral  form of  $ \fqm \, $}. 
The presentation shows that  $ \gerfqm $  is a  {\sl quantization\/} 
of  $ F\big[(M_n)_{\Z}\big] $,  i.e.~$ \; \gerfunom := {\big( \gerfqm
\big)}_1 \cong F\big[(M_n)_\Z\big] \; $  as Poisson bialgebras (where
on left-hand side we consider the standard Poisson structure inherited
from  $ \gerfqm \, $),  the isomorphism being given by  $ \, t_{i,j}
{\big|}_{q=1} \cong \bar{t}_{i,j} \, $  for all  $ i, j $.  
   Finally let  $ \ell \, $,  $ \varepsilon \, $,  $ \Zeps $  be as
in \S 1.4,  $ \; \gerfem := {\big(\gerfqm\big)}_\varepsilon \, $. 
Then there is a bialgebra 
 embedding   
$ \, \gerFr_{M_n}^{\,\Z} \,
\colon \, F\big[(M_n)_\Z\big] \longhookrightarrow \gerfem \; $ 
$ \big(\, \bar{t}_{i,j} \mapsto t_{i,j}^{\,\ell}
{\big|}_{q=\varepsilon} \,\big) \, $,  \, the 
{\it quantum Frobenius morphism\/}  for  $ M_n \, $.   
%
%

\vskip4pt

   Besides the presentation given above, another description of 
$ \gerfqg $  is available.  In fact, using a characterization as
algebra of matrix coefficients (cf.~[APW], Appendix),  $ \fqm $ 
embeds into  $ \, \uqgl^* \, $  as a dense subalgebra (w.r.t.~the
weak topology).  Hence the natural evaluation pairing  $ \,
\big(\, \langle f, u \rangle \mapsto f(u) \,\big) \, $  restricts
to a perfect (\,=\, non-degenerate) Hopf pairing  $ \; \langle \ ,
\ \rangle \, \colon \fqm \times \uqgl \loongrightarrow \Qq \; $ 
(see e.g.~[No], \S 1.3, for details).  Then, by the analysis in [DL],
\S 4 (adapted to  $ M_n \, $),  one has  $ \; \gerfqm \, = \, \left\{
f \in \fqm \;\Big|\; \big\langle f, \geruqgl \big\rangle \subseteq
\Zqqm \right\} \; $.   
                                             \par   
   All this gives the idea, like in [Ga1], \S 4.3, for semisimple
algebraic groups, to  {\sl define}   
 \vskip-6pt   
   $$  \calfqm \; := \; \left\{ f \in \fqm \;\Big|\; \big\langle f,
\caluqgl \big\rangle \subseteq \Zqqm \right\} \quad .  $$
 \vskip-5pt  
   The arguments in [Ga1], Proposition 5.12,  
%
%
 prove also that  $ \; \Q \cdot \calfqm \; $  is
a  $ \Qqqm $--integral  form of  $ \, \fqm \, $.  Moreover, the
analysis therein together with [Ga2], \S 6, proves that  $ \;
\Q \cdot \calfqm \; $  is a  {\sl quantization\/}  of  $ \, U
\big({\gergl_n}^{\!*}\big) \, $,  \, i.e.~$ \, {\big( \Q \cdot
\calfqm \big)}_1 \cong U \big( {\gergl_n}^{\!*} \big) \, $  as
co-Poisson bialgebras (taking on left-hand side the co-Poisson
structure inherited from  $ \; \Q \cdot \calfqm \, $).  
                                          \par
   Finally, let  $ \ell \, $,  $ \varepsilon $  and  $ \Zeps $  be as
in \S 1.4, and set  $ \; \calfem := {\big(\calfqm\big)}_\varepsilon
\, $.  Then again the arguments in [Ga1], \S 7, show that there is
an epimorphism  
 \vskip-5pt  
  $$  \calFr_{M_n}^{\,\Qeps} \, \colon \, \Qeps \otimes_{\Zeps}
\calfem \llongtwoheadrightarrow \, \Qeps \otimes_\Z U_\Z \big(
{\gergl_n}^{\!\!*} \big) \, = \, \Qeps \otimes_\Q U \big(
{\gergl_n}^{\!\!*} \big)   \eqno (2.1)  $$
 \vskip-3pt  
\noindent   
 of  $ \Qeps $--bialgebras,  which will be called (again) the 
{\sl quantum Frobenius morphism\/}  for  $ M_n \, $.

\vskip7pt

   All the previous construction can be repeated for  $ {GL}_n $ 
instead of  $ M_n \, $  (see again [No], \S 1.1, and references
therein).  One defines  $ \; \fqgl := \big(\fqm\big)\big[{D_q}^{\!
-1}\big] \, $,  \; the extension of  $ \fqm $  by a formal inverse
to  $ D_q \, $;  \,thus the presentation of  $ \fqm $  induces one
of  $ \fqgl $  as well.  This  $ \fqgl $  is a Hopf algebra: the
coproduct and counit on the  $ t_{i,j} $'s  are given by the
formul{\ae}  for  $ \fqm $  (which is then a sub-bialgebra) plus
the formul{\ae}  saying that  $ {D_q}^{\! -1} $  is group-like,
and the antipode by  $ \; S(t_{i,j}) = {(-1)}^{i+j} \,
\text{\sl det}_q \! \left(\! {\big(t_{k,\ell}\big)}_{k, \ell = 1,
\dots, n}^{k \not= i; \ell \not= j} \right) D_q^{-1} $,
$ \; S\big({D_q}^{\! -1}\big) = D_q \, $.  
                                       \par  
%
%
   We define  $ \gerfqgl $  to be the unital  $ \Zqqm $--subalgebra 
of  $ \fqgl $  generated by  $ \; t_{i,j} \; $  (for all  $ i $, 
$ j \, $)  and  $ \, {D_q}^{\! -1} \, $.  This is
%
%
 a  {\sl  $ \Zqqm $--integral  form of  $ \fqgl \, $},  as a
Hopf algebra, and it is a  {\sl quantization\/}  of  $ F \big[
({GL}_n)_{\,\Z} \big] $,  with  $ \, t_{i,j}{\big|}_{q=1} \cong
\bar{t}_{i,j} \, $  (for all  $ i, j \, $)  and  $ \, D_q^{-1}
{\big|}_{q=1} \cong D^{-1} \, $.  
                                               \par    
   It is known (cf.~[Ga4],  Theorem 2.1{\it (a)},  and references
therein), that  $ \gerfqgl $  is free over  $ \Zqqm $,  with (PBW-like)
basis the set of ordered (w.r.t.~any total order) monomials   
%
%
 \vskip-13pt  
  $$  B_{GL_n} \! := \! \Big\{ {\textstyle \prod_{i<j}} t_{i,j}^{\tau_{i,j}}
{\textstyle \prod_{k=1}^n} t_{k,k}^{\tau_{k,k}} {\textstyle \prod_{i>j}}
t_{i,j}^{\tau_{i,j}} \, D_q^{-m} \,\Big|\, \tau_{i,j}, m \in \N \;\,
\forall \, i, j \, ; \, \min \! \big( \{ \tau_{i,i} \big\}_{1 \leq i
\leq n} \! \cup \{m\} \big) = 0 \Big\}  $$
 \vskip-5pt  
%
%

\vskip4pt

   The perfect Hopf pairing  $ \; \langle \ , \ \rangle \, \colon \fqm
\times \uqgl \loongrightarrow \Qq \, $  extends (on left-hand side) to 
$ \fqgl \, $;  \, then one deduces that  $ \; \gerfqgl \, = \, \left\{\,
f \in \fqgl \;\Big|\; \big\langle f, \geruqgl \big\rangle \subseteq
\Zqqm \right\} \, $  (just like for  $ \fqm \, $).  This gives the
idea   --- as in [Ga1], \S 4.3 ---   to  {\sl define}   
 \vskip-3pt
  $$  \calfqgl \; := \; \left\{\, f \in \fqgl \;\Big|\; \big\langle f,
\caluqgl \big\rangle \subseteq \Zqqm \right\} \quad .  $$
 \vskip-3pt
   \indent   Again, one proves that  $ \; \Q \cdot \calfqgl \; $ 
is a  $ \Qqqm $--integral  form of  $ \, \fqgl \, $  (as a Hopf
algebra), and also that  $ \; \Q \cdot \calfqgl \; $  is a  {\sl
quantization\/}  of  $ U\big({\gergl_n}^{\!*}\big) \, $.   
%
%
                                         \par   
   Finally, for  $ \ell $  and  $ \varepsilon $  as in \S 1.4,
let  $ \; \gerfegl := {\big(\gerfqgl\big)}_\varepsilon \, $ 
and  $ \; \calfegl := {\big(\calfqgl\big)}_\varepsilon \, $. 
Then there are two  {\it quantum Frobenius morphisms\/}  for 
$ {GL}_n $   --- both extending the corresponding ones for 
$ M_n $  ---   namely  $ \; \gerFr_{{GL}_n}^{\,\Z} \, \colon
\, F\big[({GL}_n)_{\,\Z}\big] \cong \gerfunogl \longhookrightarrow
\gerfegl \, $,  \; a Hopf algebra monomorphism (given by  $ \;
\bar{t}_{i,j} \mapsto t_{i,j}^{\,\ell}{\big|}_{q=\varepsilon}
\, $,  \,  $ D^{-1} \mapsto D_q^{-\ell} \; $),   
%
%
 \, and a Hopf algebra epimorphism  
 \vskip-3pt
  $$  \calFr_{{GL}_n}^{\,\Qeps} \, \colon \, \Qeps \otimes_{\Zeps}
\calfegl \llongtwoheadrightarrow \, \Qeps \otimes_\Z \calfunogl
\, \cong \, \Qeps \otimes_\Q U\big({\gergl_n}^{\!*}\big) 
\quad .   \eqno (2.2)  $$      

\vskip7pt

   As to  $ {SL}_n $  (see [{\it loc.~cit.}]),  one defines  $ \;
\fqsl := \fqgl \Big/ \! \big( D_q - 1\big) \, \cong \, \fqm \Big/
\! \big( D_q - 1 \big) \, $;  \; here  $ \, \big(D_q - 1\big) \, $ 
is the two-sided ideal of $ \fqgl $  or of  $ \fqm $  generated by
the central element  $ \, D_q - 1 \; $:  \; it is clearly a Hopf ideal,
so  $ \fqsl $  is a Hopf algebra on its own.  Explicitly,  $ \fqsl $ 
admits the same presentation as  $ \fqm $  or  $ \fqgl $  but with
the additional relation  $ \; D_q = 1 \; $  (in either case); its
Hopf structure is given by the same formul{\ae},  but for setting 
$ \; {D_q}^{\! -1} = 1 \; $.
                                       \par  
   Like before, we define  $ \gerfqsl $  as the unital 
$ \Zqqm $--subalgebra  of  $ \fqsl $  generated by the 
$ t_{i,j} $'s:  it has  {\sl the same presentation\/} 
as  $ \fqsl $  {\sl but over\/}  $ \Zqqm \, $,  \, and
is a Hopf subalgebra, hence a  $ \Zqqm $--integral  form
of  $ \fqsl \, $.  In addition we have isomorphisms  $ \;
\gerfqsl \, \cong \, \gerfqgl \Big/ \! \big( D_q - 1 \big)
\, \cong \, \gerfqm \Big/ \! \big(D_q - 1\big) \; $  of
Hopf  $ \Zqqm $--algebras,  and  $ \gerfqsl $  is a  {\sl quantization\/} 
of  $ F \big[ ({SL}_n)_\Z \big] $,  i.e.~$ \; \gerfunosl := {\big(
\gerfqsl \big)}_1 \cong F \big[({SL}_n)_\Z\big] \, $.   
%
%
 It is proved in [Ga4],  Theorem 2.1{\it (a)},  that  $ \gerfqsl $ 
is free over  $ \Zqqm $  with (PBW-like) basis the set    
 \vskip-9pt  
  $$  B_{{SL}_n}  \; := \;  \Big\{ {\textstyle \prod_{i<j}}
t_{i,j}^{\tau_{i,j}} {\textstyle \prod_{k=1}^n} t_{k,k}^{\tau_{k,k}}
{\textstyle \prod_{i>j}} t_{i,j}^{\tau_{i,j}} \;\Big|\; \tau_{i,j}
\in \N \;\, \forall \; i, j \, ; \, \min \{\tau_{k,k}\}_{k=1,\dots,n}
= 0 \,\Big\}  $$     
 \vskip-1pt  
\noindent   
 of ordered monomials (w.r.t.~any total order).   

\vskip4pt

   Like for  $ \fqm $  and  $ \fqgl \, $,  \, also  $ \fqsl $  can be
naturally embedded as a dense subalgebra of  $ \, \uqsls := {\uqsl}^*
\, $.  Then the natural evaluation pairing between  $ {\uqsl}^* $ 
and  $ \uqsl $  restricts to a perfect Hopf pairing  $ \; \langle \ ,
\ \rangle \, \colon \fqsl \times \uqsl \loongrightarrow \Qq \, $;  \,
then (like for  $ \fqm $  and  $ \fqgl \, $)  one finds that  $ \;
\gerfqsl \, = \, \left\{\, f \in \fqsl \;\Big|\; \big\langle f,
\geruqsl \big\rangle \subseteq \Zqqm \right\} \, $.  As in [Ga1],
\S 4.3, this leads us to  {\sl define}   
 \vskip-7pt
  $$  \calfqsl \; := \; \left\{\, f \in \fqsl \;\Big|\; \big\langle f,
\caluqsl \big\rangle \subseteq \Zqqm \right\} \quad .  $$
 Then  $ \; \Q \cdot \calfqsl \; $  is a  $ \Qqqm $--integral 
form of  $ \, \fqsl \, $  (as a Hopf algebra), and it is a 
{\sl quantization\/}  of  $ U\big({\gersl_n}^{\!*}\big) \, $. 
Also, a bialgebra epimorphism  $ \, \calfqm
\relbar\joinrel\relbar\joinrel\twoheadrightarrow
\calfqsl \, $  and a Hopf algebra epimorphism  $ \, \calfqgl
\relbar\joinrel\relbar\joinrel\twoheadrightarrow \calfqsl \, $ 
exist, both dual to  $ \, \caluqsl \longhookrightarrow \caluqgl \, $.   

                                         \par   
   Finally, for any root of 1 of odd order  $ \ell \, $,  say 
$ \varepsilon \, $,  \, there are two  {\it quantum Frobenius
morphisms\/}  for  $ {SL}_n \, $:  \, a Hopf algebra monomorphism 
$ \; \gerFr_{{SL}_n}^{\,\Z} \, \colon \, F\big[({SL}_n)_{\,\Z}\big]
\cong \gerfunosl \llonghookrightarrow \gerfesl \; $  (given by  $ \;
\bar{t}_{i,j} \mapsto t_{i,j}^{\,\ell}{\big|}_{q=\varepsilon} \, $), 
\, and a Hopf algebra epimorphism  
 \vskip-3pt
  $$  \calFr_{{SL}_n}^{\,\Qeps} \, \colon \, \Qeps \otimes_{\Zeps}
\calfesl \llongtwoheadrightarrow \, \Qeps \otimes_\Z \calfunosl
\, \cong \, \Qeps \otimes_\Q U\big({\gersl_n}^{\!*}\big) 
\quad .   \eqno (2.3)  $$      

\vskip11pt

  {\bf 2.5  $ \uqgls $,  $ \uqsls $,  integral forms and quantum
Frobenius morphisms.} \, It is shown in [Ga1], \S\S 5--6, that the
linear dual  $ {\uqsl}^* $  of  $ \uqsl $  can be seen again as a
quantum group on its own: namely, we can set  $ \, \uqsls := {\uqsl}^*
\, $,  \, such a notation being motivated by the fact that  $ \uqsls $ 
stands for the Lie bialgebra  $ {\gersl_n}^{\!\!*} $  exactly like 
$ \uqsl $  stands for  $ \gersl_n \, $.  Indeed,  $ \uqsls $  bears
a natural structure of  {\sl topological\/}  Hopf  $ \Qq $--algebra, 
and has two integral forms  $ \geruqsls $  and  $ \caluqsls $  which
play for  $ \uqsls $  the same r\^{o}le as  $ \geruqsl $  and 
$ \caluqsl $  for  $ \uqsl \, $.  We shall here mainly consider
the first integral form (and less the second), after [Ga1], \S\S
5--6, and apply the same procedure    
   \hbox{to construct  $ \uqgls $  and  $ \geruqgls $  too.}   
\vskip3pt

   First we define  $ \, \H^g_q \, $  as the unital associative 
$ \Qq $--algebra  with generators  $ \; F_1, \dots, F_{n-1} $, 
$ \Lambda_1^{\pm 1}, \dots, \Lambda_{n-1}^{\pm 1}, \Lambda_n^{\pm 1},
E_1, \dots, E_{n-1} \; $  and relations   
  $$  \displaylines { 
   \hfill   \Lambda_i \, \Lambda_i^{-1} \, = \, 1 \, = \, \Lambda_i^{-1}
\, \Lambda_i \, ,  \qquad  \Lambda_i^{\pm 1} \, \Lambda_j^{\pm 1} \, =
\, \Lambda_j^{\pm 1} \, \Lambda_i^{\pm 1}  \qquad   \hfill  \forall \;
i, j  \qquad  \cr
   \hfill   \Lambda_i \, F_j \, \Lambda_i^{-1} \, = \, q^{\delta_{i,j} - \delta_{i,j+1}} \, F_j \, ,  \qquad  \Lambda_i \, E_j \, \Lambda_i^{-1} \,
= \, q^{\delta_{i,j} - \delta_{i,j+1}} \, E_j  \qquad   \hfill  \forall \;
i, j \qquad  \cr
   \hfill   E_i \, F_j \, - \, F_j \, E_i \, = \, 0  \qquad \qquad   \hfill  \forall \; i, j  \qquad  \cr  
 }  $$   
  $$  \displaylines { 
   \hfill   \quad  E_i \, E_j \, = \, E_j \, E_i \, ,  \qquad  F_i \, F_j
\, = \, F_j \, F_i   \hfill  \forall \;\; i, j : \vert i - j \vert > 1 \,
\phantom{.} \;  \cr  
   \hfill   E_i^2 \, E_j \, - \, \left( q + q^{-1} \right) \, E_i \, E_j \,
E_i \, + \, E_j \, E_i^2 \, = \, 0   \hfill  \forall \;\; i, j : \vert
i - j \vert = 1 \, \phantom{.} \;  \cr  
   \hfill   F_i^2 \, F_j \, - \, \left( q + q^{-1} \right) \, F_i \, F_j \,
F_i \, + \, F_j \, F_i^2 \, = \, 0   \hfill  \forall \;\; i, j : \vert i -
j \vert = 1 \, . \;  \cr }  $$   
   \indent   Consider also the elements  $ \; L_h^{\pm 1} :=
\prod_{k=1}^h \Lambda_k^{\pm 1} \, $  (for all  $ \, h = 1,
\dots, n \, $),  and denote  $ \H^s_q $  the   
   unital  $ \Qq $--subalgebra\footnote{The corresponding notation
in [Ga1], \S 6, is  $ \, \H^P_\varphi \, $,  \, and similarly for
integral forms.}  
 of  $ \H^g_q $  generated by 
$ \; \big\{ F_i \, , L_i^{\pm 1}, E_i \big\}_{i=1,\dots,n-1} \; $.  From
the presentation of  $ \H^g_q $  one deduces one of  $ \H^s_q $  too. 
Note also that both  $ \H^g_q $  and  $ \H^s_q $  contain quantum root
vectors  $ E_{i,j} $  and  $ F_{j,i} $  (for all  $ \, i < j \, $) 
as in \S 2.2, and the sets of PBW-like ordered monomials  $ \; B^g_*
:= \Big\{\! \prod_{i<j} E_{i,j}^{\eta_{i,j}} \prod_{k=1}^n
\Lambda_k^{\lambda_k} \prod_{i<j} F_{j,i}^{\varphi_{j,i}} \,\Big|\,
\varphi_{j,i}, \lambda_k, \eta_{i,j} \in \N \,\Big\} \, $ 
for  $ \H^g_q $  and  $ \; B^s_* := \Big\{\! \prod_{i<j}
E_{i,j}^{\eta_{i,j}} \prod_{k=1}^{n-1} L_k^{l_k}
\prod_{i<j} F_{j,i}^{\varphi_{j,i}} \,\Big|\, \varphi_{j,i},
l_k, \eta_{i,j} \in \N \,\Big\} \, $    
       \hbox{for  $ \H^s_q $  are bases over  $ \Qq \, $.}  
                                            \par    
   One defines  $ \uqsls $  as a suitable completion of  $ \H^s_q \, $, 
so that  $ \uqsls $  is a topological  $ \Qq $--algebra  topologically
generated by  $ \H^s_q \, $,  and  $ B^s_* $  is a $ \Qq $--basis  of 
$ \uqsls $  in topological sense.  Then  $ \uqsls $  is also a topological
Hopf  $ \Qq $--algebra (see [Ga1], \S 6).  The same construction makes
sense with  $ \H^g_q $  instead of  $ \H^s_q $  and yields the definition
of  $ \uqgls $,  a topological Hopf algebra with  $ B^s_* $  as
(topological)  $ \Qq $--basis.  By construction  $ \H^s_q $  is a
subalgebra of  $ \H^g_q $  but also a quotient via  $ \, \H^g_q \Big/
\big( L_n - 1 \big) \cong \H^s_q \, $.  Similarly  $ \uqsls $  is a
topological Hopf subalgebra of  $ \uqgls $  but also a quotient via 
$ \, \uqgls \Big/ \big( L_n - 1 \big) \cong \uqsls_q \, $.   
                                                \par    
   Let  $ \; \langle \,\ , \ \rangle \, \colon \, \uqsls \times \uqsl
\llongrightarrow \Qq \; $  be the natural evaluation pairing, given by 
$ \, \langle f, u \rangle := f(u) \, $  for all  $ \, u \in \uqsl \, $, 
$ \, f \in \uqsls := {\uqsl}^* \, $.  Using it, one defines (see [Ga1],
\S 6)  $ \, \geruqsls := \Big\{\, \eta \in \uqsls \;\Big|\; \big\langle 
\eta \, , \, \caluqsl \big\rangle \subseteq \Zqqm \text{\sl \; and  $ \, \Phi (\gamma) \, $  holds} \,\Big\} \, $,  \, where  $ \, \Phi(\gamma) \, $ 
is a suitable ``growth condition'' for elements in  $ \uqsls \, $.  One proves that  $ \geruqsls $  is a  $ \Zqqm $--integral  form
(as a topological Hopf subalgebra) of  $ \uqsls \, $.  To describe 
$ \geruqsls $, let  $ \, \gerH^s_q \, $  be the  $ \Zqqm $--subalgebra 
of  $ \H^s_q $  generated by all  $ F_i^{(m)} $'s,  $ E_i^{(m)} $'s, 
$ L_i^{\pm 1} $'s  and  $ \Big(\! {{L_i \, ; \, c} \atop n} \!\Big) $'s 
($ \, m \in \N \, $,  $ c \in \Z \, $,  $ 1 \leq i < n \, $):  \, then 
$ \, \frak{B}^s_* := \! \Big\{\! \prod_{i<j} \! E_{i,j}^{\,(\eta_{i,j})}
\prod_{k=1}^{n-1} \! \Big(\! {{L_k \, ; \, 0} \atop l_k} \!\Big) \,
L_k^{-\text{\it Ent}(l_k/2)} \prod_{i<j} \!F_{j,i}^{\,(\varphi_{j,i})}
\,\Big|    
 \allowbreak     
   \Big|\, \eta_{i,j}, l_k, \varphi_{j,i} \! \in \! \N \Big\} \, $  is
a basis of  $ \gerH^s_q $  over  $ \Zqqm $,  and  $ \geruqsls $  is the
topological closure of  $ \gerH^s_q \, $,  \, so that  $ \frak{B}^s_* $ 
is a topological  $ \Zqqm $--basis  of  $ \geruqsls $.   
                                                \par    
   Similarly, starting from the pairing  $ \; \langle \,\ , \ \rangle \,
\colon \, \uqgls \times \uqgl \llongrightarrow \Qq \; $  one can define 
$ \, \geruqgls := \Big\{\, \eta \in \uqgls \;\Big|\; \big\langle \eta
\, , \, \caluqgl \big\rangle \subseteq \Zqqm  \text{\sl \; and  $ \,
\Psi(\gamma) \, $  holds} \,\Big\} \, $,  \, where  $ \, \Psi(\gamma) \, $ 
is a suitable ``growth condition'' in  $ \uqgls $.  Then  $ \geruqgls $ 
is a (topological)  $ \Zqqm $--integral  form of  $ \uqgls \, $,  \, which
can be described explicitly: everything goes like for  $ \geruqsls $  but
replacing the symbols  $ L_k $  with the  $ \Lambda_h $'s  ($ h=1, \dots,
n $)  and the  $ \Zqqm $--basis  $ \frak{B}^s_* $  with  $ \, \frak{B}^g_*
:= \Big\{\! \prod_{i<j} E_{i,j}^{\,(\eta_{i,j})} \prod_{h=1}^n \!
\Big(\! {{\Lambda_h \, ; \, 0} \atop \lambda_h} \!\Big) \, \Lambda_h^{-
\text{\it Ent}(\lambda_h/2)} \prod_{i<j} F_{j,i}^{\, (\varphi_{j,i})}
\,\Big|\, \eta_{i,j}, \lambda_h, \varphi_{j,i} \! \in \! \N \Big\} \, $. 
We call  $ \gerH^g_q $  the  $ \Zqqm $--span  of the latter, and 
$ \geruqgls $  is just the completion of  $ \gerH^g_q \, $.  
           Again by construc-\break
 tion  $ \gerH^s_q $  is a subalgebra of  $ \gerH^g_q $  but also a
quotient (restricting  $ \, \H^g_q \Big/ \! \big( L_n \! - \! 1 \big)
\cong \H^s_q \, $  to  $ \gerH^g_q \, $),  \, so  $ \geruqsls $  is a
topological Hopf subalgebra of  $ \geruqgls \, $,  and also quotient
of the latter.   
                                                \par    
   We can describe  $ \gerH^s_q $  rather explicitly: it is the unital
associative  $ \Zqqm $--algebra  with generators  $ \, F_i^{(m)} $, 
$ E_i^{(m)} $,  $ L_i^{\pm 1} $,  $ \Big(\! {{L_i \, ; \, c} \atop m}
\!\Big) \, $  (for  $ \, m \in \N \, $,  $ c \in \Z \, $,  $ i=1, \dots,
n-1 \, $)  and relations   
 \vskip-11pt   
  $$  \displaylines{ 
   {\textstyle \prod\limits_{s=1}^m} \! \big( q^s \! - \! 1 \big)
{{\! L_k ; c} \choose m} \! = \! {\textstyle \prod\limits_{s=1}^m}
\big( q^{1-s+c} L_k -
1 \big) ,  \hskip7pt  L_k L_k^{-1} \! = \! 1 \! = \! L_k^{-1} L_k \, , 
\hskip7pt  {{\! L_h ; c} \choose m} L_k^{\pm 1} \! = \! L_k^{\pm 1}
{{\! L_h ; c} \choose m}  \cr
   \hbox{\sl relations\/  {\rm (1.2)}  for all}  \;\, X \! \in
\! \{L_k\}_{k=1, \dots, n-1} \; ,  \phantom{\Big|}  \hskip5pt 
\hbox{\sl relations\/  {\rm (1.1)}  for all}  \;\, X \! \in \!
\{F_i,E_i\}_{i=1,\dots,n-1}  \cr
   L_h^{\pm 1} F_k^{(m)} = \, q^{\pm \delta_{h,k} m} F_k^{(m)} L_h^{\pm 1}
\; ,  \phantom{\Big|}  \hskip9pt  E_i^{(r)} F_j^{(s)} = \, F_j^{(s)} E_i^{(r)} \; ,  \hskip11pt  L_h^{\pm 1} E_k^{(m)} = \, q^{\pm \delta_{h,k} m} E_k^{(m)} L_h^{\pm 1}  \cr
   {{L_h \, ; \, c} \choose t} \, E_k^{(m)} \, = \, E_k^{(m)} {{L_h \, ; \,
c + \delta_{h,k} \, m} \choose t} \; ,  \qquad  {{L_h \, ; \, c} \choose t} \,
F_k^{(m)} \, = \, F_k^{(m)} {{L_h \, ; \, c + \delta_{h,k} \, m} \choose t} 
\; .  \cr }  $$   
The Hopf structure of  $ \geruqsls $  then can be described explicitly, but we do not need it.   
                                                \par    
   Similarly,  $ \gerH^g_q $  is the unital associative  $ \Zqqm $--algebra  with generators  $ \, F_i^{(m)} $,  $ E_i^{(m)} $,  $ \Lambda_k^{\pm 1} $, 
$ \Big(\! {{\Lambda_k \, ; \, c} \atop m} \!\Big) \, $  (for  $ \, m \in
\N \, $,  $ c \in \Z \, $,  $ i=1, \dots, n-1 $,  $ k=1, \dots, n \, $) 
and relations   
 \vskip-11pt   
  $$  \displaylines{ 
   {\textstyle \prod\limits_{s=1}^m} \! \big( q^s \! - \! 1 \big)
{{\! \Lambda_k ; c} \choose m} \! = \! {\textstyle \prod\limits_{s=1}^m}
\big( q^{1-s+c} \Lambda_k - 1 \big) ,  \hskip7pt  \Lambda_k \Lambda_k^{-1}
\! = \! 1 \! = \! \Lambda_k^{-1} \Lambda_k \, ,  \hskip7pt 
{{\! \Lambda_h ; c} \choose m} \Lambda_k^{\pm 1} \! = \!
\Lambda_k^{\pm 1} {{\! \Lambda_h ; c} \choose m}  \cr
   \hbox{\sl relations\/  {\rm (1.2)}  for all}  \;\,
X \in \{\Lambda_k\}_{k=1, \dots, n} \; ,  \phantom{\Big|}  \hskip5pt 
\hbox{\sl relations\/  {\rm (1.1)}  for all}  \;\,
X \in \{F_i,E_i\}_{i=1,\dots,n-1}  \cr
   E_i^{(r)} F_j^{(s)} = \, F_j^{(s)} E_i^{(r)} \; ,  \phantom{\Big|}  \hskip15pt  \Lambda_h^{\pm 1} \, Y_k^{(m)} = \, q^{\pm (\delta_{h,k} - \delta_{h,k+1}) m} \, Y_k^{(m)} \Lambda_h^{\pm 1}  \hskip11pt  \forall \;\, Y \in \{F,E\}  \cr
   {{\Lambda_h \, ; \, c} \choose t} \, Y_k^{(m)} \, = \, Y_k^{(m)}
{{\Lambda_h \, ; \, c + (\delta_{h,k} - \delta_{h,k+1}) \, m} \choose t}  \hskip13pt  \forall \;\, Y \in \{F,E\} \; .  \cr }  $$   
The Hopf structure of  $ \geruqgls $  can also be given explicitly, yet we do not need that.   

\vskip2pt   

   $ \geruqsls $  is a  {\sl quantization\/}  of  $ U_\Z \big(
{\gersl_n}^{\!\!*} \big) $,  for  $ \, \geruunosls \! := \! {\big(
\geruqsls \big)}_1 \! \cong U_\Z \big({\gersl_n}^{\!\!*}\big) \, $ 
as co-Poisson Hopf algebras, with on left-hand side the co-Poisson
structure inherited from  $ \geruqsls \, $.  In terms of generators
(notation of \S 1.2) this reads  $ \, F_i^{(m)}{\Big|}_{q=1} \! \cong
\text{f}_i^{\,(m)} \, $,  $ \, \Big(\! {{L_i \, ; \, 0} \atop m} \!\Big)
{\Big|}_{q=1} \! \cong \Big(\! {\text{h}_i \atop m} \!\Big) \, $,  $ \,
L_i^{\pm 1}{\Big|}_{q=1} \! \cong 1 \, $,  $ \, E_i^{(m)}{\Big|}_{q=1}
\! \cong \text{e}_i^{(m)} \, $  for all  $ \, i=1, \dots, n-1 \, $ 
and  $ \, m \in \N \, $.   Similarly  $ \, \geruunogls \! := \! {\big(
\geruqgls \big)}_1 \! \cong U_\Z \big({\gergl_n}^{\!\!*}\big) \, $ 
as co-Poisson Hopf algebras, with  $ \, F_i^{(m)}{\Big|}_{q=1} \! \cong
\text{f}_i^{\,(m)} $,  $ \, \Big(\! {{\Lambda_h \, ; \, 0} \atop m}
\!\Big){\Big|}_{q=1} \! \cong \Big(\! {\text{l}_h \atop m} \!\Big) \, $, 
$ \, \Lambda_h^{\pm 1}{\Big|}_{q=1} \! \cong 1 \, $,  $ \, E_i^{(m)}
{\Big|}_{q=1} \! \cong \text{e}_i^{(m)} \, $  for  $ \, m \in \N \, $, 
$ \, i=1, \dots, n-1 \, $,  $ \, h=1, \dots, n \, $.  The map  $ \;
\geruqgls \!\relbar\joinrel\twoheadrightarrow \geruqsls \, $  then is
a quantization of the natural epimorphism  $ \; {\gergl_n}^{\!\!*}
\!\relbar\joinrel\twoheadrightarrow {\gersl_n}^{\!\!*} \, $.   
                                          \par    
   Finally, let  $ \ell $  and  $ \varepsilon $  be as in \S 1.4.  Set 
$ \, \geruesls := {\big(\geruqsls\big)}_\varepsilon \, $  and  $ \,
\gerH^s_\varepsilon := {\big(\gerH^s\big)}_\varepsilon \, $.  Then
(cf.~[Ga1], \S 7.7) the embedding  $ \, \gerH^s_\varepsilon
\longhookrightarrow \geruesls \, $  is an isomorphism, thus  $ \,
\geruesls = \gerH^s_\varepsilon \, $.  Similarly (with like notation) 
$ \, \geruegls = \gerH^g_\varepsilon \, $.  Also, there are Hopf algebra epimorphisms   
  $$  \eqalignno{ 
   \gerFr_{{\gersl_n}^{\!\!*}}^{\,\Zeps} \, \colon \, \geruesls =
\gerH^s_\varepsilon \, \llongtwoheadrightarrow \, \Zeps \otimes_\Z
\gerH^s_1  &  = \Zeps \otimes_\Z \geruunosls \, \cong \, \Zeps
\otimes_\Z U_\Z\big({\gersl_n}^{\!\!*}\big)  \qquad   &  (2.4)  \cr     
   \gerFr_{{\gergl_n}^{\!\!*}}^{\,\Zeps} \, \colon \, \geruegls =
\gerH^g_\varepsilon \, \llongtwoheadrightarrow \, \Zeps \otimes_\Z
\gerH^g_1  &  = \Zeps \otimes_\Z \geruunogls \, \cong \, \Zeps
\otimes_\Z U_\Z\big({\gergl_n}^{\!\!*}\big)  \qquad   &  (2.5)  \cr }  $$   
defined by  $ \, X_i^{(s)}\Big\vert_{q=\varepsilon} \!\!\! \mapsto
\text{x}_i^{\, (s / \ell)} \, $,  $ \, \Big(\! {{Y_j \, ; \, 0} \atop s}
\!\Big) \! \Big\vert_{q=\varepsilon} \!\!\! \mapsto \Big(\! {{\text{k}_j}
\atop {s / \ell}} \!\Big) \, $  if  $ \, \ell \Big\vert s \, $,  $ \,
X_i^{(s)} \Big\vert_{q=\varepsilon} \!\!\! \mapsto 0 \, $,  $ \, \Big(\!
{{Y_j \, ; \, 0} \atop s} \!\Big) \! \Big\vert_{q=\varepsilon} \!\!\!
\mapsto 0 \, $  if  $ \, \ell \hbox{$ \not\big\vert $}  s \, $,  \, and 
$ \, Y_j^{\pm 1} \Big\vert_{q=1} \!\!\! \mapsto 1 \, $,  \, with  $ \,
\big( X, \text{x} \big) \in \big\{\! \big(F,\text{f}\,\big), \big( E,
\text{e} \big) \!\big\} \, $,  and  $ \, \big(Y,\text{k}\big) := \big(
L, \text{h} \big) \, $  in the  $ \gersl_n $  case,  $ \, \big( Y,
\text{k} \big) := \big( \Lambda, \text{g} \big) \, $  in the 
$ \gergl_n $  case.  These are   
         \hbox{{\it quantum Frobenius morphism\/}  for 
$ {\gersl_n}^{\!\!*} $  and  $ {\gergl_n}^{\!\!*} $.}   
                                           \par   
   The epimorphism  $ \; \pi_q \, \colon \, \geruqgls
\longtwoheadrightarrow \geruqsls \, $  of topological
Hopf  $ \Zqqm $--algebras  mentioned above is compatible
with the quantum Frobenius morphisms, that is to say 
$ \, \pi_1 \circ \gerFr_{{\gergl_n}^{\!\!*}}^{\,\Zeps} =
\gerFr_{{\gersl_n}^{\!\!*}}^{\,\Zeps} \circ \pi_\varepsilon
\, $  (where  $ \pi_1 $  and  $ \pi_\varepsilon $  have
the obvious meaning).   
                                                \par   
   We remark that  $ \, \gerH^s_q = \big\{\, \eta \in \H^s_q \,\big|\,
\big\langle \eta \, , \, \caluqsl \big\rangle \subseteq \Zqqm \big\}
\, $,  \, and similarly for  $ \gerH^g_q \, $.  Then one can consider
([Ga1], \S\S 6--7)  $ \, \calH^s_q := \big\{\, \eta \in \H^s_q \,\big|\,
\big\langle \eta \, , \, \geruqsl \big\rangle \subseteq \Zqqm \big\}
\, $.  This is the  $ \Zqqm $--subalgebra  of  $ \H^s_q $  generated
by  $ \, \big\{\, \ebar_{i,j} \, , L_h^{\pm 1} \, , \fbar_{j,i}
\big\}_{1 \leq i < j \leq n}^{1 \leq h \leq n-1} \, $,  \, with a
PBW-like basis over  $ \Zqqm $  given by the ordered monomials in
these generators.  The unrestricted  $ \Zqqm $--integral  form 
$ \caluqsls $  is the (suitable) topological completion of 
$ \calH^s_q $  inside  $ \uqsls $;  one also has another  {\it
quantum Frobenius morphism}  $ \; \calFr_{\gergl_n^*}^{\,\Z} \,
\colon \, \calU_1(\gerg^*) \, \longhookrightarrow \, \caluegs \, $, 
\; defined on generators as an  ``$ \ell $-th  power operation''. 
Similarly holds for  $ \, \calH^g_q := \big\{\, \eta \! \in \!
\H^g_q \,\big|\, \big\langle \eta \, , \, \geruqgl \big\rangle
\subseteq \Zqqm \!\big\} \, $  and  $ \caluqgls \, $,  with 
$ \Lambda_k $'s  instead of  $ L_h $'s,  and a suitable
quantum Frobenius  $ \, \calFr_{\gersl_n^*}^{\,\Z} \, $.   

\vskip0,5pt   

  {\it  $ \underline{\hbox{\it Warnings}} $:}  {\it (a)\/}  In [Ga1]
integral forms are over  $ \Qqqm $,  but now we work over  $ \Zqqm $. 
Some results here for  $ \calfqsl $  are improvements of those
for  $ \, \Qqqm \otimes_{\Z\,[q,q^{-1}]} \calfqsl \, $  in [Ga1],
and similarly for  $ {GL}_n \, $  (see \S\S 4.8--10).  For  $ \uqgs $ 
too, the results in [Ga1] are over  $ \Qqqm $:  yet the arguments 
therein also apply over  $ \Zqqm $,  giving the results of \S 2.5.   
                                                      \par  
   {\it (b)} \, In this work we deal with integral forms  $ \calF_q[X] $ 
of  $ F_q[X] $  (for  $ \, X \in \big\{ M_n \, , {GL}_n \, , {SL}_n
\big\} \, $),  passing through integral forms  $ \frak{U}_q\big(
\gerg^* \big) $  of  $ \uqgs $  (for  $ \, \gerg \in \big\{ \gergl_n
\, , \gersl_n \big\} \, $).  The same methods and ideas apply also
to the forms  $ \gerF_q[X] $,  passing through the forms  $ \calU_q
\big(\gerg^*\big) $:  this yields new proofs of the results for them
and their quantum Frobenius morphisms mentioned in \S\S 2.4.   

\vskip0,87truecm

\centerline {\bf \S\; 3 \ Embedding quantum function algebras into
dual quantum groups. }   

\vskip6pt

  {\bf 3.1 Construction of  $ \uqgs $  and the natural evaluation
pairing.} \, The description of  $ \, \uqgs := {\uqg}^* \, $,  for 
$ \, \gerg \in \big\{ \gergl_n , \gersl_n \big\} \, $,  given above
comes from the following construction (cf.~[Ga1], \S\S 5--6).  Let 
$ U_q(\gerb_+) $  and  $ U_q(\gerb_-) $  be the  {\sl quantum Borel
subalgebras\/}  of  $ \uqg \, $:  \, namely,  $ U_q(\gerb_+) \, $, 
resp.~$ U_q(\gerb_-) \, $,  is the subalgebra of  $ \uqg $  obtained
discarding the  $ F_i $'s,  resp.~the  $ E_i $'s,  from the set of
generators; similarly one defines integral forms  $ \gerU_q(\gerb_\pm) $ 
and  $ \calU_q(\gerb_\pm) $  of  $ U_q(\gerb_\pm) $  as well.  Then 
$ \uqg $  as a coalgebra is a quotient of  $ \, U_q(\gerb_+) \otimes
U_q(\gerb_-) $.  Therefore  $ \, \uqgs := {\uqg}^* \, $  is a subalgebra
of  $ \, {\big( U_q(\gerb_+) \otimes U_q(\gerb_-) \big)}^* \! \cong
{U_q(\gerb_+)}^* \widehat\otimes \, {U_q(\gerb_-)}^* \, $,  \, where 
$ \, \widehat\otimes \, $  denotes topological tensor product (completion
of  $ \otimes $  w.r.t.~weak topology).  In [Ga1], \S 5, it was observed
that  $ {U_q(\gerb_\pm)}^* $  as an algebra is the completion of  $ U_q
(\gerb_\mp) \, $,  \, so that  $ \, {\big( U_q(\gerb_+) \otimes U_q
(\gerb_-) \big)}^* \! \cong U_q(\gerb_-) \, \widehat\otimes \,\, U_q
(\gerb_+) \, $.  Then  $ \uqgs $  can be described as topologically
generated by elements  $ \, F_i \otimes 1 \, $,  $ \, T_k^{\mp 1}
\otimes T_k^{\pm 1} \, $,  $ \, 1 \otimes E_i \in U_q(\gerb_-) \,
\widehat\otimes \,\, U_q(\gerb_+) \cong {\big( U_q(\gerb_+) \otimes
U_q(\gerb_-) \big)}^* \, $  (with  $ \, T_k^{\pm 1} \in \big\{
G_j^{\pm 1} \big\}_{1 \leq j \leq n} \, $  if  $ \, \gerg = \gergl_n
\, $  and  $ \, T_k^{\pm 1} \in \big\{ K_h^{\pm 1} \big\}_{1 \leq h
\leq n-1} \, $  if  $ \, \gerg = \gersl_n \, $).   
                                                 \par   
   Now we observe that  $ \, U_q(\gerb_+) \cong U_q(\gerb_-) \, $  as
algebras, via  $ \, E_i \mapsto F_i \, $,  $ \, T_k^{\pm 1} \mapsto
T_k^{\mp 1} \, $.  Then we have also  $ \, {\big( U_q(\gerb_+) \otimes
U_q(\gerb_-) \big)}^* \! \cong U_q(\gerb_-) \, \widehat\otimes \,\,
U_q(\gerb_+) \cong U_q(\gerb_+) \, \widehat\otimes \,\, U_q(\gerb_-)
\, $.  Thus  $ \, \uqgs := {\uqg}^* \, $  can be described as
topologically generated by elements  $ \, E_i \otimes 1 \, $, 
$ \, T_k^{\pm 1} \otimes T_k^{\mp 1} \, $,  $ \, 1 \otimes F_i  
 \allowbreak    
   \in U_q(\gerb_+) \, \widehat\otimes \,\, U_q(\gerb_-) \, $,  \, which
by abuse of notation we shall call  $ \, E_i \, $,  $ \, T_k^{\pm 1} \, $ 
and  $ \, F_i \, $.  The relations among them are those for  $ E_i \, $, 
$ T_k^{\pm 1} $  (in  $ \big\{\, G_j^{\pm 1} \,\big\}_{1 \leq j
\leq n} $  or in  $ \big\{\, K_h^{\pm 1} \,\big\}_{1 \leq h \leq n-1} $ 
according to whether  $ \, \gerg = \gergl_n \, $  or  $ \, \gerg =
\gersl_n \, $) and  $ F_i $  in \S 2.5.  Using this description for 
$ \uqgs \, $,  the evaluation pairing  $ \; \langle \,\ , \ \rangle
\, \colon \uqgs \times \uqg \longrightarrow \Qq \; $  is uniquely
determined by its values on PBW bases.  It is described explicitly
in [Ga1] (\S 2.4 and \S 6), but before identifying  $ \, U_q(\gerb_+)
\cong U_q(\gerb_-) \, $  and  $ \, U_q(\gerb_-) \cong U_q(\gerb_+)
\, $.   
%
%
 Once we adopt these last identifications, which map  $ \; E_{i,j}
\mapsto {(-q)}^{i-j+1} F_{j,i} \; $  and  $ \; F_{j,i} \mapsto
{(-q)}^{j-i-1} E_{i,j} \; $  (for  $ \, i < j \, $),  the results
in [Ga1], \S 6, eventually read   
 \vskip-11pt
  $$  \hbox{ $ \eqalign{ 
         {\textstyle \bigg\langle 
   \prod\limits_{i<j} \! E_{i,j}^{\,(\eta_{i,j})}
\prod\limits_{r=1}^l \! \Big(\! {{Y_r ; \, 0} \atop \chi_r} \!\Big)
\, Y_r^{-\text{\it Ent}(\chi_r/2)}
\prod\limits_{i<j} \! F_{j,i}^{\,(\varphi_{j,i})} \, , \;
   \prod\limits_{i<j} \! \ebar_{i,j}^{\,e_{i,j}}
\prod\limits_{s=1}^l \! Z_s^{\!\phantom{|}z_s}
\prod\limits_{i<j} \! \fbar_{j,i}^{\,f_{j,i}} \bigg\rangle } \; = 
\qquad  &   \hfill  \cr   
   \hfill   = \; {\textstyle \prod\limits_{i<j}} {(-1)}^{\eta_{i,j}}
\, \delta_{f_{j,i},\varphi_{j,i}} \, \delta_{e_{i,j},\eta_{i,j}}
\, {(-q)}^{(j-i-1) (\varphi_{j,i} - \eta_{i,j})} \,
{\textstyle \prod\limits_{t=1}^l} {z_t \choose \chi_t}_{\!\!q} \,
q^{- z_t \text{\it Ent}\,(\chi_t/2)}  &  \cr } $ }   \eqno (3.1)  $$   
where  $ \big( Y, Z, l \big) $  stands for  $ \big( L, K, n-1 \big) $ 
in the  $ \gersl_n $  case and for  $ \big( \Lambda, G, n \big) $  in
the  $ \gergl_n $  case.   

\vskip11pt

  {\bf 3.2 The embedding  $ \, \xi \, \colon \, \fqg \longhookrightarrow
\uqgs \, $  and its restrictions to integral forms.} \, Let  $ \, G \in
\big\{ {SL}_n, {GL}_n \big\} \, $  and  $ \, \gerg := \Cal{L}\text{\it
ie}\,(G) \, $.  Definitions embed  $ \fqg $  into  $ \, {\uqg}^*
=: \uqgs \, $,  \, via a monomorphism  $ \; \xi \, \colon \, \fqg
\longhookrightarrow \uqgs \; $  of topological Hopf  $ \Qq $--algebras:
they also imply  $ \, \calfqg = \xi^{-1} \big(\geruqgs\big) \, $, 
\, thus  $ \xi $  restricts to a monomorphism  $ \; \widehat{\xi} \,
\colon \, \calfqg \longhookrightarrow \geruqgs \, $,  \, and similarly
to  $ \, \widetilde{\xi} \, \colon \, \gerfqg \! \longhookrightarrow
\! \caluqgs \, $  too.  These verify  $ \, \fqg = \xi^{-1}(\H_q^x) \, $, 
$ \, \calfqg \, = \, \widehat{\xi}^{-1}\big(\gerH_q^x\big) \, $  and 
$ \, \gerfqg \, = \, \widetilde{\xi}^{-1}\big(\calH_q^x\big) \, $ 
(with  $ \, x = s \, $  or  $ \, x = g \, $  according to the type
of  $ G \, $).  As  $ \fqm $  embeds into  $ \fqgl $   --- and
similarly for integer forms ---   one has also bialgebra embeddings 
$ \; \xi \, \colon \fqm \longhookrightarrow \uqgls \, $,  $ \;
\widehat{\xi} \, \colon \, \calfqm \longhookrightarrow \geruqgls
\; $  and  $ \; \widetilde{\xi} \, \colon \, \gerfqm
\longhookrightarrow \caluqgls \; $.   
                                            \par
   Moreover, by construction  $ \widehat{\xi} $  and  $ \widetilde{\xi} $ 
are compatible with specializations and quantum Frobenius morphisms,
{i.e.}  $ \Big( \text{id}_{\Qeps} \otimes_{\Z} \widehat{\xi}
{\,\big|}_{q=1} \Big) \circ \, \calFr_G^{\,\Qeps} = \Big(
\text{id}_{\Qeps} \otimes_{\Zeps} \gerFr_{\gerg^*}^{\,\Zeps}
\Big) \circ \Big( \text{id}_{\Qeps} \otimes_{\Zeps} \widehat{\,\xi}
{\big|}_{q=\varepsilon} \Big) \; $  for  $ \widehat{\xi} $,  and
similarly for  $ \widetilde{\xi} \, $  (for any root of unity 
$ \varepsilon $  of odd order).  Finally, all the embeddings 
$ \xi $  and their restrictions to integral forms are compatible
(in the obvious sense) with the epimorphisms  $ \, \uqgls {\buildrel
\pi \over {\relbar\joinrel\relbar\joinrel\twoheadrightarrow}\,}
\uqsls \, $,  $ \, \fqgl \, {\buildrel \pi \over
{\relbar\joinrel\relbar\joinrel\twoheadrightarrow}}
\, \fqsl \, $  or  $ \, \fqm \, {\buildrel \pi \over
{\relbar\joinrel\relbar\joinrel\twoheadrightarrow}}
\, \fqsl \, $  and their restrictions to integral
forms, see Proposition 3.7{\it (b)}.  All this follows
from definitions (see also [Ga1], \S 6).   

\vskip4pt

   The following gives an explicit description of all the
embeddings mentioned above:   

\vskip11pt

\proclaim{Proposition 3.3} The embeddings  $ \; \widetilde{\xi}
\, \colon \, \gerfqm \longhookrightarrow \caluqgls \, $, 
$ \; \widetilde{\xi} \, \colon \, \gerfqgl \longhookrightarrow
\caluqgls \; $,  $ \; \xi \, \colon \, \fqm \longhookrightarrow
\uqgls \; $  and  $ \; \xi \, \colon \, \fqgl \longhookrightarrow
\uqgls \; $  are uniquely determined by   
  $$  \eqalignno{ 
   t_{i,j}  \; \mapsto \;  &  -{(-q)}^{j-i-1} \, \ebar_{i,j} \,
\Lambda_j \, - \, {(-q)}^{j-i-2} \, {\textstyle \sum\limits_{k=j+1}^n}
\ebar_{i,k} \, \Lambda_k \, \fbar_{k,j}  \qquad  &  
\forall \;\; i < j  \cr  
   t_{l,l}  \; \mapsto \;  &  \; \Lambda_l \, - \, q^{-2} \,
{\textstyle \sum\limits_{k=l+1}^n} \ebar_{l,k} \, \Lambda_k
\, \fbar_{k,l}  \qquad  &   \forall \;\; i = l = j  \cr  
   t_{i,j}  \; \mapsto \;  &  +{(-q)}^{j-i-1} \, \Lambda_i
\, \fbar_{i,j} \, - \, {(-q)}^{j-i-2} \, {\textstyle
\sum\limits_{k=i+1}^n} \ebar_{i,k} \, \Lambda_k \,
\fbar_{k,j}  \qquad  &   \forall \;\; i > j  \cr }  $$   
With  $ \, \Lambda_k^{\pm 1} = L_k^{\pm 1} L_{k-1}^{\mp 1} \, $, 
this describes  $ \, \widetilde{\xi} \, \colon \gerfqsl \hookrightarrow
\caluqsls \, $  and  $ \, \xi \, \colon \fqsl \hookrightarrow
\uqsls \, $.   
\endproclaim  

\demo{Proof} It is enough to prove the claims for the 
$ \widetilde{\xi} $'s:  the rest follows by scalar extension. 
                                      \par   
   Let  $ \, G := {GL}_n \, $,  \, let $ \gerF_q[B_+] $  and 
$ \gerF_q[B_-] $  be the quotient Hopf algebras of  $ \gerfqg $ 
obtained factoring out the generators  $ t_{i,j} $  with  $ \,
i \! > \! j \, $  and  $ \, i \! < \! j \, $  respectively.  Let 
$ \, \pi_+ \, \colon \, \gerfqg \relbar\joinrel\twoheadrightarrow
\gerF_q[B_+] \, $  and  $ \, \pi_- \, \colon \, \gerfqg
\relbar\joinrel\twoheadrightarrow \gerF_q[B_-] \, $  be
the corresponding epimorphisms.   
                                      \par   
   In [DL], \S 4.3, an embedding  $ \, \widetilde{\xi} \, \colon \,
\gerfqg \longhookrightarrow \caluqgs \, $  is given by the composition   
  $$  \gerfqg \, @>{\;\;\Delta\;\;}>> \, \gerfqg \otimes \gerfqg 
\, @>{\;\; \pi_+ \otimes \, \pi_- \;\;}>> \, \gerF_q[B_+] \otimes
\gerF_q[B_-] \, @>{\;\; \vartheta_+ \otimes \, \vartheta_- \;\;}>>
\, \calU_q(\gerb_-) \otimes \calU_q(\gerb_+)  $$   
whose image is contained in the subalgebra of  $ \, \calU_q(\gerb_-)
\otimes \calU_q(\gerb_+) \, $  which is isomorphic (through the
identifications explained in \S 3.1) to the subalgebra  $ \calH^g_q $ 
of  $ \caluqgs $  considered at the end of \S 2.5.  The last step in
this composition is given by 
%
%
 isomorphisms  $ \, \vartheta_+ \, \colon
\, \gerF_q[B_+] \cong \calU_q(\gerb_-) \, $  and  $ \, \vartheta_- \,
\colon \, \gerF_q[B_-] \cong \calU_q(\gerb_+) \, $.  We modify this
construction (as sketched in \S 3.1) by changing the last step with 
$ \, \gerF_q[B_+] \otimes \gerF_q[B_-] \, @>{\;\; \vartheta_+ \otimes
\, \vartheta_- \;\;}>> \, \calU_q(\gerb_+) \otimes \calU_q(\gerb_-)
\, $,  where  $ \, \vartheta_+ \, \colon \, \gerF_q[B_+] \cong \calU_q
(\gerb_+) \, $  and  $ \, \vartheta_- \, \colon \, \gerF_q[B_-] \cong
\calU_q(\gerb_-) \, $  are the algebra isomorphisms given by
  $$  \matrix 
     \hskip-27pt  \vartheta_+ \, : \hskip-3pt  &  t_{i,k} \, \mapsto
\, +{(-q)}^{k-i} \Lambda_k^{+1} \ebar_{i,k}  \quad  \big(\, \forall
\;\, i < k \,\big) \, ,  \quad  &  t_{l,l} \, \mapsto \, \Lambda_l^{+1} 
\quad  \big(\, \forall \;\, l \,\big)  \\   
     \hskip-27pt  \vartheta_- \, : \hskip-3pt  &  t_{k,j} \, \mapsto
\, -{(-q)}^{j-k} \fbar_{k,j} \, \Lambda_k^{-1}  \quad  \big(\, \forall
\;\, k > j \,\big) \, ,  \quad  &  t_{l,l} \, \mapsto \, \Lambda_l^{-1} 
\quad  \big(\, \forall \;\, l \,\big)  \endmatrix  $$    
   \indent   Therefore, explicit computation gives (using notation 
$ \, \delta_{i \sim j} := 1 \, $  or  $ \, \delta_{i \sim j} := 0
\, $  according to whether  $ \, i \sim j \, $  or not, for any
relation  $ \sim \, $)  
  $$  \displaylines{ 
   \quad  \big( (\vartheta_+ \otimes \vartheta_-) \circ (\pi_+ \otimes
\pi_-) \circ \Delta \big)(t_{i,j}) \; = \;   
%
%
   \delta_{i=j} \, \Lambda_i^{+1} \otimes \Lambda_i^{-1} \, - \,
\delta_{i>j} \, {(-q)}^{j-i} \, \Lambda_i^{+1} \otimes \fbar_{i,j}
\, \Lambda_i^{-1} \, +   
 \hfill   
        \cr   
   \hfill   + \;\, \delta_{i<j} \, {(-q)}^{j-i} \, \Lambda_j^{+1}
\, \ebar_{i,j} \otimes \Lambda_j^{-1} \; - \, {\textstyle
\sum_{k = (i \vee j) + 1}^n} {(-q)}^{j-i} \, \Lambda_k^{+1} \,
\ebar_{i,k} \otimes \fbar_{k,j} \, \Lambda_k^{-1}  \cr }  $$   
Exploiting the commutation relations in  $ \calU_q(\gerb_\pm) $  and
using notation  $ \, E_{i,\ell} := E_{i,\ell} \otimes 1 \, $,  $ \,
\Lambda_\ell^{\pm 1} := \Lambda_\ell^{\pm 1} \otimes \Lambda_\ell^{\mp 1} 
\, $  and  $ \, F_{\ell,j} := 1 \otimes F_{\ell,j} \, $  as in \S 3.1,
the above result     
%
%
 gives the claim.   \hskip35pt \hfill \qed\break   
\enddemo   

\vskip-3pt

  {\bf 3.4 Remarks:}  {\it (a)} \, The embedding considered
in [GR] for  $ \, n = 2 \, $  is  
  $$  \widetilde{\xi} \; \colon \quad 
t_{1,1} \, \mapsto \, \Lambda_1 - \fbar \Lambda_2 \ebar \; , 
\quad  t_{1,2} \, \mapsto \, - \fbar \Lambda_2 \; ,  \quad 
t_{2,1} \, \mapsto \, + \Lambda_2 \ebar \; ,  \quad 
t_{2,2} \, \mapsto \Lambda_2 \; ,  $$   
after the recipe of [DL].  In the sequel we exploit some results
from [GR], but these only use intrinsic properties of  $ F_q[{SL}_2]
\, $,  independent of any embedding of the latter into  $ U_q \big(
\gersl_2 \big) \, $.  Thus the present work can be applied (in a non
self-contradicting way) for  $ \, n = 2 \, $  too.   
                                                \par    
   {\it (b)} \, In the following to simplify notation we shall 
identify  $ \fqm $  and  $ \fqg $  with their isomorphic images
in  $ \uqgs $  via the embeddings  $ \xi \, $,  \, and the same
for integral forms.   
%
%

\vskip11pt

  {\bf 3.5 From inclusions to identities.} \, Let  $ \, G := {SL}_n
\, $.  By [DL], \S 4, or [Ga1], Theorem 5.14, one knows that the
monomorphisms  $ \, \widetilde{\xi} \, \colon \, \gerF_q[G]
\longhookrightarrow \caluqgs \, $,  $ \, \xi \, \colon \, F_q[G]
\longhookrightarrow \uqgs \, $  and  $ \, \widehat{\xi} \, \colon
\, \calF_q[G] \longhookrightarrow \geruqgs \, $  extend to
isomorphisms (of algebras)  $ \; \widetilde{\xi} \, \colon \,
\gerF_q[G] \big[ \phi^{-1} \big] \, {\buildrel \cong \over
{\lhook\joinrel\relbar\joinrel\relbar\joinrel\twoheadrightarrow}}
\, \caluqgs \, $,  $ \; \xi \, \colon \,
F_q[G]\big[\phi^{-1}\big] \, {\buildrel \cong \over
{\lhook\joinrel\relbar\joinrel\relbar\joinrel\twoheadrightarrow}}
\, \uqgs \; $  and  $ \; \widehat{\xi} \, \colon \,
\calF_q[G]\big[\phi^{-1}\big] \, {\buildrel \cong \over
{\lhook\joinrel\relbar\joinrel\relbar\joinrel\twoheadrightarrow}}
\, \geruqgs \, $,  \; for some  $ \, \phi \in \gerF_q[G] \, $.   
%
%
                                       \par   
   We shall now improve this result, with an independent approach,
and extend it to  $ \, G = {GL}_n \, $  and  $ \, M_n \, $. 
Mimicking [DL], \S 4,   
%
%
   the would-be element  $ \, \phi \in \gerfqgl \, $  should be given by  
  $$  \phi \, = \, T_1 \, T_2 \cdots T_{n-1} \, T_n \, = \, ( \Lambda_1
\cdots \Lambda_n) \, (\Lambda_2 \cdots \Lambda_n) \cdots (\Lambda_{n-1}
\, \Lambda_n) \, \Lambda_n \, = \, \Lambda_1 \, \Lambda_2^{\,2} \,
\Lambda_3^{\,3} \cdots \Lambda_{n-1}^{\,n-1} \, \Lambda_n^{\,n}  $$   
where  $ \, T_l := \Lambda_l \, \Lambda_{l+1} \, \cdots \Lambda_{n-1}
\, \Lambda_n \, \in \, \calH_q^g \, $  (see \S 2.5)  for  $ \, l = 1,
2, \dots, n \, $,  and we identify quantum function algebras with
their images in  $ \uqgls $  as in Proposition 3.3.  The key step
is      

\vskip11pt   

\proclaim{Lemma 3.6} \, Let's identify  $ \fqm $  with its copy 
$ \xi\big(\fqm\big) $  inside  $ \uqgls \, $.  Then   
  $$  \phi \, := \, T_1 \, T_2 \cdots T_{n-1} T_n \, = \, \Lambda_1
\Lambda_2^{\,2} \cdots \Lambda_{n-1}^{\,n-1} \Lambda_n^{\,n} \; \in
\; \gerfqm \;\; \big( \subset \calfqm \subset \fqm \,\big)  $$   
\endproclaim  

\demo{Proof} Exploiting carefully the  formul{\ae}  of Proposition 3.3,
one easily proves that   
  $$  T_l \, := \; \Lambda_l \, \Lambda_{l+1} \, \cdots \Lambda_{n-1}
\, \Lambda_n \, \in \, \gerfqm   \qquad  \hbox{for all  $ \, l = 1, 2,
\dots, n, n+1 \, $}   \eqno (3.2)  $$   
by a simple induction, whence the claim follows at once.   \hskip135pt \hfill \qed\break   
\enddemo   

\vskip-3pt   

   Now using this lemma (and Corollary 3.8 later on) we get our
``improved result'':   
%
%
%
%
 \eject   

\proclaim{Theorem 3.7}  The embeddings  $ \widetilde{\xi} \, $, 
$ \widehat{\xi} $  and  $ \xi $  extend to the following
identifications:   
  $$  \matrix  
   \text{\it (a)}  &  \hskip13pt  \gerfqm\big[\phi^{-1}\big] \, =
\; \calH_q^g \;\; ,  & \quad  \calfqm\big[\phi^{-1}\big] \, = \;
\gerH_q^g \;\; ,  & \quad  \fqm\big[\phi^{-1}\big] \, = \; \H_q^g 
&  \hskip5pt  \\  
   \text{\it (b)}  &  \phantom{\Big|}  \hskip11pt  \gerfqgl\big[\phi^{-1}\big] \, = \; \calH_q^g \;\; ,  & \quad  \calfqgl\big[\phi^{-1}\big]
\, = \; \gerH_q^g \;\; ,  & \quad  \fqgl\big[\phi^{-1}\big] \, = \; \H_q^g 
&  \hskip5pt  \\  
   \text{\it (c)}  &  \phantom{\big|}  \hskip11pt  \gerfqsl\big[\phi^{-1}\big] \, =
\; \calH_q^s \;\; ,  & \quad  \calfqsl\big[\phi^{-1}\big] \, = \;
\gerH_q^s \;\; ,  & \quad  \fqsl\big[\phi^{-1}\big] \, = \; \H_q^s 
&  \hskip5pt  
      \endmatrix  $$   
where in the last row we write again  $ \, \phi \, $  for the image of 
$ \, \phi \in \! \fqm \, $  via  $ \, \fqm \, {\buildrel \pi \over {\relbar\joinrel\relbar\joinrel\twoheadrightarrow}} \, \fqsl \, $.   
\endproclaim   

\demo{Proof}  {\it (a)} \, It is clear that having  $ \, T_1 \, ,
T_2 \, , \dots , T_{n-1} \, , T_n \in \gerfqm \subseteq \gerfqm
\big[\phi^{-1}\big] \, $  together with  $ \, \phi^{-1} = T_1^{-1}
\, T_2^{-1} \cdots T_{n-1}^{\,-1} \, T_n^{-1} \in \gerfqm
\big[\phi^{-1}\big] \; $  imply   
  $$  T_l^{-1} = \Lambda_l^{-1} \, \Lambda_{l+1}^{\,-1} \,
\cdots \Lambda_{n-1}^{\,-1} \, \Lambda_n^{-1} \, \in \,
\gerfqm\big[\phi^{-1}\big]  \qquad  \forall \;\; l = 1,
2, \dots, n, n+1 \; .   \eqno (3.3)  $$   
   \indent   From  formul{\ae}  (3.2--3) we easily find that 
$ \; \Lambda_1^{\pm 1} \, $,  $ \, \Lambda_2^{\pm 1} \, $, 
$ \, \cdots \, $,  $ \, \Lambda_{n-1}^{\pm 1} \, $,  $ \,
\Lambda_n^{\pm 1} \in \gerfqm\big[\phi^{-1}\big] \, $.  Finally,
using this and the  formul{\ae}  for the  $ t_{i,j} $'s  in
Proposition 3.3 we eventually find also  $ \; \ebar_{i,j} \, $, 
$ \, \fbar_{i,j} \in \gerfqm\big[\phi^{-1}\big] \; $  for all 
$ \, 1 \leq i < j \leq n \, $.  All this says that  $ \gerfqm
\big[\phi^{-1}\big] $  {\sl contains all generators of}  $ \,
\calH_q^g \, $,  so (as  $ \, \gerfqm \big[\phi^{-1}\big] \subseteq
\calH_q^g \, $  by construction) we conclude that  $ \, \gerfqm
\big[\phi^{-1}\big] = \, \calH_q^g \, $,  \, q.e.d.  By scalar
extension it follows that  $ \, \fqm\big[\phi^{-1}\big] = \,
\H_q^g \, $  too.   
                                                 \par   
   Finally, recall that  $ \, \calfqm = \fqm \cap \, \gerH_q^g \; $. 
Then given  $ \, f \in \gerH_q^g \subset \H_q^g = \fqm \big[
\phi^{-1} \big] \, $,  \, there is  $ \, k \in \N \, $  such that  $ \,
f \phi^{\,k} \in \fqm \cap \phi^{\,k} \calfqm \subseteq \fqm \cap
\gerH_q^g = \calfqm \, $,  \, whence  $ \, f \in \calfqm \big[
\phi^{-1} \big] \, $.  Thus  $ \, \gerH_q^g \subseteq \calfqm
\big[\phi^{-1}\big] \, $,  \, and the converse is clear.  This
proves  {\it (a)}.   
                                                 \par   
   {\it (b)} \, By construction we have  $ \; \fqm \longhookrightarrow
\fqgl \longhookrightarrow \H_q^g \; $,  \; which clearly induces  $ \,
\fqm\big[\phi^{-1}\big] \longhookrightarrow \fqgl\big[\phi^{-1}\big]
\longhookrightarrow \H_q^g \; $.  But then by  {\it (a)\/}  we get
also  $ \, \fqgl\big[\phi^{-1}\big] = \H_q^g \, $.  The same
argument works for integral forms as well.  
                                                 \par   
   {\it (c)} \, By construction we have  $ \,\H_q^s = \pi\big(\H_q^g\big)
\, $  and  $ \, \fqsl = \pi\big(\fqm\big) \, $,  \, where in both cases 
$ \pi $  denotes the natural epimorphism of \S 3.2.  Now, by Corollary
3.8 below   --- whose proof needs only claim  {\it (a)\/}  above, so
we are  {\sl not\/}  using a circular argument! ---  we have  
 \vskip-5pt  
  $$  \fqsl\big[\phi^{-1}\big] = \, \pi\big(\fqm\big)
\big[\pi(\phi)^{-1}\big] = \, \pi\Big(\fqm\big[\phi^{-1}\big]\Big)
= \, \pi\big(\H_q^g\big) = \, \H_q^s \; .  $$   
 \vskip-3pt  
The same argument works for integral forms too, thus proving  {\it (c)}.  
\hskip87pt \hfill \qed\break   
\enddemo   

\vskip-5pt   

   As a byproduct of  Theorem 3.7{\it (a)},  we have the following
consequence:   
%
%

\vskip9pt   

\proclaim{Corollary 3.8} \, Let's identify  $ \fqm $  with its
copy  $ \, \xi\big(\fqm\big) $  inside  $ \uqgls \, $.  Then   
 \vskip-7pt  
  $$  D_q \, = \, \Lambda_1 \Lambda_2 \cdots \Lambda_{n-1}
\Lambda_n \, =: \, T_1  $$  
 \vskip1pt  
   \indent   Therefore, the various embeddings  $ \xi $  of quantum
function algebras into dual quantum groups, and their restrictions
to integral forms (cf.~\S 3.2), are compatible with the canonical
epimorphisms  $ \pi \, $,  \, and their restrictions to integral
forms (cf.~\S 3.2), namely  $ \, \pi \circ \xi = \xi \circ \pi \, $. 
In other words,  $ \, \fqgl \, {\buildrel \pi \over
{\relbar\joinrel\relbar\joinrel\twoheadrightarrow}}
\, \fqsl \, $  and  $ \, \fqm \, {\buildrel \pi \over
{\relbar\joinrel\relbar\joinrel\twoheadrightarrow}}
\, \fqsl \, $  are the restrictions of  $ \, \uqgls {\buildrel
\pi \over {\relbar\joinrel\relbar\joinrel\twoheadrightarrow}\,}
\uqsls \, $,  \, and similarly for their restrictions to integral
forms.   
\endproclaim   

\demo{Proof} Let  $ \, Z(A) \, $  denote the centre of any algebra 
$ A \, $.  We know that  $ \, D_q \in Z\big(\fqm\big) \, $;  \, by 
Theorem 3.7{\it (a)},  this gives also  $ \, D_q \in Z \big( \fqm
\big[\phi^{-1}\big] \big) = Z\big(\H_q^g\big) \, $.  On the other
hand, one sees (via direct computation) that  $ \, Z\big(\H_q^g\big)
= \Qq \big[\, T_1, T_1^{\,-1} \big] \, $,  \, where  $ \, T_1 :=
\Lambda_1 \Lambda_2 \cdots \Lambda_{n-1} \Lambda_n \, $.  Therefore 
$ \, D_q = P\big(T_1,T_1^{\,-1}\big) \, $,  \, for some Laurent
polynomial  $ P \, $. 
                                                   \par   
   Now, since  $ D_q $  and  $ T_1^{\pm 1} $  are both group-like,
it must necessarily be  $ \, D_q = T_1^{\,z} \, $,  \, for some 
$ \, z \in \Z \, $.  Finally, the pairing with  $ \uqgl $  gives 
$ \, \big\langle D_q \, , G_i \big\rangle = q \, $  and  $ \,
\big\langle T_1^{\,z} \, , G_i \big\rangle = q^z \, $  (for all 
$ \, i = 1, \dots, n \, $)  for any  $ \, z \in \Z \, $.  This
forces  $ \, z = 1 \, $,  \, hence  $ \, D_q = T_1 \; $.  
\hskip121pt \hfill \qed   
\enddemo   
%
%
 \eject   

\centerline {\bf \S\; 4 \ The structure of  $ \, \calfqm \, $,  $ \,
\calfqgl \, $  and  $ \, \calfqsl \, $,}   
\centerline {\bf specializations and quantum Frobenius epimorphisms. }

\vskip13pt

  {\bf 4.1 Summary.} \, In this section we present our main results. 
First, we prove that  $ \calfqm $  is a free  $ \Zqqm $--module, 
providing a basis of Poincar\'e-Birkhoff-Witt type; we give a
presentation by generators and relations, and we show that  $ \calfqm $ 
is a  $ \Zqqm $--subbialgebra,  hence a  $ \Zqqm $--integral  form, of 
$ \fqm \, $.  Also, we prove that the specialization of  $ \calfqm $ 
at  $ \, q = 1 \, $  is just  $ U_\Z\big({\gergl_n}^{\!\!*} \big) $; 
in particular, it is a  {\sl Hopf\/}  $ \Z $--algebra,  which is true
also for any other specialization at a root of 1 but needs other
arguments to be proved (see Corollary 5.7).  Finally, we show that
the quantum Frobenius morphism (2.1) for  $ {\gergl_n}^{\!\!*} $ 
is defined over  $ \Zeps \, $,  and we describe it in terms of
the previously mentioned presentation.  

  All this has direct consequences regarding  $ \calfqgl $  and 
$ \calfqsl \, $.  Namely, we provide a  $ \Zqqm $--spanning  set
for  $ \calfqgl $,  we give a presentation by generators and relations,
and we show that  $ \calfqgl $  is a  $ \Zqqm $--integral  form, as a
Hopf algebra, of  $ \fqgl \, $.  As a consequence, we prove that the
specialization of  $ \calfqgl $  at  $ \, q = 1 \, $  is  $ U_\Z
\big({\gergl_n}^{\!*} \big) $.  Another result is that for any
root of unity  $ \varepsilon $  of odd order we have an isomorphism 
$ \, \calfegl \cong \calfem \, $  as  $ \Zeps $--bialgebras: 
in particular,  $ \calfem $  is a Hopf  $ \Zeps $--algebra. 
Finally, we show that the quantum Frobenius morphism (2.2) for 
$ {\gergl_n}^{\!*} $  is defined over  $ \Zeps \, $,  and we
describe it in terms of the above mentioned presentation.  The same
results also are proved, up to details, for  $ \calfqsl $  too.   
                                           \par
   The results about  $ \calfqgl $  are proved from those about 
$ \calfqm $,  via a monomorphism  $ \; \calfqm \longhookrightarrow
\calfqgl \, $  (over  $ \Zqqm \, $)  induced by the monomorphism  $ \;
\fqm \longhookrightarrow \fqgl \, $  (over  $ \Qq \, $).  Instead, all
results about  $ \calfqsl $  come out of those for  $ \calfqm $  via
the epimorphism  $ \; \calfqm \longtwoheadrightarrow \calfqm \Big/
\big(D_q - 1\big) \, \cong \, \calfqsl \, $  induced by the
natural map (a  $ \Zqqm $--bialgebra  epimorphism!)  $ \;
\fqm \longtwoheadrightarrow \fqm \Big/ \big(D_q - 1\big)
\, \cong \, \fqsl \; $.  
 \vskip1pt   
   We need more notation.  Set  $ \; \t_{i,j} := \big( q - q^{-1}
\big)^{-1} t_{i,j} \; $  for all  $ \, i \not= j \, $,  \, and 
$ \; \t_{i,j}^{\,(m)} \! := \t_{i,j}^{\; m} \! \Big/ [m]_q! \; $ 
for all  $ \, m \in \N \, $,  like in \S 1.3.  Next result will
be basic in the following:     

\vskip11pt

\proclaim{Lemma 4.2}  For all  $ \, m, k \in \N \, $  and  $ \,
c \in \Z \, $,  \, the embedding  $ \; \xi \, \colon \, \fqm
\longhookrightarrow \uqgls \; $  gives   
  $$  \eqalign{ 
   \xi\Big(\t_{i,j}^{\,(m)}\Big) \;  =  &  \; {(-1)}^{m(j-i+1)} \,
q^{m(j-i-2) - {m \choose 2}} {\textstyle \sum\limits_{\sum_{k=0}^{n-j}
e_k = m}}  \hskip-7pt  q^{e_0 + \sum_{s=1}^{n-j} {e_s \choose 2}} \,
{(-1)}^{-e_0} \, {\big( q \! - \! q^{-1} \big)}^{m - e_0} \times  \cr  
   \times \, {[e_1]}_q!  &  \cdots {[e_{n-j}]}_q! \cdot
E_{i,j}^{\,(e_0)} E_{i,j+1}^{\,(e_1)} \cdots E_{i,n}^{\,(e_{n-j})}
\, \Lambda_j^{e_0} \Lambda_{j+1}^{e_1} \cdots \Lambda_n^{e_{n-j}}
\, F_{j+1,j}^{\,(e_1)} \cdots F_{n,j}^{\,(e_{n-j})}  \hskip2pt 
\big( i \! < \! j \big)  \cr   
   \xi\bigg(\! \bigg( {{t_{h,h} \, ; c} \atop k}  &  \bigg) \!\bigg) \;
= \; {\textstyle \sum\limits_{r=0}^k} \, q^{r(c-k-2) - {r \choose 2}}
\, {\big( q - q^{-1} \big)}^r {\textstyle \sum\limits_{\sum_{s=1}^{n-h}
e_s = r}}  \hskip-9pt  q^{\sum_{s=1}^{n-h} {e_s \choose 2}} \,
{[e_1]}_q! \cdots {[e_{n-h}]}_q! \, \times   \hfill  \cr   
   &  \quad  \hfill   \times E_{h,h+1}^{\,(e_1)} \cdots E_{h,n}^{\,
(e_{n-h})} \, \left\{ {{\Lambda_h \, ; c-r} \atop {k \, , r}} \right\} \,
\Lambda_{h+1}^{\,e_1} \cdots \Lambda_n^{\,e_{n-h}} \, F_{h+1,h}^{\,(e_1)}
\cdots F_{n,h}^{\,(e_{n-h})}  \hskip4pt 
\big(\, \forall \, h \big)  \cr
   \xi\Big(\t_{i,j}^{\,(m)}\Big) \; =  &  \; {(-1)}^{m(j-i+1)} \,
q^{m(j-i-2) - {m \choose 2}} {\textstyle \sum\limits_{\sum_{k=0}^{n-i}
e_k = m}}  \hskip-7pt  q^{m e_0 + \sum_{s=1}^{n-i} {e_s \choose 2}}
\, {\big( q - q^{-1} \big)}^{m - e_0} \times  \cr  
   \times \, {[e_1]}_q!  &  \cdots {[e_{n-i}]}_q! \cdot E_{i,i+1}^{\,(e_1)}
\cdots E_{i,n}^{\,(e_{n-i})} \, \Lambda_i^{e_0} \Lambda_{i+1}^{e_1}
\cdots \Lambda_n^{e_{n-i}} \, F_{i,j}^{\,(e_0)} F_{i+1,j}^{\,(e_1)}
\cdots F_{n,j}^{\,(e_{n-i})}  \hskip7pt  \big( i \! > \! j \big) 
\cr }  $$   
\endproclaim   

\demo{Proof}  Let's look at  $ \, \xi\Big(\t_{i,j}^{\,(m)}\Big) \, $ 
for  $ \, i < j \, $.  By Proposition 3.3 we have   
  $$  \xi\Big(\t_{i,j}^{\,(m)}\Big) \; = \; {\left( -{(-q)}^{j-i-1}
\, E_{i,j} \, \Lambda_j \, - \, {(-q)}^{j-i-2} \, {\textstyle
\sum_{k=j+1}^n} \big( q - q^{-1} \big) \, E_{i,k} \, \Lambda_k
\, F_{k,j} \right)}^{(m)} \; .  $$
 \eject   
\noindent   
 The right-hand side above fulfills the hypotheses of  Lemma
1.5{\it (b)\/}  with  $ \, p = q \, $,  \, and the claim follows
at once by brute computation.  The case  $ \, i > j \, $  is
entirely similar.   
                                       \par   
   As for  $ \, \displaystyle{ \xi\left(\! {{t_{h,h} \, ; c} \choose k}
\!\right) } $,  \, Proposition 3.3 again gives  
  $$  \xi\left(\! {{t_{h,h} \, ; c} \choose k} \!\right) \; = \;
{{\xi(t_{h,h}) \, ; c} \choose k} \; = \; \Bigg( {{\Lambda_h \, - \,
{\big( q - q^{-1} \big)}^2 \, q^{-2} \, {\textstyle \sum_{s=h+1}^n}
\, E_{h,s} \, \Lambda_s \, F_{s,h} \;\; ; \;\, c} \atop k} \Bigg)  $$   
Applying  Lemma 1.6{\it (a)},  with  $ \; x := \Lambda_h \; $  and 
$ \; w := - q^{-2} \, {\textstyle \sum_{s=h+1}^n} \, E_{h,s} \,
\Lambda_s \, F_{s,h} \, $,  \; we get   
  $$  \xi\left(\! {{t_{h,h} \, ; c} \choose k} \!\right) \; = \;
{\textstyle \sum\limits_{r=0}^k} \; q^{r(c-k)} {\big( q - q^{-1}
\big)}^r \, q^{-2r} \cdot {\Big( {\textstyle \sum_{s=h+1}^n} \, E_{h,s}
\, \Lambda_s \, F_{s,h} \Big)}^{(r)} \, \left\{ {{\Lambda_h \, ; c}
\atop {k \, , r}} \right\}  $$   
Now we apply  Lemma 1.5{\it (b)\/}  to expand  $ \, {\Big(
{\textstyle \sum_{s=h+1}^n} \, E_{h,s} \, \Lambda_s \, F_{s,h}
\Big)}^{(r)} $.  Eventually, one rearranges factors using the
commutation relations in  $ \uqgls $,  and the result is there.  
\hskip61pt \hfill \qed\break   
\enddemo   

\vskip-19pt

\proclaim{Theorem 4.3}
 \vskip3pt  
   (a) \,  $ \calfqm $  is a free  $ \, \Zqqm $--module,  with basis the
set of ordered (PBW-like) monomials
 \vskip-13pt  
  $$  \Cal{B}_{M_n} := \, \bigg\{\, \Cal{M}_{\underline{\tau}} :=
{\textstyle \prod\limits_{i<j}} \, \t_{i,j}^{\,(\tau_{i,j})} \,
{\textstyle \prod\limits_{k=1}^n} {{t_{k,k} \, ; 0} \choose \tau_{k,k}}
\, {\textstyle \prod\limits_{i>j}} \, \t_{i,j}^{\,(\tau_{i,j})}
\;\;\bigg|\;\; \underline{\tau} = \! {\big( \tau_{i,j} \big)}_{i,j=1}^n
\in M_n(\N) \,\bigg\}  $$
 \vskip-5pt  
\noindent   
 where monomials are ordered w.r.t.~any total order of the pairs 
$ (i,j) $  with  $ \, i \not= j \, $.  Similarly, any other
set obtained from  $ \Cal{B}_{M_n} $  via permutations of
factors is a  $ \, \Zqqm $--basis  as well.   
 \vskip3pt  
   (b) \,  $ \calfqgl $  is the  $ \, \Zqqm $--span  of
the set of ordered (PBW-like) monomials   
 \vskip-11pt   
  $$  \Cal{S}_{{GL}_n} := \, \bigg\{\, {\textstyle \prod\limits_{i<j}}
\, \t_{i,j}^{\,(\tau_{i,j})} \, {\textstyle \prod\limits_{k=1}^n}
{{t_{k,k} \, ; 0} \choose \tau_{k,k}} \, {\textstyle \prod\limits_{i>j}}
\, \t_{i,j}^{\,(\tau_{i,j})} \, {D_q}^{\!-\delta} \;\;\bigg|\;\;
\underline{\tau} = \! {\big( \tau_{i,j} \big)}_{i,j=1}^n \in
M_n(\N) \, , \, \delta \in \N \,\bigg\} \, .  $$   
 \vskip-5pt   
\noindent   
 where monomials are ordered w.r.t.~any total order of the pairs  $ (i,j) $ 
with  $ \, i \not= j \, $.  Similarly, any other set obtained from 
$ \Cal{S}_{{GL}_n} $  via permutations of factors (of the monomials
in  $ \Cal{S}_{{GL}_n} $)  is a  $ \, \Zqqm $--spanning set for  $ \,
\calfqgl \, $  as well.  Moreover, if  $ \, f \in \calfqgl \, $  then 
$ f $  can be expanded into a  $ \Zqqm $--linear  combination of elements
of  $ \Cal{S}_{{GL}_n} $  which all bear the same exponent  $ \delta \, $; 
\, similarly for the other spanning sets mentioned above.   
 \vskip3pt  
   (c) \,  $ \calfqsl $  is the  $ \, \Zqqm $--span  of the
set of ordered (PBW-like) monomials
 \vskip-3pt   
  $$  \Cal{S}_{{SL}_n} := \, \bigg\{\, {\textstyle \prod\limits_{i<j}}
\, \t_{i,j}^{\,(\tau_{i,j})} \, {\textstyle \prod\limits_{k=1}^n}
{{t_{k,k} \, ; 0} \choose \tau_{k,k}} \, {\textstyle \prod\limits_{i>j}}
\, \t_{i,j}^{\,(\tau_{i,j})} \;\;\bigg|\;\; \underline{\tau} = \! {\big(
\tau_{i,j} \big)}_{i,j=1}^n \in M_n(\N) \,\bigg\} \, .  $$  
 \vskip-1pt   
\noindent   
 where monomials are ordered w.r.t.~any total order of the pairs  $ (i,j) $ 
with  $ \, i \not= j \, $.  Similarly, any other set obtained from 
$ \Cal{S}_{{SL}_n} $  via permutations of factors (of the monomials
in  $ \Cal{S}_{{SL}_n} $)  is a  $ \, \Zqqm $--spanning set for 
$ \, \calfqsl \, $  as well.   
\endproclaim  

\demo{Proof} {\it (a)} \, We know from \S 2.4 that  $ B_{M_n} $  is a  $ \Qq $--basis 
of  $ \fqm \, $.  Then it is clear that the same is true for  $ \Cal{B}_{M_n}
\, $,  \, hence its  $ \Zqqm $--span  is a free  $ \Zqqm $--submodule  of 
$ \fqm \, $.   
                                                 \par   
   Thanks to Lemma 4.2, each of the factors  $ \t_{i,j}^{\,(\tau_{i,j})} $ 
or  $ {\textstyle \prod_{k=1}^n} \Big(\! {{t_{k,k} \, ; \, 0} \atop
\tau_k} \!\Big) $  of the monomials in  $ \Cal{B}_{M_n} $  is mapped
by  $ \; \xi \, \colon \, \fqm \longhookrightarrow \uqgls \, $  into 
$ \gerH_q^g \, $.  Since the latter is a subalgebra, each of the
monomials also is mapped into  $ \gerH_q^g \, $.  Thus (cf.~\S 3.2)
we deduce that  $ \Cal{B}_{M_n} $  is contained in  $ \calfqm \, $, 
\, hence the same holds for its (free)  $ \Zqqm $--span.   
%
%
 \eject   
   Consider in  $ \gerH_q^g $  the unique  $ \Zqqm $--algebra 
grading such that the generic monomial in the PBW basis  $ \frak{B}^g_* $ 
of  $ \gerH_q^g \, $,  \, say  $ \; \prod\limits_{i<j} E_{i,j}^{\,
(\eta_{i,j})} \, \prod\limits_{h=1}^n \Big(\! {{\Lambda_h \, ; \; 0}
\atop \lambda_h} \Big) \, \Lambda_h^{-\text{\it Ent}(\lambda_h/2)}
\, \prod\limits_{i<j} F_{j,i}^{\,(\varphi_{j,i})} \; $  (cf.~\S 2.5),
has degree  $ \; {\textstyle \sum\limits_{i<j}} \, \big( \eta_{i,j}
+ \varphi_{i,j} \big) \; $.  
%
%
Then Lemma 4.2 reads   
  $$  \eqalignno{ 
   \xi\Big(\t_{i,j}^{\,(m)}\Big) \; = \;  &  {(-1)}^{m(j-i+1)} \,
q^{m(j-i-2) - {m \choose 2}} q^m \, {(-1)}^{-m} \cdot E_{i,j}^{\,(m)}
\, \Lambda_j^m \; + \; \text{\it h.d.t.}  \hskip35pt   &  
\forall \;\; i < j  \cr   
   \xi\bigg(\! \bigg( {{t_{h,h} \, ; 0} \atop k} \bigg) \!\bigg)
\;  &  = \; \left\{ {{\Lambda_h \, ; 0} \atop {k \, , 0}} \right\}
+ \; \text{\it h.d.t.} \; = \; {{\Lambda_h \, ; 0} \choose k}
+ \; \text{\it h.d.t.}  \hskip35pt   &   \forall \;\; h  \cr  
   \xi\Big(\t_{i,j}^{\,(m)}\Big) \; = \;  &  {(-1)}^{m(j-i+1)} \,
q^{m(j-i-2) - {m \choose 2}} q^{m^2} \cdot \Lambda_i^m \, F_{i,j}^{\,(m)}
\; + \; \text{\it h.d.t.}  \hskip35pt   &   \forall \;\; i > j  \cr }  $$  
where  ``{\it h.d.t.\/}''  stands for  ``{\it higher degree terms\/}''. 
It follows that   
  $$  \displaylines{ 
   \quad  \xi\big(\Cal{M}_{\underline{\tau}}\big) \; = \;
{\textstyle \prod\limits_{i<j}} \, \xi\Big(\t_{i,j}^{\,(\tau_{i,j})}\Big)
\, {\textstyle \prod\limits_{k=1}^n} \xi\bigg(\! {{t_{k,k} \, ; 0}
\choose \tau_{k,k}} \!\bigg) \, {\textstyle \prod\limits_{i>j}}
\, \xi\Big(\t_{i,j}^{\,(\tau_{i,j})}\Big) \; =   \hfill  \cr   
%
%
   \hfill   = \; \varepsilon \, q^\zeta \, {\textstyle \prod\limits_{i<j}}
E_{i,j}^{\,(\tau_{i,j})} \, {\textstyle \prod\limits_{i<j}}
\Lambda_j^{\tau_{i,j}} \, {\textstyle \prod\limits_{k=1}^n}
{{\Lambda_k \, ; 0} \choose \tau_k} \, {\textstyle \prod\limits_{i<j}}
\Lambda_i^{\tau_{i,j}} \, {\textstyle \prod\limits_{i<j}}
F_{i,j}^{\,(\tau_{i,j})} \; + \; \text{\it h.d.t.}  \quad  \cr }  $$   
for some sign  $ \, \varepsilon = \pm 1 \, $  and some integers 
$ \, z, \zeta \in \Z \, $.  
                                             \par
   Now pick  $ \, f \in \calfqm \, $.  Since  $ \Cal{B}_{M_n} $
is a  $ \Qq $--basis  of  $ \fqm \, $,  there exists a unique
expansion  $ \; f = \sum_{\underline{\tau} \in M_n(\N)}
\chi_{\underline{\tau}} \, \Cal{M}_{\underline{\tau}} \; $  for some 
$ \, \chi_{\underline{\tau}} \in \Qq \, $,  \, almost all zero: our
goal is to show that indeed they belong to  $ \Zqqm \, $.  Let  $ \, \tau_0
:= \min \Big\{ \sum_{i \not= j} \tau'_{i,j} \,\big|\, \underline{\tau} \in
M_n(\N) \, , \, \chi_{\underline{\tau}} \not= 0 \Big\} \, $.  Then the
previous analysis yields   
  $$  \displaylines{ 
   \;\;  \xi(f) \, = \,   
%
%
                   {\textstyle \sum_{\underline{\tau} :
\sum_{i \not= j} \tau_{i,j} = \tau_0}} \, \chi_{\underline{\tau}}
\, \varepsilon \, q^\zeta \, {\textstyle \prod\limits_{i<j}}
E_{i,j}^{\,(\tau_{i,j})} \, {\textstyle \prod\limits_{i<j}}
\Lambda_j^{\tau_{i,j}} \, {\textstyle \prod\limits_{k=1}^n}
{{\Lambda_k \, ; 0} \choose \tau_k} \, {\textstyle \prod\limits_{i<j}}
\Lambda_i^{\tau_{i,j}} \, {\textstyle \prod\limits_{i<j}}
F_{i,j}^{\,(\tau_{i,j})} \; + \; \text{\it h.d.t.}  \;  \cr }  $$   
Pick now  $ \, \prod\limits_{i<j} \! \ebar_{i,j}^{\,\tau'_{i,j}}
\prod\limits_{k=1}^n \! G_k^{\!\phantom{|}\tau'_{k,k}}
\prod\limits_{i<j} \! \fbar_{j,i}^{\,\tau'_{j,i}} \, $,  \,
a monomial in the PBW-like basis  $ \Cal{B}^g $  of  $ \caluqgl $ 
(see \S 2.3).  Then the above expansion of  $ \xi(f) $  along with
formula (3.1) yields  
  $$  \displaylines{ 
   \bigg\langle f \, , \, {\textstyle \prod\limits_{i<j}}
\ebar_{i,j}^{\,\tau'_{i,j}} {\textstyle \prod\limits_{k=1}^n}
G_k^{\!\phantom{|}\tau'_{k,k}} {\textstyle \prod\limits_{i<j}}
\fbar_{j,i}^{\,\tau'_{j,i}} \bigg\rangle \; = \;
\bigg\langle \xi(f) \, , \, {\textstyle \prod\limits_{i<j}}
\ebar_{i,j}^{\,\tau'_{i,j}} {\textstyle \prod\limits_{k=1}^n}
G_k^{\!\phantom{|}\tau'_{k,k}} {\textstyle \prod\limits_{i<j}}
\fbar_{j,i}^{\,\tau'_{j,i}} \bigg\rangle \; =   \hfill  \cr   
   = \hskip-5pt \sum_{\Sb  \tau_{i,j} = \tau'_{i,j}  \\   
                \forall \; i \not= j  \endSb} \hskip-5pt  
\chi_{\underline{\tau}} \, \varepsilon \, q^\zeta \,
\bigg\langle \, {\textstyle \prod\limits_{i<j}} E_{i,j}^{(\tau'_{i,j})}
{\textstyle \prod\limits_{i<j}} \Lambda_j^{\tau'_{i,j}}
{\textstyle \prod\limits_{k=1}^n} \! {{\Lambda_k \, ; 0}
\choose \tau_{k,k}} {\textstyle \prod\limits_{i<j}}
\Lambda_i^{\tau'_{i,j}} {\textstyle \prod\limits_{i<j}}
F_{i,j}^{(\tau'_{i,j})} , \, {\textstyle \prod\limits_{i<j}} 
\ebar_{i,j}^{\,\tau'_{i,j}} {\textstyle \prod\limits_{k=1}^n}
G_k^{\!\phantom{|}\tau'_{k,k}} {\textstyle \prod\limits_{i<j}}
\fbar_{j,i}^{\,\tau'_{j,i}} \!\bigg\rangle \, =  \cr }  $$   
 \vskip-29pt   
  $$  \hskip15pt   = \, {\textstyle \sum\limits_{\tau_{i,j} = \tau'_{i,j}
\, ,  \; \forall \, i \not= j}} \hskip5pt  \chi_{\underline{\tau}} \,
\eta \, q^\beta \cdot {\textstyle \prod\limits_{r,s=1}^n} \bigg\langle
{{\Lambda_r \, ; 0} \choose \tau_{r,r}} , \, G_k^{\!\phantom{|}
\tau'_{s,s}} \bigg\rangle \, = \, {\textstyle \sum_{\tau_{i,j}
= \tau'_{i,j} \, , \; \forall \, i \not= j}}  \hskip5pt 
\chi_{\underline{\tau}} \, \eta \, q^\beta \cdot {\textstyle
\prod\limits_{k=1}^n} {\tau'_{k,k} \choose \tau_{k,k}}  $$   
for some  $ \, \eta \in \{+1,-1\} \, $  and  $ \, \beta \in \Z \, $. 
Since  $ \, f \in \fqm \, $,  \, the last term belongs to  $ \Zqqm \, $.  
                                        \par   
  Now choose any  $ \underline{\tau}' $  inside  $ \, \Big\{ \underline{\tau}
\in \! M_n(\N) \,\big|\, \sum_{i \not= j} \tau_{i,j} = \tau_0 \, , \, \chi_{\underline{\tau}} \not= 0 \Big\} \, $  which is minimal (for the restriction of the product order in  $ \, M_n(\N) \cong \N^{n^2} \, $).  Then  $ \, \prod_{k=1}^n \Big({\tau'_{k,k} \atop \tau_{k,k}}\Big) = \delta_{\underline{\tau}',\underline{\tau}} \, $,  \, and so  $ \, \chi_{\underline{\tau}'} \, \eta \, q^\beta = \Big\langle f \, , \, \prod_{i<j} \ebar_{i,j}^{\,\tau'_{i,j}} \prod_{k=1}^n G_k^{\!\phantom{|}\tau'_{k,k}} \prod_{i<j} \fbar_{j,i}^{\,\tau'_{j,i}}
\Big\rangle \in \Zqqm \, $,  \, hence  $ \, \chi_{\underline{\tau}'}
\in \Zqqm \, $.   
                                      \par
   Thus we get that  $ \, \chi_{\underline{\tau}'} \in \Zqqm \, $  for
the  $ \, \underline{\tau}' \in M_n(\N) \, $  specified above.  Since we already proved that  $ \, \Cal{B}_{M_n} \subseteq \calfqm \, $, 
\, it follows that  $ \, f' := f - \chi_{\underline{\tau}'} \, \Cal{M}_{\underline{\tau}'} \in \calfqm \, $.  Now, by construction
the expansion of  $ f' $  as a  $ \Qq $--linear  combination of elements
of  $ \Cal{B}_{M_n} $  has one non-trivial summand less than  $ f \, $:  \, then we can apply again the same argument, and iterate till we find that all coefficients  $ \chi_{\underline\tau} $  in the original expansion of  $ f $  do belong to  $ \Zqqm \, $.    
                                      \par
   Finally, the last observation about other bases is clear.   
 \vskip3pt   
   {\it (b)} \, By claim  {\it (a)},  every monomial of type  $ \;
\Cal{M}_{\underline{\tau}} := {\textstyle \prod\limits_{i<j}}
\, \t_{i,j}^{\,(\tau_{i,j})} \, {\textstyle \prod\limits_{k=1}^n}
\Big(\! {{t_{k,k} \, ; \, 0} \atop \tau_{k,k}} \!\Big) \,
{\textstyle \prod\limits_{i>j}} \, \t_{i,j}^{\,(\tau_{i,j})}
\; $  is  $ \Zqqm $--valued  on  $ \caluqgl \, $,  i.e., it takes
values in  $ \Zqqm $  when paired with  $ \caluqgl \, $.  On the
other hand, the same is true for  $ \, {D_q}^{\!-\delta} \, $  ($ \,
\forall \; \delta \in \N \, $),  \, because  $ \; {D_q}^{\!-\delta}
\in \gerfqgl \, $  and  $ \, \gerfqgl \subseteq \calfqgl \, $  since 
$ \, \geruqgl \supseteq \caluqgl \, $.  Finally,  
  $$  \left\langle \Cal{M}_{\underline{\tau}} \, {D_q}^{\!\! -\delta} ,
\, \caluqgl \! \right\rangle = \left\langle \Cal{M}_{\underline{\tau}}
\otimes {D_q}^{\!\! -\delta} , \, \Delta\big(\caluqgl\big) \!
\right\rangle \subseteq \left\langle \Cal{M}_{\underline{\tau}} ,
\, \caluqgl \right\rangle \cdot \left\langle {D_q}^{\!\! -\delta},
\, \caluqgl \right\rangle  $$   
and the last product on right-hand side belongs to  $ \Zqqm \, $. 
Thus  $ \; \Cal{M}_{\underline{\tau}} \, {D_q}^{\!\! -\delta} \,
\in \, \calfqgl \, $,  \, for all  $ \, \underline{\tau} \in
M_n(\N) \, $,  $ \, \delta \in \N \, $,  and so the 
$ \Zqqm $--span  of  $ \Cal{S}_{{GL}_n} $  is contained
in  $ \calfqgl \, $.   
                                       \par
   Conversely, let  $ \, f \in \calfqgl \, $.  Then there exists  $ \,
\delta \in \N \, $  such that  $ \, f \, {D_q}^{\!\delta} \in \fqm \, $. 
In addition,  $ \; \left\langle f \, {D_q}^{\!\delta}, \, \caluqgl
\right\rangle \, = \, \left\langle f \otimes {D_q}^{\!\delta}, \,
\Delta\big(\caluqgl\big) \right\rangle \, \subseteq \left\langle f,
\, \caluqgl \right\rangle \cdot \left\langle {D_q}^{\!\delta}, \,
\caluqgl \right\rangle \, \subseteq \, \Zqqm \; $  because  $ \, f,
{D_q}^{\!\delta} \in \calfqgl \, $.  Thus  $ \; f \, {D_q}^{\!\delta}
\in \calfqm \, $;  \; then, by claim  {\it (a)},  $ \, f \, {D_q}^{\!\delta}
\, $  belongs to the  $ \Zqqm $--span  of  $ \, \Cal{B}_{M_n} \, $, 
\, whence claim  {\it (b)\/}  follows at once.   
                                      \par
   Finally, the last observation about other spanning sets
is self-evident.   
 \vskip3pt   
   {\it (c)} \, The projection epimorphism  $ \; \fqm \,{\buildrel \pi
\over {\llongtwoheadrightarrow}}\, \fqm \Big/ \big(D_q - 1\big) \,
\cong \, \fqsl \, $  maps  $ \Cal{B}_{M_n} $  onto  $ \Cal{S}_{{SL}_n}
\, $.  It follows directly from definitions that  $ \, \pi \big( \calfqm
\big) \subseteq \calfqsl \, $,  \, hence in particular (thanks to claim 
{\it (a)\/})  the  $ \Zqqm $--span  of  $ \Cal{S}_{{SL}_n} $  is contained
in  $ \calfqsl \, $.  
                                            \par   
   Conversely, let  $ \, f \in \calfqsl \, $.  Then 
$ \, \left\langle f \, , \, {\textstyle \prod_{i<j}} \,
\ebar_{i,j}^{\,\eta_{i,j}} {\textstyle \prod_{k=1}^{n-1}} \,
K_k^{\kappa_k} {\textstyle \prod_{i<j}} \, \fbar_{j,i}^{\,
\varphi_{j,i}} \right\rangle \in \Zqqm \, $   
for every monomial  $ {\textstyle \prod_{i<j}} \,
\ebar_{i,j}^{\,\eta_{i,j}} {\textstyle \prod_{k=1}^{n-1}} \,
K_k^{\kappa_k} {\textstyle \prod_{i<j}} \, \fbar_{j,i}^{\,
\varphi_{j,i}} $  in the  $ \Zqqm $--basis  $ \Cal{B}^s $  of 
$ \caluqsl $  (cf.~\S 2.3).  Moreover,  $ \, {\textstyle
\sum_{(f)}} \, f_{(1)} \otimes f_{(2)} = \Delta(f) \in \Delta
\big(\fqsl\big) \subseteq \fqsl \otimes \fqsl \, $,  \, hence   
  $$  \displaylines{ 
   \left\langle f \, , \, L_n^\lambda \cdot {\textstyle \prod_{i<j}}
\, \ebar_{i,j}^{\,\eta_{i,j}} {\textstyle \prod_{k=1}^{n-1}} \,
K_k^{\kappa_k} {\textstyle \prod_{i<j}} \, \fbar_{j,i}^{\,
\varphi_{j,i}} \right\rangle \, = \, \bigg\langle \Delta(f)
\, , \, L_n^\lambda \otimes {\textstyle \prod\limits_{i<j}} \,
\ebar_{i,j}^{\,\eta_{i,j}} {\textstyle \prod\limits_{k=1}^{n-1}} \,
K_k^{\kappa_k} {\textstyle \prod\limits_{i<j}} \, \fbar_{j,i}^{\,
\varphi_{j,i}} \bigg\rangle \, =   \hfill  \cr   
   \hfill   = \, {\textstyle \sum_{(f)}} \Big\langle f_{(1)}
\, , \, L_n^\lambda \Big\rangle \cdot \left\langle f_{(2)} \, , \,
{\textstyle \prod_{i<j}} \, \ebar_{i,j}^{\,\eta_{i,j}} {\textstyle
\prod_{k=1}^{n-1}} \, K_k^{\kappa_k} {\textstyle \prod_{i<j}}
\, \fbar_{j,i}^{\,\varphi_{j,i}} \right\rangle  \cr }  $$   
   \indent   On the other hand, the natural embedding  $ \, \uqsl
\longhookrightarrow \uqgl \, $  has a canonical section  $ \, \uqgl
\,{\buildrel \text{\it pr} \over {\longtwoheadrightarrow}}\, \uqsl \, $ 
whose kernel is generated by  $ \big( L_n - 1 \big) $.  In fact, one
has an algebra isomorphism  $ \, \uqgl \cong \big(\Zqqm\big) \big[L_n,
L_n^{-1}\big] \otimes_{\Z\,[q,q^{-1}]} \uqsl \, $.  Dually, the
natural epimorphism  $ \, \fqgl \longtwoheadrightarrow \fqsl \, $ 
(or  $ \, \fqm \longtwoheadrightarrow \fqsl \, $  as well) has a
canonical section  $ \, \fqsl \longhookrightarrow \fqgl \, $  (or 
$ \, \fqsl \longhookrightarrow \fqm \, $  as well)  given by  $ \,
\phi \mapsto \phi \circ \text{\it pr} \, $,  \, for every  $ \, \phi
\in \fqsl \, $.  In particular this yields  $ \, \Big\langle \phi \, ,
\, L_n^{\,\lambda} \Big\rangle \equiv \Big\langle \phi \circ \text{\it
pr} \, , \, L_n^{\,\lambda} \Big\rangle = \Big\langle \phi \, , \,
\text{\it pr}\big(L_n^{\,\lambda}\big) \Big\rangle = \big\langle
\phi \, , 1 \big\rangle = \epsilon(\phi) \, $  for all  $ \, \phi
\in \fqsl \, $  and  $ \, \lambda \in \Z \, $.  This and the
previous analysis give   
  $$  \displaylines{ 
   \bigg\langle f \, , \, L_n^\lambda \cdot {\textstyle
\prod\limits_{i<j}} \, \ebar_{i,j}^{\,\eta_{i,j}} {\textstyle
\prod\limits_{k=1}^{n-1}} \, K_k^{\kappa_k} {\textstyle
\prod\limits_{i<j}} \, \fbar_{j,i}^{\,\varphi_{j,i}} \bigg\rangle
\, = \, \bigg\langle \, {\textstyle \sum\limits_{(f)}} \; \epsilon(f_{(1)})
\cdot f_{(2)} \; , \, {\textstyle \prod\limits_{i<j}} \, \ebar_{i,j}^{\,
\eta_{i,j}} {\textstyle \prod\limits_{k=1}^{n-1}} \, K_k^{\kappa_k}
{\textstyle \prod\limits_{i<j}} \, \fbar_{j,i}^{\,\varphi_{j,i}}
\bigg\rangle \, =   \hfill  \cr   
   \hfill   = \, \left\langle f \, , \, {\textstyle \prod_{i<j}}
\, \ebar_{i,j}^{\,\eta_{i,j}} {\textstyle \prod_{k=1}^{n-1}} \,
K_k^{\kappa_k} {\textstyle \prod_{i<j}} \, \fbar_{j,i}^{\,
\varphi_{j,i}} \right\rangle \, \in \, \Zqqm  \cr }  $$   
Now observe that  $ \, L_n^\lambda \cdot \prod\limits_{i<j} \,
\ebar_{i,j}^{\,\eta_{i,j}} \prod\limits_{k=1}^{n-1} \, K_k^{\kappa_k}
\prod\limits_{i<j} \, \fbar_{j,i}^{\,\varphi_{j,i}} = \prod\limits_{i<j}
\, \ebar_{i,j}^{\,\eta_{i,j}} \prod\limits_{k=1}^{n-1} \, K_k^{\kappa_k}
L_n^\lambda \prod\limits_{i<j} \, \fbar_{j,i}^{\,\varphi_{j,i}} \, $, 
\, be\-cause  $ L_n $  is central in  $ \uqgl \, $,  and the latter
monomials  (for all  $ \, \eta_{i,j} \, , \varphi_{j,i} \in \N \, $ 
and all  $ \, \kappa_k \, , \lambda \in \Z \, $)  clearly form a 
$ \Zqqm $--basis  of  $ \caluqgl \, $,  \, as  $ \Cal{B}^g $  is a
similar basis (see \S 2.2, and make an easy comparison, through
definitions).  But then we proved that any  $ f $  in  $ \calfqsl $ 
takes values in  $ \Zqqm $  onto elements of a  $ \Zqqm $--basis 
of  $ \caluqgl \, $,  hence also  $ \, f \in \calfqgl \, $  via the
embedding  $ \, \fqsl \longhookrightarrow \fqgl \, $  and  $ \, f
\in \calfqm \, $  via  $ \, \fqsl \longhookrightarrow \fqm \, $. 
By claim  {\it (a)},  the last fact implies that  $ f $  (in 
$ \fqm \, $)  is in the  $ \Zqqm $--span  of  $ \Cal{B}_{M_n}
\, $,  \, hence also that  $ \, f = \pi(f) \, $  is in the 
$ \Zqqm $--span  of  $ \pi(\Cal{B}_{M_n}) = \Cal{S}_{M_n} \, $.   
                                       \par   
   Finally, the last remark about other bases is clear.  
\hskip151pt \hfill \qed\break   
\enddemo

\vskip-5pt

   {\bf 4.4 Remarks.} \, {\it (a)} {\sl Exactly the same argument
as above}   --- just using monomials in the generators  $ t_{i,j} $ 
instead of the  $ \t_{i,j}^{(m)} $'s  and the  $ \Big(\! {{t_{k,k}
\, ; \, 0 \,} \atop \tau} \!\Big) $'s  ---   {\sl gives a new,
independent proof of the identity  $ \; \gerfqm \, = \, \left\{
f \in \fqm \;\Big|\; \big\langle f, \geruqgl \big\rangle
\subseteq \Zqqm \right\} \; $}  mentioned in \S 2.4, and,
similarly, also the identities    
$ \; \gerfqgl \, = \, \left\{ f \in \fqgl \;\Big|\; \big\langle f, \geruqgl
\big\rangle \subseteq \Zqqm \right\} \; $  (see \S 2.4) and
$ \; \gerfqsl \, = \, \left\{ f \in \fqsl \;\Big|\; \big\langle f,
\geruqsl \big\rangle \subseteq \Zqqm \right\} \, $  (see \S 2.4).   
                                      \par   
   {\it (b)} \,  Clearly the  $ \Zqqm $--spanning  set
$ \Cal{S}_{{GL}_n} $  {\sl is not\/}  a  $ \Zqqm $--basis 
of  $ \calfqgl \, $.  For instance, expanding  $ D_q $  as
a  $ \Zqqm $--linear  combination of elements of the basis 
$ \Cal{B}_{M_n} $  (of  $ \calfqm \, $),  \, the relation 
$ \, D_q \, {D_q}^{\!-1} = 1 \, $  (inside  $ \, \calfqgl
\, $)  yields a non-trivial  $ \Zqqm $--linear  dependence
relation among elements of  $ \Cal{S}_{{GL}_n} \, $.   

   {\it (c)} \,  The  $ \Zqqm $--spanning  set  $ \Cal{S}_{{SL}_n} $ 
{\sl is definitely not\/}  a  $ \Zqqm $--basis  of  $ \calfqsl \, $: 
e.g., for  $ n \! = \! 2 $  the relation  $ \; t_{1,1} \, t_{2,2} -
q \, t_{1,2} \, t_{2,1} - 1 = 0 \, \; $  
%
%
yields the   
%
%
 relation (in  $ \calF_q[SL_2] \, $)   
 \vskip-10pt   
  $$  {{t_{1,1} \, ; 0} \choose 1} + {{t_{2,2} \, ; 0} \choose 1} +
\, (q-1) {{t_{1,1} \, ; 0} \choose 1} {{t_{2,2} \, ; 0} \choose 1} -
\, (1+q) \big( q - q^{-1} \big) \; \t_{1,2}^{\,(1)} \; \t_{2,1}^{\,(1)}
\; = \; 0   $$
                                      \par   
   {\it (d)} \,  Theorem 4.3 has the following immediate consequence,
whose proof is straightforward:   

\vskip11pt

\proclaim{Corollary 4.5}  For every  $ \, X \in \{M, GL\} \, $,  \,
let  $ \big( D_q - 1 \big) $  be the two-sided ideal of  $ F_q[X_n] $ 
generated by  $ \big( D_q - 1 \big) $,  and let  $ \, \Cal{D}(X_n) :=
\big( D_q - 1 \big) \cap \calF_q[X_n] \, $,  \, a two-sided ideal
of  $ \calF_q[X_n] \, $.  
                 \hfill\break
   \indent (a) \; The epimorphism  $ \; \fqm \,{\buildrel \pi \over
{\llongtwoheadrightarrow}}\, \fqm \Big/ \big( D_q - 1 \big) \, \cong
\, \fqsl \; $  restricts to an epimorphism  $ \; \calfqm \,{\buildrel
\pi \over {\llongtwoheadrightarrow}}\, \calfqm \Big/ \Cal{D}(M_n) \,
\cong \, \calfqsl \; $  of  $ \, \Zqqm $--bialgebras.
                 \hfill\break
   \indent (b) \; The epimorphism  $ \; \fqgl \,{\buildrel \pi \over
{\llongtwoheadrightarrow}}\, \fqgl \Big/ \big( D_q - 1 \big) \, \cong
\, \fqsl \; $  restricts to an epimorphism  $ \; \calfqgl \,{\buildrel
\pi \over {\llongtwoheadrightarrow}}\, \calfqgl \Big/ \Cal{D}(GL_n)
\, \cong \, \calfqsl \; $  of Hopf  $ \, \Zqqm $--algebras.  
\hskip1pt \hfill  $ \square $   
\endproclaim

\vskip11pt

\proclaim{Theorem 4.6}   
                                      \par   
   (a) \,  $ \calfqm $  is the unital associative
$ \Zqqm $--algebra  with generators  $ \; \displaystyle{
\t_{i,j}^{\,(h)}} \, $,  $ \, \displaystyle{{{t_{\ell,\ell} \, ; c}
\choose k}} \, $,  \; for all  $ \, i \, , j \, , \ell \in \{1,\dots,n\} \, $,  $ \, i \not= j \, $,  $ \, h, k \in \N \, $,  $ \, c \in \Z \, $, 
\; and relations   
  $$  \displaylines{
 \text{\bf [qDP]}  \hskip5pt
   \hfill   \hbox{\sl relations\/ {\rm (1.1)}  for \ }  X \in
\big\{\, \t_{i,j} \,\big|\, i, j \in \{1,\dots,n\},
\, i \not= j \,\big\}   \hfill  \cr   
 \text{\bf [qBC]}  \hskip5pt
   \hfill   \hbox{\sl relations\/  {\rm (1.2)}  for \ }  X \in
\big\{\, \t_{\ell,\ell} \,\big|\, \ell = 1, \dots, n \,\big\}  
\hfill  \cr   
 \text{\bf [H-V.1]}  \hskip5pt
   \hfill   \t_{i,j}^{(h)} \, \t_{i,k}^{(f)} \,\; = \;\;
q^{h f} \, \t_{i,k}^{(f)} \, \t_{i,j}^{(h)} \;\;\; , 
\hskip25pt  \t_{j,i}^{(h)} \, \t_{k,i}^{(f)} \, = \;
q^{h f} \, \t_{k,i}^{(f)} \, \t_{j,i}^{(h)}  
\hfill  \forall \;\; i \not= j < k \not= i  \cr   
 }  $$   
  $$  \displaylines{
 \text{\bf [H-V.2]}  \hskip5pt
   \hfill   {{t_{i,i} \, ; c} \choose k} \, \t_{i,j}^{(h)} \,\; = \;\;
\t_{i,j}^{(h)} \, {{t_{i,i} \, ; c+h} \choose k} \;\; ,  \hskip11pt 
{{t_{i,i} \, ; c} \choose k} \, \t_{j,i}^{(h)} \, = \;
\t_{j,i}^{(h)} \, {{t_{i,i} \, ; c+h} \choose k}
\hfill  \quad  \forall \;\; i < j  \cr   
 \text{\bf [H-V.3]}  \hskip5pt
   \hfill   {{t_{i,i} \, ; c} \choose k} \, \t_{i,j}^{(h)} \,\; = \;\;
\t_{i,j}^{(h)} \, {{t_{i,i} \, ; c-h} \choose k} \;\; ,  \hskip11pt 
{{t_{i,i} \, ; c} \choose k} \, \t_{j,i}^{(h)} \, = \;
\t_{j,i}^{(h)} \, {{t_{i,i} \, ; c-h} \choose k}
\hfill  \quad  \forall \;\; i > j  \cr   
 \text{\bf [CD.1]}  \hskip5pt
   \hfill   \hskip29pt  {{t_{\ell,\ell} \, ; c} \choose k} \,
\t_{i,j}^{(h)} \,\; = \;\; \t_{i,j}^{(h)} \, {{t_{\ell,\ell}
\, ; c} \choose k}   \hfill  \quad  \forall \;\; i < \ell < j 
\;\;\; \text{or} \;\;\;  i > \ell > j  \cr   
 \text{\bf [CD.2]}  \hskip5pt
   \hfill   \t_{i,j}^{(h)} \, \t_{\ell,k}^{(f)} \,\; =
\;\; \t_{\ell,k}^{(f)} \, \t_{i,j}^{(h)}   \hfill  \quad 
\forall \;\; j \not= i < \ell \not= k < j  \cr   
%
 \text{\bf [D.1]}  \hskip9pt  \text{for all \ }  \; i < \ell \, , 
\; j < k \; ,  \text{\ with \ }  \;  \big|\{i,j,\ell,k\}\big| = 4 \; ,  
\hfill  \cr  
   \t_{\ell,k}^{(f)} \, \t_{i,j}^{(h)} \, = \;
{\textstyle \sum\limits_{s=0}^{h \wedge f}} \, {(-1)}^s \,
q^{{{s+1} \choose 2} - s \, (h+f-s)} \, \big( q - q^{-1} \big)^s
\, {[s]}_q! \, \cdot \, \t_{i,j}^{(h-s)} \, \t_{i,k}^{(s)} \,
\t_{\ell,j}^{(s)} \, \t_{\ell,k}^{(f-s)}  \cr   
%
 \text{\bf [D.2]}  \hskip9pt  \text{for all \ }  i < j \; , 
\text{\ with \ }  C(h,k,r,s,p) := p \, ((h+k) - (r+s)) -
{\textstyle \Big(\! {p \atop 2} \!\Big)} \; ,   \hfill  \cr  
   {{t_{j,j} \, ; k} \choose r} \, {{t_{i,i} \, ; h} \choose s}
=  \,{\textstyle \sum\limits_{p=0}^{r \wedge s}} \, {(-1)}^p
\, q^{C(h,k,r,s,p)} \, \big( q - q^{-1}
\big)^p \, {[p]}_q! \, \t_{i,j}^{(p)} \, \left\{\! {t_{i,i}
\, ; \, h \! - \! p} \atop {s \, , \, p} \!\right\} \,
\left\{\! {t_{j,j} \, ; \, k \! - \! p} \atop {r \, , \, p}
\!\right\} \, \t_{j,i}^{(p)}  \cr  
%
%
 \text{\bf [D.3$\boldkey{+}$]}  \hskip9pt  \text{for all \ } 
i < j < k \; ,  \text{\ with \ }  A(h,f,r,s) := {\textstyle
\Big(\! {{r+1} \atop 2} \!\Big)} + {\textstyle \Big(\!
{s \atop 2} \!\Big)} - r \, (h+f-r) \; ,   \hfill  \cr  
   \t_{j,k}^{(f)} \, \t_{i,j}^{(h)} \, = \;
{\textstyle \sum\limits_{r=0}^{h \wedge f}} \;
{\textstyle \sum\limits_{s=0}^r} \, {(-1)}^r \, q^{A(h,k,r,s,p)}
\, {(q-1)}^s \, {(s)}_q! \, {r \choose s}_q \, \t_{i,j}^{(h-r)} \,
\t_{i,k}^{(r)} \, \t_{j,k}^{(f-r)} \, {{t_{j,j} \, ; \, f-r \,}
\choose s}  \cr   
   \hfill   =  \! {\textstyle \sum\limits_{r=0}^{h \wedge f}} \;
{\textstyle \sum\limits_{s=0}^r} \, {(-1)}^r \, q^{A(h,f,r,s)} \,
{(q \! - \! 1)}^s \, {(s)}_q! \, {r \choose s}_q \, \t_{i,j}^{(h-r)}
\, \t_{i,k}^{(r)} \, {{t_{j,j} \, ; \, 0} \choose s} \,
\t_{j,k}^{(f-r)}   \hskip41pt  \cr   
%
%
 \text{\bf [D.3$\boldkey{-}$]}  \hskip9pt  \text{for all \ } 
i < j < \ell \; ,  \text{\ with \ }  A(h,f,r,s) := {\textstyle
\Big(\! {{r+1} \atop 2} \!\Big)} + {\textstyle \Big(\!
{s \atop 2} \!\Big)} - r \, (h+f-r) \; ,   \hfill  \cr  
   \t_{\ell,j}^{(f)} \, \t_{j,i}^{(h)} \, = \;
{\textstyle \sum\limits_{r=0}^{h \wedge f}} \;
{\textstyle \sum\limits_{s=0}^r} \, {(-1)}^r
\, q^{A(h,f,r,s)} \, {(q-1)}^s \, {(s)}_q! \,
{r \choose s}_q \, {{t_{j,j} \, ; \, h-r} \choose s} \,
\t_{j,i}^{(h-r)} \, \t_{\ell,i}^{(r)} \, \t_{\ell,j}^{(f-r)}  \cr   
   \hfill \hskip1pt   =  \! {\textstyle \sum\limits_{r=0}^{h \wedge f}}
\; {\textstyle \sum\limits_{s=0}^r} \, {(-1)}^r \, q^{A(h,f,r,s)} \,
{(q \! - \! 1)}^s \, {(s)}_q! \, {r \choose s}_q \, \t_{j,i}^{(h-r)}
\, {{t_{j,j} \, ; \, 0} \choose s} \, \t_{\ell,i}^{(r)} \,
\t_{\ell,j}^{(f-r)}   \hskip41pt  \cr   
%
%
 \text{\bf [D.4$\boldkey{+}$]}  \hskip9pt  \text{for all \ } 
\;  i < \ell < j  \;\;  \text{\ or \ } \; i < j < \ell   \hfill  \cr   
   \t_{\ell,j}^{(f)} \, {{t_{i,i} \, ; c} \choose k} \, = \;
{\textstyle \sum\limits_{s=0}^{f \wedge k}} \, q^{{{s+1} \choose 2}
- s \, (f-k+c)} \, \big( q - q^{-1} \big)^s \, {[s]}_q! \, \cdot \,
\left\{ {t_{i,i} \, ; \, c} \atop {k \, , \, s} \!\right\} \,
\t_{i,j}^{(s)} \, \t_{\ell,i}^{(s)} \, \t_{\ell,j}^{(f-s)}  \cr   
 \text{\bf [D.4$\boldkey{-}$]}  \hskip9pt  \text{for all \ } 
\;  i < h < \ell  \;\;  \text{\ or \ } \; h < i < \ell   \hfill  \cr   
   {{t_{\ell,\ell} \, ; c} \choose k} \, \t_{i,h}^{(f)} \; = \;
{\textstyle \sum\limits_{s=0}^{f \wedge k}} \, q^{{{s+1} \choose 2}
- s \, (f-k+c)} \, \big( q - q^{-1} \big)^s \, {[s]}_q! \, \cdot
\, \t_{i,h}^{(f-s)} \, \t_{i,\ell}^{(s)} \, \t_{\ell,h}^{(s)} \,
\left\{ {t_{\ell,\ell} \, ; \, c} \atop {k \, , \, s} \!\right\} 
\cr }  $$   
%
%
   Moreover, $ \calfqm $  is a  $ \Zqqm $--bialgebra,
whose bialgebra structure is given by   
  $$  \displaylines{ 
   \Delta \left( \!\! {{t_{i,i} \, ; c} \choose k} \!\! \right) \, =
\; {\textstyle \sum\limits_{r=0}^k} \; q^{r (c-k-1)} \big( q - q^{-1}
\big)^r \, {\textstyle \sum\limits_{s=0}^r} \; q^{{s+1} \choose 2}
{r \choose s}_{\!q} \; \times   \hfill  \cr   
   \times \, {\textstyle \sum\limits_{h=0}^{k-r}} \,
q^{h (c + 2 \, r + s - k - 1)} \big( q - q^{-1} \big)^h 
%
%
 {\textstyle \sum\limits_{\sum_{j<i} e_j = \, h}} \;
%
%
 {\textstyle \sum\limits_{\sum_{j>i} e_j = \, r}} \;
{\textstyle \prod\limits_{\ell = 1}^n} \hskip-3pt
{\phantom{\big|}}^{\widehat{i}} \,
[e_{\ell}]_q! \; q^{\sum_{a=1}^n \hskip-7pt
{\phantom{|}}^{\widehat{i}} \, {{e_a + 1} \choose 2} -
{{h + 1} \choose 2} - {{r + 1} \choose 2}} \, \times  \cr   
   \hfill   \times \; {\textstyle \prod\limits_{\ell = 1}^{i-1}}
\Big( \t_{i,\ell}^{(e_\ell)} \otimes \t_{\ell,i}^{(e_\ell)} \Big)
\cdot \left\{\! {\t_{i,i} \otimes \t_{i,i} \, ; \, c+r+s-h} \atop
{k-r \, , \, h} \!\right\} \cdot {\textstyle \prod\limits_{\ell =
i+1}^n} \Big( \t_{i,\ell}^{(e_\ell)} \otimes \t_{\ell,i}^{(e_\ell)}
\Big) \; =  \cr
 }  $$   
  $$  \displaylines{
   \qquad \hskip25pt   = \; {\textstyle \sum\limits_{r+s \leq k}} \;
{\textstyle \sum\limits_{h=0}^r} \; {\textstyle \sum\limits_{\sum_{j<i}
e_j = \, s}} \; {\textstyle \sum\limits_{\sum_{j>i} e_j = \, r}} \;
\big( q - q^{-1} \big)^{r+s} \; {r \choose h}_{\!q} \; {\textstyle
\prod\limits_{\ell = 1}^n} \hskip-3pt {\phantom{\big|}}^{\widehat{i}}
\, [e_{\ell}]_q! \; \times   \hfill  \cr   
   \times \; q^{{{h+1} \choose 2} - {{r+1} \choose 2} -
{{s+1} \choose 2} + (c-k-1)(r+s) + (2\,r+h) s + \sum_{a=1}^n
\hskip-7pt {\phantom{|}}^{\widehat{i}} \, {{e_a + 1} \choose 2}}
\; \times  \hskip39pt  \cr   
   \hfill   \times \; {\textstyle \prod\limits_{\ell = 1}^{i-1}}
\Big( \t_{i,\ell}^{(e_\ell)} \otimes \t_{\ell,i}^{(e_\ell)} \Big)
\cdot \left\{\! {\t_{i,i} \otimes \t_{i,i} \, ; \, c+r+h-s} \atop
{k-r \, , \, s} \!\right\} \cdot {\textstyle \prod\limits_{\ell =
i+1}^n} \Big( \t_{i,\ell}^{(e_\ell)} \otimes \t_{\ell,i}^{(e_\ell)}
\Big)  \hskip11pt   \forall \;\, i  \cr   
   \Delta \left( \t_{i,j}^{(h)} \right) \, = {\textstyle
\sum\limits_{\sum_s e_s = h}} q^{\sum_{r=1}^n {e_r \choose 2}
- {h \choose 2}} \, \big( q - q^{-1} \big)^{h - e_i - e_j}
\, {\textstyle \prod\limits_{k=1}^n} \hskip-3pt
{\phantom{\big|}}^{\widehat{i},\widehat{j}} \, {[e_k]}_{q\,}!
\, {\textstyle \sum\limits_{r=0}^{e_i}} \, {\textstyle
\sum\limits_{s=0}^{e_j}} \, q^{{r \choose 2} + {s \choose 2}}
\, {(r)}_{q\,}! \, {(s)}_{q\,}! \; \times   \hfill  \cr   
   \hfill   \times \; {e_i \choose r}_{\!q} \, {e_j \choose s}_{\!q}
\, {(q \! - \! 1)}^{r+s} \cdot {\textstyle \prod\limits_{k=1}^{i-1}}
\t_{i,k}^{(e_k)} \, {{t_{i,i} \, ; \, 0} \choose e_i} \,
{\textstyle \prod\limits_{k=i+1}^n} \t_{i,k}^{(e_k)} \, \otimes
\, {\textstyle \prod\limits_{k=1}^{j-1}} \t_{k,j}^{(e_k)} \,
{{t_{j,j} \, ; \, 0} \choose e_j} \, {\textstyle \prod\limits_{k=j+1}^n}
\! \t_{k,j}^{(e_k)}  \hskip9pt   \forall \; i \not= \! j  \cr   
   \hfill   \epsilon \left( \!\! {{\t_{\ell,\ell} \, ; c} \choose k}
\!\! \right) = {c \choose k}_{\!q} \;\; ,  \hskip7pt
\epsilon \left( \t_{i,j}^{(h)} \right) = 0   \hskip13pt  \forall \;\,
h \, , k \in \N \, ,  \; c \in \Z \, ,  \; \ell, i, j = 1, \dots, n \, , 
\; i \not= j  \cr }  $$
(notation of\/ \S 1.3) where terms like  $ \; \Big(\! {{x \otimes x
\, ; \, \sigma} \atop t} \!\Big) \; $  (with  $ \, x \in \{a,d\,\} \, $, 
$ \, \sigma \in \Z \, $  and  $ \, t \in \N_+ \, $)  read   
  $$  \displaylines{
   {{x \otimes x \, ; 2 \tau} \choose t} \, = \, \sum_{r+s=\nu} q^{-sr}
{{x \, ; \tau} \choose r} \! \otimes {{x \, ; \tau} \choose s} \, x^r
\, = \, \sum_{r+s=\nu} q^{-rs} x^s {{x \, ; \tau} \choose r} \! \otimes
{{x \, ; \tau} \choose s}  \cr
   {{x \otimes x \, ; 2 \tau \! + \! 1} \choose t} \! = \!\!\!
\sum_{r+s=\nu} \!\!\! q^{-(1-s)r} {{x \, ; \tau} \choose r} \otimes
{{x \, ; \tau \! + \! 1} \choose s} \, x^r = \!\!\! \sum_{r+s=\nu}
\!\!\! q^{-(1-r)s} x^s {{x \, ; \tau \! + \! 1} \choose r} \!
\otimes {{x \, ; \tau} \choose s}  \cr }  $$
according to whether  $ \sigma $  is even  ($ \, = 2 \tau $)  or
odd  ($ \, = 2 \tau + 1 $),  and consequently for  $ \; \left\{\!
{x \otimes x \, ; \, \sigma} \atop {t \, , \, \ell} \!\!\right\} \, $.
                                     \hfill\break
   \indent   In particular,  $ \calfqm $  is a  $ \Zqqm $--integral
form (as a bialgebra) of  $ \fqm \, $.
 \vskip4pt   
   (b) \,  $ \calfqgl $  is the unital associative  $ \Zqqm $--algebra 
with generators  $ \; \displaystyle{ \t_{i,j}^{\,(h)}} \, $,  $ \,
\displaystyle{{{t_{\ell,\ell} \, ; c} \choose k}} \, $,  $ \,
{D_q}^{\!-1} \, $,  \; for all  $ \, i \, , j \, , \ell \in
\{1,\dots,n\} \, $,  $ \, i \not= j \, $,  $ \, h, k \in \N
\, $,  $ \, c \in \Z \, $,  \; enjoying the same relations
as in claim (a) plus the additional relations
 \vskip-6pt   
  $$  \displaylines{
   \hfill   {D_q}^{\!-1} \, {{t_{\ell,\ell} \, ; c} \choose k}
\, = \, {{t_{\ell,\ell} \, ; c} \choose k} \, {D_q}^{\!-1} \; , 
\quad  {D_q}^{\!-1} \, \t_{i,j}^{\,(h)} \, = \, \t_{i,j}^{\,(h)}
\, {D_q}^{\!-1}  \hfill  \forall \;\, i, j, \ell, h, k, c \;\;
(i \not= j)  \cr
   {\textstyle \sum\limits_{\sigma \in \Cal{S}_n}} \hskip-1pt
{(-q)}^{l(\sigma)} t_{1,\sigma(1)} \, t_{2,\sigma(2)} \cdots
t_{n,\sigma(n)} \cdot {D_q}^{\!-1} \, = \, 1 \, = \, {D_q}^{\!-1}
\cdot {\textstyle \sum\limits_{\sigma \in \Cal{S}_n}} \hskip-1pt
{(-q)}^{l(\sigma)} t_{1,\sigma(1)} \, t_{2,\sigma(2)} \cdots
t_{n,\sigma(n)}  \cr }  $$
 \vskip-2pt   
\noindent   
 where any  $ \, t_{i,j} $  reads  $ \; t_{i,j} := ( q - q^{-1} \big)
\, \t_{i,j}^{\,(1)} \; $  if  $ \, i \not= j \, $  and  $ \; t_{i,j}
:= 1 + (q \! - \! 1) \displaystyle{{t_{i,i} \, ; 0} \choose 1}
\, $  if  $ \; i \! = \! j \, $.   
                                               \par   
   Moreover, $ \calfqgl $  is a Hopf  $ \, \Zqqm $--algebra,  with
Hopf algebra structure given by the  formul{\ae}  in claim (a)
for  $ \Delta $  and  $ \epsilon $  plus the  {\sl implicit} 
formul{\ae}  for the antipode
  $$  \displaylines{
   S \left( \!\! {{t_{\ell,\ell} \, ; \, c} \choose k} \!\! \right)
= \, \left( {{\text{\sl det}_q \Big(\! \big(t_{i,j}\big)_{i, j = 1,
\dots, n}^{i, j \not= \ell} \Big) \, {D_q}^{\!-1} \; ; \; c} \atop k}
\right)  \cr   
   S \left( \t_{i,j}^{\,(h)} \right) \, = \, {(-1)}^{(i+j)h} \,
{\big( q - q^{-1} \big)}^{-h} \cdot \text{\sl det}_q \! \left(\!
{\big(t_{r,s}\big)}_{r,s=1,\dots,n}^{r \not= i; \, s \not= j}
\right)^{\! h} \cdot D_q^{-h} \Big/ {[h]}_q!  \cr }  $$      
(for all  $ \, \ell, i, j = 1, \dots, n \, $,  $ \, h, k \in \N \, $, 
$ \, c \in \Z \, $,  \; with  $ \, i \not= j \, $),  \; and   
  $$  \Delta\big({D_q}^{\!-1}\big) \, = \, {D_q}^{\!-1} \otimes
{D_q}^{\!-1} \; ,  \qquad  \varepsilon\big({D_q}^{\!-1}\big) \, =
\, 1 \; ,  \qquad  S\big({D_q}^{\!-1}\big) \, = \, D_q \; .  $$   
   \indent   In particular,  $ \calfqgl $  is a  $ \Zqqm $--integral
form (as a Hopf algebra) of  $ \fqgl \, $.
%
%
 \eject  
   (c) \,  $ \calfqsl $  is generated, as a unital associative 
$ \Zqqm $--algebra,  by generators as in claim (a).  These generators
enjoy all relations in claim (a), plus some additional relations
(springing out of the relation  $ \; D_q \, = \, 1 \; $  in 
$ \fqsl \, $).   
                                      \par   
   Moreover,  $ \calfqsl $  is a Hopf  $ \, \Zqqm $--algebra,
whose Hopf algebra structure is given by the formul{\ae}  in
claim (a) for  $ \Delta $  and  $ \epsilon \, $,  \, and by the
formul{\ae}  in (b), but taking into account  $ \, D_q = 1 \, $ 
\, for the antipode.   
                                                   \par  
   In particular,  $ \calfqsl $  is a  $ \, \Zqqm $--integral
form (as a Hopf algebra) of  $ \fqsl \, $.
\endproclaim

\demo{Proof} {\it (a)} \, By Theorem 4.3, the set of elements considered
in the statement generates  $ \calfqm $  (and even more, as for binomial
coefficients we can restrict to take only those with  $ \, r = 0 = s
\, $).    
 \vskip5pt   
   As for relations, we set an alphanumerical key on their left, which
refer to the type of relation, reminding what they arise from.  
 \vskip3pt   
   {\bf [qBC]}  is a reminder for ``quantum binomial coefficients'';
similarly  {\bf [qDP]}  stand for ``quantum divided powers''.  It is
clear by definitions that all these relations do hold in  $ \calfqm
\, $.
 \vskip3pt   
   {\bf [H-V.n]}  (for  $ \, n \! = \! 1 $,  $ 2 $,  $ 3 \, $)  is a
reminder for ``horizontal-vertical'': indeed, these are the relations
induced   --- inside  $ \calfqm $  ---   by the ``horizontal'' or
``vertical'' relations holding true among matrix entries of the 
$ q $--matrix  $ \, {\big( t_{i,j} \big)}_{i,j=1,\dots,n} \, $  of
the generators of  $ \fqm $  which lie on the same row or the same
column.  Namely,  $ \; t_{i,j} \, t_{i,k} = q \, t_{i,k} \, t_{i,j} \; $ 
(for  $ \, j < k \, $)  is a ``horizontal'' relation, and  $ \; t_{i,k}
\, t_{h,k} = q \, t_{h,k} \, t_{i,k} \; $  (for  $ \, i < h \, $)  is a
``vertical'' one.  It is clear again by construction that both kind of
relations 
   \hbox{induce the corresponding relations in  $ \calfqm $  as in the
claim.}   
 \vskip3pt   
   {\bf [CD.m]}  (for  $ \, m \! = \! 1 $,  $ 2 \, $)  stands for
``counter-diagonal''.  These are the relations induced   --- inside 
$ \calfqm $  ---   by the relations in  $ \fqm $  involving a couple of
generators which, as entries of the  $ q $--matrix  of generators, are
in ``counter-diagonal'' mutual position.  Namely, we mean relations of
type  $ \; t_{i,l} \, t_{j,k} = t_{j,k} \, t_{i,l} \; $  with  $ \, i
< j \, $  and  $ \, k < l \, $,  \, that are plain commutation relations. 
Then both  {\bf [CV.1]}  and  {\bf [CV.2]}  are trivially true in 
$ \calfqm \, $.    
 \vskip3pt   
   Finally,  {\bf [D.m\,(${\boldsymbol\pm}$)]}  (for  $ \, m \! = \! 1 $, 
$ \dots $,  $ 4 \, $)  are labels for ``diagonal'' relations, i.e.~all
those relations   --- among the generators in the claim ---   which are
induced by relations in  $ \fqm $  involving a couple of entries in the 
$ q $--matrix  $ \, {\big( t_{i,j} \big)}_{i,j=1,\dots,n} \, $  which
stand in mutual diagonal position, that is relations of type  $ \;
t_{i,k} \, t_{j,l} - \, t_{j,l} \, t_{i,k} = \left( q - q^{-1} \right)
\, t_{i,l} \, t_{j,k} \; $  with  $ \, i < j \, $  and  $ \, k < l \, $. 
But now, this single type of relation in  $ \fqm $  gives rise to several
types of relations in  $ \calfqm \, $.  Indeed, once one fixes the
(``mutually diagonal'') positions, say  $ (i,k) $  and  $ (j,l) $, 
we can single out (within in the  $ q $--matrix)  the unique rectangle 
$ R $  which has one diagonal with  $ (i,k) $  and  $ (j,l) $  as
vertices.  Then different types of relations occur depending on the
position of the four corners of  $ R $  w.r.t~the main diagonal of
the  $ q $--matrix.  In the sequel, we call ``diagonal'' the corners
of  $ R $  which lie on the top-left to bottom-right diagonal of 
$ R \, $,  and ``counter-diagonal'' those on the top-right to
bottom-left diagonal of  $ R \, $.   
                                                  \par   
   Relations of type  {\bf [D.1]}  occur when no corner of  $ R $  lies
on the main diagonal.   
                                                  \par   
   Relations  {\bf [D.2]}  show up  when both diagonal corners of  $ R $ 
are on the main diagonal.   
                                                  \par   
   Relations of type  {\bf [D.3${\boldsymbol\pm}$]}  occur when one of
the counter-diagonal corners of  $ R $  lies on the main diagonal.  If
it is the bottom-left one, we have relation  {\bf [D.3$\boldkey{+}$]}, 
with the  ``$ \boldkey{+} $''  sign to remind us that the rectangle 
$ R $  lies above the main diagonal.  If instead it is the
top-right one, we have relation  {\bf [D.3$\boldkey{-}$]},  where 
``$ \boldkey{-} $''  reminds that  $ R $  lies below the diagonal.    
                                                  \par   
   Similarly, we have relations of type  {\bf [D.4${\boldsymbol\pm}$]} 
when one of the diagonal corners of  $ R $  lies on the main diagonal. 
If it is the top-left one, we have relation  {\bf [D.4$\boldkey{+}$]}, 
while if it is the bottom-right one, we have relation  {\bf
[D.4$\boldkey{-}$]}  instead.    
                                                  \par   
   In any case, these ``diagonal relations'' are the only ones which
are not immediate from definitions.  We shall now prove all of them
by suitable induction arguments.   
 \vskip4pt   
   To simplify notations, we set the following terminology for the
corners of  $ R \, $:  
 \vskip-11pt   
  $$  a \, := \, t_{i,j} \;\; ,  \quad  b \, := \, t_{i,k} \;\; , 
\quad  c \, := \, t_{\ell,j} \;\; ,  \quad  d \, := \, t_{\ell,k}
\;\; ,   \hskip21pt  \text{so that}  \hskip20pt   
   R \; = \; \hbox{$ \matrix
       a  &  \!\!\cdots\!\!  &  b  \\  
       \vdots  &  {}  &  \vdots  \\  
       c  &  \!\!\cdots\!\!  &  d  \endmatrix $} 
\hskip11pt   \eqno (4.1)  $$    
 \vskip-1pt   
\noindent   
 where some of the indices may also coincide: they do iff  $ a $  or 
$ d $  is on the main diagonal.  These generators do generate a copy
of  $ F_q[M_2] $  inside  $ F_q[M_n] \, $,  \, hence the identities
we are going to prove actually are identities in  $ F_q[M_2] \, $.   
 \vskip4pt   
   In terms of (4.1), formula  {\bf [D.1]}  reads (for all  $ \;
\big|\{i,j,\ell,k\}\big| = 4 \, $,  $ \, i < \ell \, $,  $ \, j
< k \, $)   
  $$  d^{(f)} \, a^{(h)} \, = \;
{\textstyle \sum\limits_{s=0}^{h \wedge f}} \, {(-1)}^s \,
q^{{{s+1} \choose 2} - s \, (h+f-s)} \, \big( q - q^{-1} \big)^s
\, {[s]}_q! \, \cdot \, a^{(h-s)} \, b^{(s)} \, c^{(s)}
\, d^{(f-s)} \quad ;   \eqno (4.2)  $$   
by the identities  $ \; x^{(m)} = x^m \big/ {[m]}_q! \; $  (for all 
$ m $)  this formula is equivalent, inside  $ \fqm \, $,  \, to   
  $$  d^f \, a^h \, = \; {\textstyle \sum\limits_{s=0}^{h \wedge f}} \,
{(-1)}^s \, q^{{{s+1} \choose 2} - s \, (h+f-s)} \, \big( q - q^{-1}
\big)^s \, {\textstyle \Big[ {h \atop s} \Big]}_q \, {\textstyle
\Big[ {f \atop s} \Big]}_q \, {[s]}_q! \, \cdot \, a^{h-s} \, b^s
\, c^s \, d^{f-s}   \eqno (4.3)  $$   
so we shall prove the latter.  The proof is inductive: first on  $ h $ 
for  $ \, f = 1 \, $,  \, and later on  $ \, f < h \, $.  Then one can
use the  $ \Qq $--algebra  anti-automorphism of  $ F_q[M_2] $  given
by  $ \; a \mapsto d \, $,  $ \, b \mapsto c \, $,  $ \, c \mapsto b
\, $  and  $ \, d \mapsto a \, $,  \, to handle the opposite case,
that is for  $ \, f \geq h \, $.   
                                              \par   
   For  $ \, h = f = 1 \, $  formula (4.3) is true: it is the
relation  $ \; a \, d \, - \, d \, a \, = \, \big( q - q^{-1}
\big) \, b \, c \; $.  Then set  $ \, f = 1 \, $,  \, and assume 
$ \; d \, a^h \, = \; a^h \, d \, - \, q^{1-h} \, \big( q - q^{-1}
\big) \, {\textstyle \Big[ {h \atop 1} \Big]}_q \, a^{h-1} \, b \,
c \; $,  \; by induction on  $ h \, $.  Then   
  $$  \displaylines{ 
   d \, a^{h+1} \, = \; \big( d \, a^h \big) \, a \; = \;
\Big( a^h \, d \, - \, q^{1-h} \, \big( q - q^{-1} \big) \,
{\textstyle \Big[ {h \atop 1} \Big]}_q \, a^{h-1} \, b \, c \Big)
\, a \; =   \hfill  \cr  
%
%
   = \; a^h \, \big( a \, d \, - \, \big( q - q^{-1} \big)
\, b \, c \big) \, - \, q^{1-h-2} \, \big( q - q^{-1} \big) \,
{\textstyle \Big[ {h \atop 1} \Big]}_q \, a^h \, b \, c \; =  
\hfill  \cr  
   = \; a^{h+1} \, d \, - \, \big( q - q^{-1} \big) \Big(
1 + \, q^{1-h-2} \, {\textstyle \Big[ {h \atop 1} \Big]}_q
\,\Big) \, a^h \, b \, c  \; = \;  a^{h+1} \, d \, - \,
\big( q - q^{-1} \big) \, q^{1-(h+1)} \, {\textstyle
\Big[ {{h+1} \atop 1} \Big]}_q \, a^h \, b \, c  
\hfill  \cr }  $$   
thanks to the identity  $ \; 1 + \, q^{1-h-2} \, {\textstyle \Big[
{h \atop 1} \Big]}_q  = \, q^{-h} \, {\textstyle \Big[ {{h+1}
\atop 1} \Big]}_q \; $.  This ends the proof for  $ \, f = 1 \, $.  
                                         \par   
   Now assuming that (4.3) hold for some  $ f $  with  $ \, 1 <
f < h \, $,  \, we prove it for  $ \, f+1 \, $.  Using shorthand
notation  $ \; \alpha^s_{f,h} := {{s+1} \choose 2} - s \, (h+f-s)
\, $,  \; we have   
  $$  \displaylines{ 
   d^{f+1} \, a^h \, = \; d \, \big( d^f \, a^h \big) \,
= \; d \, \bigg(\, {\textstyle \sum\limits_{s=0}^f} \, {(-1)}^s
\, q^{\alpha^s_{f,h}} \, \big( q - q^{-1} \big)^s \,
{\textstyle \Big[ {h \atop s} \Big]}_q \, {\textstyle
\Big[ {f \atop s} \Big]}_q \, {[s]}_q! \, \cdot \, a^{h-s}
\, b^s \, c^s \, d^{f-s} \bigg) \; =   \hfill  \cr  
   \hfill   = \; {\textstyle \sum\limits_{s=0}^f} \, {(-1)}^s
\, q^{\alpha^s_{f,h}} \, \big( q - q^{-1} \big)^s \,
{\textstyle \Big[ {h \atop s} \Big]}_q \, {\textstyle
\Big[ {f \atop s} \Big]}_q \, {[s]}_q! \, \cdot \, d
\, a^{h-s} \, b^s \, c^s \, d^{f-s} \; =  \cr  
%
%
   = \; {\textstyle \sum\limits_{s=0}^f} \,
{(-1)}^s \, q^{\alpha^s_{f,h} - 2 \, s} \, \big( q - q^{-1} \big)^s
\, {\textstyle \Big[ {h \atop s} \Big]}_q \, {\textstyle \Big[
{f \atop s} \Big]}_q \, {[s]}_q! \, \cdot \, a^{h-s} \, b^s
\, c^s \, d^{f+1-s} \; +   \hfill  \cr   
   \hfill   + \; {\textstyle \sum\limits_{s=0}^f} \, {(-1)}^{s+1}
\, q^{\alpha^s_{f,h} + 1 - h + s} \, \big( q - q^{-1} \big)^{s+1}
\, {\textstyle \Big[ {h \atop s} \Big]}_q \, {\textstyle \Big[
{f \atop s} \Big]}_q \, {\textstyle \Big[ {{h - s} \atop 1}
\Big]}_q \, {[s]}_q! \, \cdot \, a^{h-s-1} \, b^{s+1} \, c^{s+1}
\, d^{f-s}  \cr }  $$   
Now, in the last sum change summation index from  $ s $  to 
$ \, s + 1 \, $.  Then our last term reads  
  $$  \displaylines{ 
   q^{\alpha^0_{f,h}} \, a^h \, d^{f+1} \; + \; {(-1)}^{f+1} \,
q^{\alpha^f_{f,h} + f + 1 - h} \, {\textstyle \Big[ {h \atop f}
\Big]}_q \, {\textstyle \Big[ {f \atop f} \Big]}_q \, {\textstyle
\Big[ {{h-f} \atop 1} \Big]}_q \, {[f]}_q! \, \cdot \, a^{h-(f+1)}
\, b^{f+1} \, c^{f+1} \; +   \hfill  \cr   
    + \; {\textstyle \sum\limits_{s=1}^f} \, {(-1)}^s \, \big(
q - q^{-1} \big)^{s+1} \, \bigg( q^{\alpha^s_{f,h} - 2 \, s} \,
{\textstyle \Big[ {h \atop s} \Big]}_q \, {\textstyle \Big[
{f \atop s} \Big]}_q \, {[s]}_q! \; + \; q^{\alpha^{s-1}_{f,h}
+ s - h} \, \times   \qquad  \cr   
    \hfill   \times \, {\textstyle \Big[ {h \atop {s-1}} \Big]}_q
\, {\textstyle \Big[ {f \atop {s-1}} \Big]}_q \, {\textstyle \Big[
{{h-s+1} \atop 1} \Big]}_q \, {[s - \! 1]}_q! \bigg) \cdot \,
a^{h-s} \, b^s \, c^s \, d^{f+1-s} \; =   \quad  \cr  
   \hfill   = \; {\textstyle \sum\limits_{s=0}^{f+1}} \, {(-1)}^s
\, q^{\alpha^s_{f+1,h}} \, \big( q - q^{-1} \big)^s \, {\textstyle
\Big[ {h \atop s} \Big]}_q \, {\textstyle \Big[ {{f+1} \atop s}
\Big]}_q \, {[s]}_q! \, \cdot \, a^{h-s} \, b^s \, c^s \,
d^{(f+1)-s}  \cr }  $$   
thus ending the proof of  {\bf [D.1]},  thanks to the following
identities:   
  $$  \displaylines{ 
   \qquad \quad   \alpha^0_{f,h} \; = \; \alpha^0_{f+1,h} \;\; , 
\hskip65pt  \alpha^{s-1}_{f,h} + s - h \; = \;
{\textstyle \Big(\! {{s+1} \atop 2} \!\Big)} -
s \, (h+f+1-s) + f + 1 - s  \cr    
   \alpha^s_{f,h} - 2 s \; = \;
{\textstyle \Big(\! {{s+1} \atop 2} \!\Big)}
- s \, (h+f+1-s) - s \;\; ,  \hskip45pt 
  \alpha^s_{f+1,h} \; = \;
{\textstyle \Big(\! {{s+1} \atop 2} \!\Big)}
- s \, (h+f+1-s)  \cr    
   q^{-s} \, \big(\, {[f-s+1]}_q + \, q^{f+1} \, {[s]}_q \,\big)
\; = \; {[f+1]}_q  \cr }  $$   
   As to  {\bf [D.2]},  in terms of (4.1) it reads, with  $ \,
\b := b \Big/ \big( q - q^{-1} \big) \, $  and  $ \, \c :=
c \Big/ \big( q - q^{-1} \big) \, $,    
  $$  {{d \, ; k} \choose r} {{a \, ; h} \choose s} = \,
{\textstyle \sum\limits_{p=0}^{r \wedge s}} \, {(-1)}^p \,
q^{p \, ((h+k) - (r+s)) - {p \choose 2}} \, \big( q - q^{-1}
\big)^p \, {[p]}_q! \, \b^{(p)} \left\{\! {a \, ; \, h \! - \! p}
\atop {s \, , \, p} \!\right\} \left\{\! {d \, ; \, k \! - \! p}
\atop {r \, , \, p} \!\right\} \c^{(p)}  $$   
a formula which is proved in [GR], \S 3.4.  
 \vskip4pt   
   For  {\bf [D.3$\boldkey{+}$]},  set notation as in (4.1), now
with  $ \, \ell = j \, $.  Then (4.2) holds and, setting  $ \,
\boldkey{x} := x \Big/ \big( q - q^{-1} \big) \, $  for all 
$ \, x \in \{a,b,c,d\} \, $,  \, it yields   
  $$  \d^{(f)} \, \a^{(h)} \, = \;
{\textstyle \sum\limits_{s=0}^{h \wedge f}} \, {(-1)}^s \,
q^{{{s+1} \choose 2} - s \, (h+f-s)} \, \big( q - q^{-1} \big)^s
\, {[s]}_q! \, \cdot \, \a^{(h-s)} \, \b^{(s)} \, \c^{(s)}
\, \d^{(f-s)}  $$   
   \indent   Now, for all  $ \, m \in \N \, $,  the following formal
identity hold (cf.~[GR], Lemma 6.1):  
  $$  x^m  \; = \;  {\textstyle \sum\limits_{k=0}^m} \; q^{k \choose 2}
{m \choose k}_{\!q} \, {(q-1)}^k \, {(k)}_q! \, {{x \, ; \, 0} \choose k} 
\quad .  $$   
Using this identity to expand  $ \; \c^{(s)} = {\big( q - q^{-1}
\big)}^{-s} \, {[s]}_q^{\,-1} \cdot c^s \, $,  \; our last formula
turns into   
  $$  \displaylines{ 
   \d^{(f)} \, \a^{(a)} \; = \; {\textstyle
\sum\limits_{s=0}^{h \wedge f}} \, {(-1)}^s \,
q^{{{s+1} \choose 2} - s \, (h+f-s)} \, \a^{(h-s)}
\, \b^{(s)} \, \c^{(s)} \, \d^{(f-s)} \; =   \hfill  \cr   
%
    = \; {\textstyle \sum\limits_{s=0}^{h \wedge f}} \;
{\textstyle \sum\limits_{k=0}^s} \, {(-1)}^s \, q^{{{s+1} \choose 2} 
+ {k \choose 2} - s \, (h+f-s)} \, {\textstyle \Big( {s \atop k}
\Big)}_q \, {(q-1)}^k \, {(k)}_q! \cdot \a^{(h-s)} \, \b^{(s)} \,
{\textstyle \Big(\! {{c \, ; \, 0} \atop k} \!\Big)} \; \d^{(f-s)}
\; =  \cr
    \hfill   = \; {\textstyle \sum\limits_{s=0}^{h \wedge f}} \;
{\textstyle \sum\limits_{k=0}^s} \, {(-1)}^s \, q^{{{s+1} \choose 2} 
+ {k \choose 2} - s \, (h+f-s)} \, {\textstyle \Big( {s \atop k}
\Big)}_q \, {(q-1)}^k \, {(k)}_q! \cdot \a^{(h-s)} \, \b^{(s)}
\, \d^{(f-s)} \, {\textstyle \Big(\! {{c \, ; \, f-r} \atop k}
\!\Big)}  \cr }  $$   
where in the last step we used  {\bf [H-V.2]}.  This proves 
{\bf [D.3$\boldkey{+}$]},  and  {\bf [D.3$\boldkey{-}$]}  is
entirely similar.   
%
%
 \eject   
   Finally we prove  {\bf [D.4$\boldsymbol{\pm}$]},  starting with 
{\bf [D.4$\boldkey{+}$]}.  We use notation of (4.1), but with  $ \,
i = j \, $  and then using index  $ j $  instead of  $ k \, $;  \,
more in general, we set  $ \, a := t_{i,i} \, $,  $ \, \b^{(r)} :=
\t_{i,j}^{\,(r)} \, $,  $ \, \c^{(r)} := \t_{i,j}^{\,(r)} \, $, 
$ \, \d^{(r)} := \t_{\ell,j}^{\,(r)} \, $  (for  $ \, r \in \N
\, $).      
   \hbox{Then formula  {\bf [D.4$\boldkey{+}$]}  reads (with
index  $ u $  instead of  $ c \, $)}   
  $$  \d^{(f)} \, {{a \, ; \, u} \choose k}  \; = \;  {\textstyle
\sum\limits_{s=0}^{f \wedge \, k}} \, {(-1)}^s \, q^{{{s+1} \choose 2}
\, - \, s \, (f+k-u)} \, {[s\,]}_q! \, {\big( q - q^{-1} \big)}^s
\cdot \, \bigg\{ {{a \, ; \, u \! - \! 2 s} \atop {k \, , s}} \bigg\}
\; \b^{(s)} \, \c^{(s)} \, \d^{(f-s)}  $$   
which is equivalent to the commutation relation in  $ F_q[M_2] $   
  $$  d^f \, {{a \, ; \, u} \choose k}  \; = \;  {\textstyle
\sum\limits_{s=0}^{f \wedge \, k}} \, {(-1)}^s \, q^{{{s+1} \choose 2}
\, - \, s \, (f+k-u)} \, {\bigg[ {f \atop s} \!\,\bigg]}_q \cdot \,
\bigg\{ {{a \, ; \, u \! - \! 2 s} \atop {k \, , s}} \bigg\} \;
b^s \, c^s \, d^{f-s}   \eqno (4.4)  $$   
So we must prove this last identity: we do that by induction on  $ k $ 
and  $ f \, $.   
                                                       \par   
   For  $ \, k = 1 = f \, $,  \, formula (4.4) reads  $ \; d \,
{\displaystyle {{a \, ; u} \choose 1}} = {\displaystyle {{a \, ; u}
\choose 1}} \, d \, - \, q^{u-1} \, (q+1) \, b \, c \; $;  \, this
is equi\-valent to the identity   
    \hbox{$ \; d \, \big( a \, q^u \! - 1 \big) \, = \,
\big( a \, q^u \! - 1 \big) \, d - q^{u-1} \, \big( q^2 \!
- 1 \big) \, b \, c \; $,  \; which is easy to prove.}    
%
%
                                                       \par   
   For  $ \, k > 1 = f \, $,  \, formula (4.4) reads    
  $$  d \; {{a \, ; u} \choose k}  \; = \;  {{a \, ; u} \choose k} \;
d \, - \, q^{u-k} \, \bigg( {{a \, ; \, u \! - \! 2 \,} \choose {k-1}}
+ \, q \; {{a \, ; \, u \! - \! 1} \choose {k-1}} \bigg) \; b \, c  
\eqno (4.5)  $$   
which, using compact notation  $ \; {(a\, ; k , u)}_q := {(q \! - \! 1)}^k \, {(k)}_q! \, {\displaystyle {{a \, ; u} \choose k}} \; $,  \; is equivalent to   
  $$  d \; {(a \, ; k, u)}_q  \; = \;  {(a \, ; k, u)}_q \; d \, - \, q^{u-k} \, \big( q^k - 1 \big) \, {(a \, ; k \! - \! 2, u \! - \! 2)}_q \, \Big( \big(a \, q^{u-k} - 1 \big) + \, q \, \big(a \, q^{u-1} - 1 \big) \Big) \, b \, c  $$   
Making induction, the basis is proved above for  $ \, k \! = \! 1 \, $;  \, the inductive step from  $ \, k \, $  to  $ \, k \! + \! 1 \, $  is   
  $$  \displaylines{   
   d \; {(a \, ; k \! + \! 1, u)}_q  \; = \;  d \; \big( a \, q^{u-k} - 1 \big) \, {(a \, ; k, u)}_q  \; =   \hfill  \cr   
   \hfill   = \;  \big( a \, q^{u-k} \! - \! 1 \big) \, d \; {(a \, ; k, u)}_q \, - \, q^{u-(k+1)} \, \big( q^2 \! - \! 1 \big) \, b \, c \, {(a \, ; k, u)}_q  \; = \;  \big( a \, q^{u-k} \! - \! 1 \big) \, {(a \, ; k, u)}_q \, d \; -  \cr  
   \hfill   - \; q^{u-k} \, \big( q^k \! - \! 1 \big) \, \big( a \, q^{u-k} \! - \! 1 \big) \, {(a \, ; k \! - \! 2, u \! - \! 2)}_q \, \Big( \big( a \, q^{u-k} \! - \! 1 \big) + q \, \big( a \, q^{u-1} \! - \! 1 \big) \Big) \, b \, c \; -  \cr   
   \hfill   - \; q^{u-(k+1)} \, \big( q^2 \! - \! 1 \big) \, {(a \, ; k, u \! - \! 2)}_q \, b \, c  \; =  \cr   
%
%
    = \;  {(a \, ; k \! + \! 1, u)}_q \, d \, - \, q^{u-(k+1)} \, {(a \, ; k \! - \! 1, u \! - \! 2)}_q \, \times   \hfill  \cr   
   \hfill   \times \; \Big( q \, \big( q^k \! - \! 1 \big) \,
\Big(\! \big( a \, q^{u-k} \! - \! 1 \big) + q \, \big( a \, q^{u-1} \! - \! 1 \big) \!\Big) + \big( q^2 \! - \! 1 \big) \, \big( a \, q^{u-k-1} \! - \! 1 \big) \Big) \, b \, c \, =  \cr
   \hfill   = \,  {(a \, ; k \! + \! 1, u)}_q \; d \, - \, q^{u-(k+1)} \, \big( q^{k+1} \! - \! 1 \big) \, {(a \, ; k \! - \! 1, u \! - \! 2)}_q \, \Big(\! \big(a \, q^{u-(k+1)} \! - \! 1 \big) + \, q \, \big(a \, q^{u-1} \! - \! 1 \big) \!\Big) \; b \, c  \cr }  $$   
where the very last step follows from the (straightforward) identity   
  $$  \displaylines{ 
   q \, \big( q^k \! - \! 1 \big) \, \Big(\! \big( a \, q^{u-k}
\! - \! 1 \big) + q \, \big( a \, q^{u-1} \! - \! 1 \big) \!\Big)
+ \big( q^2 \! - \! 1 \big) \, \big( a \, q^{u-k-1} \! - \! 1 \big)
\Big)  \; =   \hfill  \cr   
   \hfill   = \;  \big( q^{k+1} \! - \! 1 \big) \, \Big(\!
\big(a \, q^{u-(k+1)} \! - \! 1 \big) + \, q \, \big(a \,
q^{u-1} \! - \! 1 \big) \!\Big)  \cr }  $$   
   \indent   To move next step, we need a more general formula than
(4.5), namely   
  $$  d \; \bigg\{ { {a \, ; u} \atop {k \, , h}} \bigg\}  \; = \; 
\bigg\{ { {a \, ; u} \atop {k \, , h}} \bigg\} \; d \, - \, q^{u-k+h}
\, \bigg\{ { {a \, ; u \! - \! 2} \atop {k \, , h \! + \! 1}} \bigg\}
\; b \, c   \eqno (4.6)  $$      
We prove it by induction on  $ h \, $.  The basis holds by (4.5), which
also gives the inductive step:    
  $$  \displaylines{ 
   d \; \bigg\{ {{a \, ; u} \atop {k \, , h}} \bigg\}  \; = \; 
   d \; {\textstyle \sum\limits_{s=0}^h} \; q^{{{s+1} \choose 2}} \,
{h \choose s}_{\!q} \, {{a \, ; \, u \! + \! s} \choose {k \! - \! h}} 
\; = \;  {\textstyle \sum\limits_{s=0}^h} \; q^{{{s+1} \choose 2}} \,
{h \choose s}_{\!q} \, d \, {{a \, ; \, u \! + \! s} \choose {k \! -
\! h}}  \; =   \hfill  \cr  
   \hfill   = \;  {\textstyle \sum\limits_{s=0}^h} \; q^{{{s+1}
\choose 2}} \, {h \choose s}_{\!q} \, \left( {{a \, ; \, u \! + \! s}
\choose {k \! - \! h}} \, d \, - \, q^{u+s-k+h} \, \left( {{a \, ;
\, u \! + \! s \! - \! 2} \choose {k \! - \! h \! - \! 1}} + \, q \,
{{a \, ; \, u \! + \! s \! - \! 1} \choose {k \! - \! h \! - \! 1}}
\right) \, b \, c \right) \; =  \cr   
%
%
   =  \bigg\{ {{a \, ; u} \atop {k \, , h}} \bigg\} \, d - q^{u-k+h}
\, \Bigg( {{a \, ; \, u \! - \! 2} \choose {k \! - \! h \! - \! 1}} \,
+ \, {\textstyle \sum\limits_{s=1}^h} \, q^{{s \choose 2} + s} \left(
\, q^s \, {h \choose s}_{\!q} + \, {h \choose {s \! - \! 1}}_{\!q} \,
\right) \, \times   \hfill  \cr   
   \hfill   \times \, {{a \, ; \, u \! + \! s - \! 2} \choose
{k \! - \! h \! - \! 1}} \, + \, q^{{{h+1} \choose 2} + h + 1} \,
{{a \, ; \, u \! + \! h \! - \! 1} \choose {k \! - \! h \! - \! 1}}
\Bigg) \, b \, c  \; =  \cr   
   \hfill   = \bigg\{\! { {a \, ; u} \atop {k \, , h}} \!\bigg\} \,
d \, - \, q^{u-k+h} {\textstyle \sum\limits_{s=0}^{h+1}} \; q^{{{s+1}
\choose 2}} \, {{h \! + \! 1} \choose s}_{\!q} \, {{a \, ; \, u \! +
\! s \! - \! 2} \choose {k \! - \! h \! - \! 1}} \, b \, c  \, = \, 
\bigg\{\! { {a \, ; u} \atop {k, h}} \!\bigg\} \, d \, - \, q^{u-k+h}
\bigg\{\! { {a \, ; u \! - \! 2} \atop {k \, , h \! + \! 1}} \!\bigg\}
\, b \, c  \cr }  $$   
   \indent   Finally, we assume (4.4) holds for  $ k $  and  $ f \, $, 
\, and we prove it for  $ f \! + \! 1 $  too.  By (4.6) we have   
  $$  \displaylines{ 
   d^{f+1} \, {{a \, ; \, u} \choose k}  \, =  \, d \; d^f \, {{a \, ;
\, u} \choose k}  \, = \,  {\textstyle \sum\limits_{s=0}^{f \wedge
\, k}} \, {(-1)}^s \, q^{{{s+1} \choose 2} \, - \, s \, (f+k-u)} \,
{\bigg[ {f \atop s} \!\,\bigg]}_q \, d \; \bigg\{\! {{a \, ; \, u
\! - \! 2 \, s} \atop {k \, , s}} \!\bigg\} \; b^s \, c^s \, d^{f-s} 
\, =   \hfill  \cr   
   \hfill   = \;  {\textstyle \sum\limits_{s=0}^{f \wedge \, k}} \,
{(-1)}^s \, q^{{{s+1} \choose 2} \, - \, s \, (f+k-u) \, - \, 2 \, s}
\, {\bigg[ {f \atop s} \!\,\bigg]}_q \; \bigg\{ {{a \, ; \, u \! -
\! 2 \, s} \atop {k \, , s}} \bigg\} \; b^s \, c^s \, d^{f+1-s} 
\; -   \hfill \hskip53pt  \cr   
   \hfill   - \;  {\textstyle \sum\limits_{s=0}^{f \wedge \, k}} \,
{(-1)}^s \, q^{{{s+1} \choose 2} \, - \, s \, (f+k-u) \, - \, u \,
- \, s \, - \, k} \, {\bigg[ {f \atop s} \!\,\bigg]}_q \; \bigg\{
{{a \, ; \, u \! - \! 2 \, (s \! + \! 1)} \atop {k \, , s \! +
\! 1}} \bigg\} \; b^{s+1} \, c^{s+1} \, d^{f-s}  \; =  \cr     
   = \!  {\textstyle \sum\limits_{s=0}^{(f+1) \wedge \, k}}
\hskip-7pt {(-1)}^s \, q^{{{s+1} \choose 2} - s \, (f+1+k-u)}
\left(\, q^{-s} {\bigg[ {f \atop s} \!\,\bigg]}_q + \, q^{f+1-s} \,
{\bigg[ {f \atop {s \! - \! 1}} \! \,\bigg]}_q \,\right) \, \bigg\{\!
{{a \, ; \, u \! - \! 2 \, s} \atop {k \, , s}} \!\bigg\} \, b^s \,
c^s \, d^{f+1-s}  \, =   \hfill  \cr   
     \hfill   = \;  {\textstyle \sum\limits_{s=0}^{(f+1) \wedge \, k}}
\, {(-1)}^s \, q^{{{s+1} \choose 2} \, - \, s \, (f+1+k-u)} \, {\bigg[
{{f+1} \atop s} \!\,\bigg]}_q \; \bigg\{ {{a \, ; \, u \! - \! 2 s}
\atop {k \, , s}} \bigg\} \; b^s \, c^s \, d^{f+1-s}  \cr }  $$   
which eventually proves (4.4) for  $ \, f \! + \! 1 \, $  as well.   
                                                  \par   
   The analysis above proves the identity (4.4), which is equivalent
to relation  {\bf [D.4$\boldkey{+}$]}.  The proof of relation 
{\bf [D.4$\boldkey{-}$]}  is entirely similar, just along the
same pattern.   
 \vskip5pt   
   So far we have proved that the given relations do hold in 
$ \calfqm \, $;  \, to prove the first part of the claim, we must
show that these generate the ideal of  {\sl all\/}  relations among
generators.  This amounts to show that the algebra enjoying only the
given relations is isomorphic to  $ \calfqm \, $.  In turn, this is
equivalent to the following: if  $ \Cal{B}' $  is any of the PBW--like
bases of Theorem 4.3, then the given relations are enough to expand
any product of generators as a  $ \Zqqm $--linear  combination of
the monomials in  $ \Cal{B}' \, $.   
   \hbox{But this is clear for the
very relations themselves.}   
                                        \par
   As to the bialgebra structure, everything is just a matter of
computations, via a cute use of Lemma 1.5, and partially basing
upon the case  $ \, n = 2 \, $.  We point out the main steps.   
                                        \par
   Assume  $ \, i < j \, $.  By definitions we have   
  $$  \Delta \left( \t_{i,j}^{(h)} \right) = \, \Delta \Big( \!
{\big( q - q^{-1} \big)}^{-1} \, t_{i,j} \Big)^{(h)} = {\big(
q - q^{-1} \big)}^{-h} \, \Big( {\textstyle \sum_{k=1}^n} \,
t_{i,k} \otimes t_{k,j} \Big)^{(h)} \; .   \eqno (4.7)  $$   
But  $ \, k < k' \, $  implies  $ \, (t_{i,k} \otimes t_{k,j})
\cdot (t_{i,k'} \otimes t_{k',j}) = q^2 \, (t_{i,k'} \otimes
t_{k',j}) \cdot (t_{i,k} \otimes t_{k,j}) \, $,  \, thus 
Lemma 1.5{\it (b)\/}  and (4.7) yield   
  $$  \Delta \left( \t_{i,j}^{(h)} \right) \, = \; {\big( q - q^{-1}
\big)}^{-h} \, {\textstyle \sum\limits_{e_1 + \cdots + e_n = h}}
\hskip-0pt q^{{{e_1 + 1} \choose 2} + \cdots + {{e_n + 1} \choose 2}
- {{h + 1} \choose 2}} \, y_1^{(e_1)} y_2^{(e_2)} \! \cdots y_n^{(e_n)}  
\eqno (4.8)  $$   
for  $ \, y_k := t_{i,k} \otimes t_{k,j} \, $  ($ \, k = 1, \dots,
n \, $).  But now   
  $$  \displaylines{ 
   y_k^{(e_k)} \! = \big( t_{i,k} \otimes t_{k,j} \big)^{(e_k)} \! =
{[e_k]}_q! \cdot t_{i,k}^{(e_k)} \! \otimes t_{k,j}^{(e_k)} = {\big(
q \! - \! q^{-1} \big)}^{2 \, e_k} {[e_k]}_q! \cdot \t_{i,k}^{(e_k)} \!
\otimes \t_{k,j}^{(e_k)},   \hfill   \forall \  k \notin \{i,j\}  \cr   
   y_i^{(e_i)} = \, \big( t_{i,i} \otimes t_{i,j} \big)^{(e_i)} = \,
t_{i,i}^{\,e_i} \! \otimes t_{i,j}^{(e_i)} = 
%
%
  \, {\big( q \! - \! q^{-1} \big)}^{e_i} \!
\cdot {\textstyle \sum\limits_{r=0}^{e_i}} \, q^{r \choose 2}
{e_i \choose r}_{\!q} {(q-1)}^r {(r)}_q! \cdot {{t_{i,i} \, ; \, 0}
\choose r} \otimes \t_{i,j}^{(e_i)} \; ,  \cr   
   y_j^{(e_j)} = \, \big( t_{i,j} \otimes t_{j,j} \big)^{(e_j)} = \,
t_{i,j}^{(e_j)} \! \otimes t_{j,j}^{\,e_j} =   
%
%
  \, {\big( q \! - \! q^{-1} \big)}^{e_j} \!
\cdot {\textstyle \sum\limits_{s=0}^{e_j}} \, q^{s \choose 2}
{e_j \choose s}_{\!q} {(q-1)}^s {(s)}_q! \cdot \t_{i,j}^{(e_j)}
\otimes {{t_{j,j} \, ; \, 0} \choose s} \; .  \cr }  $$  
Using these facts and the commutation relations between generators,
(4.8) turns into the formula for  $ \Delta \left( \t_{i,j}^{(h)}
\right) $  given in the claim.  The case  $ \, i> j \, $  is entirely
similar.  
                                        \par   
   As to  $ \, \Delta \Big( \! {{t_{i,i} \, ; \, c} \atop k} \! \Big)
\, $,  \, by definition of  $ \Delta $  we have  $ \; \Delta \Big( \!
\Big( {{t_{i,i} \, ; \, c} \atop k} \Big) \! \Big) \, = \, \Big(
{{\, \sum_{s=1}^n t_{i,s} \otimes t_{s,i} \; ; \; c \,} \atop k}
\Big) \; $.  To expand the latter term, we use twice Lemma 1.6 and
once Lemma 1.5.  First apply  Lemma 1.6{\it (a-2)\/}  to  $ \; x \,
:= \sum\limits_{h=1}^i t_{i,h} \otimes t_{h,i} \; $,  $ \; w \, :=
\sum\limits_{l=i+1}^n \! \t_{i,l} \otimes \t_{l,i} \; $,  $ \; t =
c \, $,  $ \, m = k \, $;  \; this gives an expansion formula in which
some  $ q $--divided  powers  $ \, w^{(r)} = {\big( \sum_{l=i+1}^n \!
\t_{i,l} \otimes \t_{l,i} \big)}^{(r)} \, $  and some terms  $ \,
\Big\{ \! {{x \; ; \; c} \atop {m \, , \, r}} \! \Big\} \, $  occur. 
Expand the latter ones as a linear combination of some  $ \, \Big( \!
{{x \; ; \; a} \atop b} \! \Big) $'s,  \, and then apply  Lemma
1.6{\it (c-2)\/}  to them with  $ \; t_{i,i} \otimes t_{i,i} \; $ 
in the r\^{o}le of  $ x $  and  $ \; \sum\limits_{h=1}^{i-1} \! \t_{i,h}
\otimes \t_{h,i} \; $  acting as  $ w \, $.  This eventually gives
an expansion formula which is a linear combination of products of 
$ q $--divided  powers  $ \, {\big( \sum_{h=1}^{i-1} \! \t_{i,h}
\otimes \t_{h,i} \big)}^{(s)} \, $,  \,  $ q $--binomial  coefficients 
$ \, \Big( \! {{t_{i,i} \otimes \, t_{i,i} \; ; \; a} \atop b} \! \Big)
\, $  and  $ q $--divided  powers  $ \, {\big( \sum_{l=i+1}^n \! \t_{i,l}
\otimes \t_{l,i} \big)}^{(r)} $,  \, in this order.  Now, the  $ \, \Big(
\! {{t_{i,i} \otimes \, t_{i,i} \; ; \; a} \atop b} \! \Big) $'s  can be
patched together into some  $ \, \Big\{ \! {{t_{i,i} \otimes \, t_{i,i}
\; ; \; \alpha} \atop {\gamma \; , \; \delta}} \! \Big\} $'s,  \, and the 
$ q $--divided  powers  $ \, {\big( \sum_{h=1}^{i-1} \! \t_{i,h} \otimes
\t_{h,i} \big)}^{(s)} \, $  and  $ \, {\big( \sum_{l=i+1}^n \! \t_{i,l}
\otimes \t_{l,i} \big)}^{(r)} \, $  can be expanded using  Lemma 1.5{\it
(b)},  and then also the obvious identities  $ \; {\big( \t_{i,s} \otimes
\t_{s\,i} \big)}^{(e)} = \, {[e]}_q! \cdot \t_{i,s}^{\,(e)} \otimes
\t_{s\,i}^{\,(e)} \; $.  Using the commutation relations between
generators, the final outcome can be written as claimed.   
                                        \par   
   Eventually,  $ \; \epsilon \Big( \! \Big( {{t_{\ell,\ell} \; ; \; c}
\atop k} \Big) \! \Big) \! = \! {\Big( {c \atop k} \Big)}_{\!q} \; $ 
and  $ \; \epsilon \Big( \t_{i,j}^{\,(h)} \Big) \! = 0 \; $  follow
from the very definitions.    
 \vskip4pt   
   {\it (b)} \,  The fact that  $ \calfqgl $  admits the presentation
above is a direct consequence of Theorem 4.3{\it (b)\/}  and of the
presentation of  $ \calfqm $  given in claim  {\it (a)\/}; in particular,
the additional relations in first line simply mean that  $ {D_q}^{-1} $
is central   --- because  $ D_q $  itself is central ---   while the
second line relation is a reformulation of the relation  $ \; D_q \,
{D_q}^{\!-1} = 1 \, $.  The statement about the Hopf structure also
follows from claim  {\it (a)\/}  and  Theorem 4.3{\it (b)\/}  and
from the  formul{\ae}  for the antipode in  $ \fqgl $  (cf.~\S 2.4),
but for the  formul{\ae}  for  $ {D_q}^{\!-1} $  which follow from
%
%
 $ \; \Delta(D_q) \, = \, D_q \otimes D_q \, $,  $ \;
\varepsilon(D_q) \, = \, 1 \, $,  $ \; S\big(D_q\big) \,
= \, {D_q}^{\!-1} \; $.  Clearly, the  formul{\ae}  for the
antipode of generators  $ \, \Big( {{t_{\ell,\ell} \, ; \; c}
\atop k} \Big) \, $  and  $ \, \t_{i,j}^{\,(h)} \, $  are  {\sl
implicit},  in that one should still expand the right-hand sides
of them in terms of the generators themselves.   
                                            \par   
   Now, the  formul{\ae}  for  $ \Delta $  and  $ \epsilon $  do
show that  $ \calfqgl $  is indeed a  $ \Zqqm $--subbialgebra 
of  $ \fqgl \, $.  For the antipode,  $ \, \Big\langle S \big(
\calfqgl \big) \, , \, \caluqgl \Big\rangle = \Big\langle \calfqgl
\, , \, S\big(\caluqgl\big) \! \Big\rangle \subseteq \Zqqm \, $,  \,
which gives  $ \, S \big( \calfqgl \big) \subseteq \calfqgl \; $. 
The claim follows.   
 \vskip4pt   
   {\it (c)} \,  This follows at once from claims  {\it (a)\/}  and 
{\it (b)},  along with Corollary 4.5.   \hskip35pt \hfill \qed\break   
\enddemo

\vskip-5pt

\noindent   
 {\bf 4.7 Remarks:}   
                                    \par    
   {\it (a)} \,  The commutation relations in Theorem 4.6 are not the
only possible ones, but several (equivalent) ones exist.  The given
ones are best suited to expand any product of the generators as a 
$ \Zqqm $--linear  combination of the elements of a PBW basis as
in Theorem 4.3.  
                                    \par    
   {\it (b)} \,  As they belong to  $ \calfqgl \, $,  one could try to
compute an explicit expression for the elements  $ \; S \Big(\! \Big(\!
{{t_{\ell,\ell} \, ; \; c} \atop k} \!\Big) \!\Big) \, $,  \; as well
as for  $ \; S \big( \t_{i,j}^{\,(h)} \big) \, $,  \;  {\sl in terms
of the generators given in Theorem 4.6{\it (b)\/}!}  However, such a
formula might be very complicated.   
                                    \par    
   {\it (c)} \,  In [GR], \S 5.6, some of the ``additional relations''
mentioned in  Theorem 4.6{\it (c)\/}  are found explicitly in case 
$ \, n = 2 \, $.   
                                    \par    
   {\it (d)} \,  Theorem 4.6 also yields presentations at  $ \, q =
1 \, $,  \, which implies the important corollary here below.  As a
consequence, since objects like  $ U_\Z \big({\gergl_n}^{\!*} \big) $ 
are usually called ``hyperalgebras'', the previous result allows us to
call  $ \calfqm \, $,  $ \calfqgl $  and  $ \calfqsl $  {\it ``quantum
hyperalgebras''}.   

\vskip11pt

\proclaim{Theorem 4.8}   
 \vskip1pt   
   {\it (a)} \,  There exists a  $ \, \Z $--bialgebra isomorphism  $ \;
\calF_1[M_n] \cong U_\Z\big({\gergl_n}^{\!*} \big) \, $  given by   
 \vskip-17pt   
  $$  \t_{i,j}^{(h)}{\Big|}_{q=1} \!\! \mapsto
{(-1)}^{h(j-i)} \, \e_{i,j}^{\,(h)} \; ,  \quad 
{{t_{s,s} \, ; \, 0} \choose k}{\bigg|}_{q=1} \!\! \mapsto
{\g_s \choose k} \; ,  \quad  \t_{j,i}^{(h)}{\Big|}_{q=1}
\!\! \mapsto {(-1)}^{h(j-i-1)} \, \f_{j,i}^{\,(h)}  $$   
 \vskip-11pt   
for all  $ \, h $,  $ s $  and  $ \, i < j \, $.  In particular 
$ \, \calF_1[M_n] \, $  is a Hopf  $ \, \Z $--algebra,  isomorphic
to  $ \, U_\Z\big({\gergl_n}^{\!*} \big) \, $.   
 \vskip3pt   
   {\it (b)} \,  There exists a Hopf  $ \, \Z $--algebra
isomorphism  $ \, \calF_1[{GL}_n] \cong U_\Z\big({\gergl_n}^{\!*}
\big) \, $,  \, which is uniquely determined by the formul{\ae} 
in claim (a).  
 \vskip3pt   
   {\it (c)} \,  There exists a Hopf  $ \, \Z $--algebra  isomorphism 
$ \,\calF_1[{SL}_n] \cong U_\Z\big({\gersl_n}^{\!*} \big) $  given
by the same formul{\ae}  as in claim (a).
\endproclaim   

\demo{Proof} The presentation of  $ \calfqm $  given in Theorem
4.5 provides a similar presentation for  $ \calF_1[M_n] \, $;  a
straightforward comparison then shows that the latter is the standard
presentation of  $ \, U_\Z\big({\gergl_n}^{\!*} \big) \, $,  \,
following the correspondence given in the claim.  To be precise in
the first presentation one also has the specialization at  $ \, q
= 1 \, $  of the  $ \Big(\! {{t_{\ell,\ell} \, ; \, r} \atop k}
\!\Big) $'s,  but these are generated by the specializations of
the  $ \Big(\! {{t_{\ell,\ell} \, ; \, 0} \atop \nu} \!\Big) $'s. 
This yields a  $ \Z $--algebra  isomorphism: a moment's check
shows that it is one of  $ \Z $--bialgebras  too.   
                                                \par   
   Up to minimal changes, the proof goes the same for claims 
{\it (b)\/}  and  {\it (c)\/}  as well.   \hskip35pt \hfill \qed\break
\enddemo

%
%
 \eject   

   The previous result regards specialization at  $ \, q = 1 \, $. 
As to specializations at roots of 1, Theorem 3.7 yields the following
important consequence:   

\vskip11pt

\proclaim{Proposition 4.9}  Let  $ \varepsilon $  be a root of unity,
of odd order, and apply notation of\/ \S 1.4.   
 \vskip3pt   
   (a) \,  The specialization  $ \calfem $  is a 
$ \, \Zeps $--bialgebra,  isomorphic to  $ \gerH_\varepsilon^g
\, $  via the specialization of the embedding  $ \; \calfqm
\longhookrightarrow \gerH_q^g \;\, $.   
 \vskip3pt   
   (b) \,  The embedding  $ \; \calfem \longhookrightarrow \calfegl \; $ 
of  $ \, \Zeps $--bial\-gebras  is an isomorphism.  In particular, 
$ \calfem $  and  $ \gerH^g_\varepsilon $  both are Hopf 
$ \Zeps $--algebras  isomorphic to  $ \calfegl \, $.   
 \vskip3pt   
   (c) \,  The specialization  $ \calfesl $  is a Hopf 
$ \, \Zeps $--algebra,  isomorphic to  $ \gerH_\varepsilon^s
\, $  via the specialization of the embedding  $ \, \calfqsl
\lhook\joinrel\rightarrow \gerH_q^s \;\, $.   
\endproclaim   

\demo{Proof}  From Theorem 3.7{\it (a)\/}  one argues  $ \,
\calfem\big[\phi^{-1}\big] = \gerH_\varepsilon^g \; $,  \; which is an
identity induced by the embedding  $ \,\; \widehat{\xi} \, : \calfem
\longhookrightarrow \gerH_\varepsilon^g \; $;  \, note that we can
see that the latter map is injective by looking at how the PBW basis
of  $ \calfqm $  given in  Theorem 4.3{\it (a)\/}  expands w.r.t.~to a
suitable PBW basis of  $ \gerH^g_q \, $,  as we did exactly in the very
proof of  Theorem 4.3{\it (a)}.  Now, for  $ \, \phi := \Lambda_1
\Lambda_2^2 \cdots \Lambda_n^n \, $  one has  $ \, \Big(\! {{\phi \, ;
\, c} \atop k} \!\Big) \in \gerH_q^g \, $  for all  $ \, k \in \N \, $, 
$ \, c \in \Z \, $.  This is easily proved using Theorem 3.1 in [DL];
otherwise, it can also be proved by a brute force computation.  Then   
  $$  \qquad  {{\phi \, ; \, c} \choose k} \, \in \, \gerH_q^g 
\; {\textstyle \bigcap} \, \fqm \, = \, \calfqm  \qquad \qquad 
\big( \, \forall \;\, k \in \N \, ,  \, c \in \Z \,\big)  $$   
Now we apply Lemma 1.8 to  $ \, x = \phi \, $,  \, which gives  $ \,
\phi^{\,\ell} \equiv 1 \mod (q - \varepsilon) \, \calfqm \, $.  Thus 
$ \, \phi^{\,\ell} = 1 \in \calfem \, $,  \, so  $ \, \phi^{-1} =
\phi^{\,\ell - 1} \in \calfem \, $  and so  $ \, \calfem = \calfem
\big[\phi^{-1}\big] = \gerH_\varepsilon^g \; $.   
                                               \par   
   The above proves claim  {\it (a)}.  Now, one can prove the identity 
$ \; \calfegl = \gerH^g_\varepsilon \, $  just like for claim  {\it
(a)},  or one can use the relations  $ \, \calfem \subseteq \calfegl
\subseteq \gerH^g_\varepsilon \; $  and  $ \, \gerH^g_\varepsilon
= \calfem \; $.  As  $ \, \calfegl \, $  is clearly a Hopf 
$ \Zeps $--algebra,  claim  {\it (b)\/}  is then proved.  Similarly, 
{\it (c)\/}  can be proved like  {\it (a)},  or deduced as a corollary
of  {\it (a)\/}  itself, of Corollary 4.5 and of Corollary 3.8. 
\hskip5pt \hfill \qed\break   
\enddemo   

\vskip-6pt

%
%
   Theorem 4.6 and 
Proposition 3.3 also yield 
%
%
 a
description of quantum Frobenius morphisms:  

\vskip11pt

\proclaim{Theorem 4.10}  Let  $ \varepsilon $  be a root of unity, of
odd order  $ \ell \, $.   
                                                \par   
   (a) \,  The quantum Frobenius morphism (2.1) for  $ \calfqm $  is
defined over  $ \Zeps \, $,  i.e.{} it restricts to an epimorphism
of  $ \, \Zeps $--bialgebras
 \vskip-6pt  
  $$  \calFr_{M_n}^{\,\Zeps} \, \colon \, \calfem
\llongtwoheadrightarrow \, \Zeps \otimes_\Z \calF_1[M_n] \,
\cong \, \Zeps \otimes_\Z U_\Z\big({\gergl_n}^{\!*}\big)  $$
 \vskip-1pt  
\noindent   
coinciding, via Theorem 4.8 and Proposition 4.9, with (2.5), and
given on generators by   
  $$  \hskip-0pt   \calFr_{M_n}^{\,\Zeps} \! : \!
 \left\{ \hskip-3pt \hbox{ $ \matrix
   \!\! (\, i \! < \! j \,)  \hskip5pt  \t_{i,j}^{(k)}\Big\vert_{q=\varepsilon} 
\hskip-7pt \mapsto \; \t_{i,j}^{(k/\ell)}\Big\vert_{q=1} \hskip-5pt
= \, {(-1)}^{(j-i) \, k/\ell} \, \e_{i,j}^{\,(k/\ell)}  \hskip3pt  \hbox{ if } \; \ell \Big\vert k \; ,  \hskip9pt \t_{i,j}^{(k)}\Big\vert_{q=\varepsilon} \hskip-5pt \mapsto \; 0 
\hskip5pt \hbox{ if } \; \ell \!\! \not\Big\vert k  \\
   \hskip-3pt  {\displaystyle \bigg( {t_{i,i} \, ; \, 0 \atop k} \bigg)}
\! \bigg\vert_{q=\varepsilon} \hskip-7pt \mapsto {\displaystyle \bigg( {t_{i,i} \, ; \, 0 \atop {k / \ell}} \bigg)} \! \bigg\vert_{q=1} \hskip-7pt = {\displaystyle \bigg(\! {\g_i \atop {k / \ell}} \!\bigg)}  \hskip11pt
\hbox{ if } \; \ell \Big\vert k \; ,  \hskip28pt 
{\displaystyle \bigg( {t_{i,i} \, ; \, 0 \atop n} \bigg)}\bigg\vert_{q=\varepsilon}  \hskip-8pt  \mapsto \, 0  \hskip5pt
\hbox{ if } \; \ell \!\! \not\Big\vert k  \\
   \! (\, i \! > \! j \,)  \hskip5pt  \t_{i,j}^{(k)}\Big\vert_{q=\varepsilon}
\hskip-7pt \mapsto \, \t_{i,j}^{(k/\ell)}\Big\vert_{q=1} \hskip-5pt =
{(-1)}^{(j-i-1) \, k/\ell} \, \f_{i,j}^{\,(k/\ell)}  \hskip3pt
\hbox{ if } \, \ell \Big\vert k \; ,  \hskip9pt \t_{i,j}^{(k)}\Big\vert_{q=\varepsilon} \hskip-7pt \mapsto \, 0 
\hskip3pt \hbox{ if } \, \ell \!\! \not\Big\vert k  
            \endmatrix $ }  \right.  $$
                                                \par   
   (b) \,  The quantum Frobenius morphism (2.2) for  $ \calfqgl $ 
is defined over  $ \Zeps \, $,  i.e.~it restricts to a Hopf 
$ \, \Zeps $--algebra  epimorphism   
 \vskip-6pt  
  $$  \calFr_{GL_n}^{\,\Zeps} \, \colon \, \calfegl
\llongtwoheadrightarrow \, \Zeps \otimes_\Z \calF_1[GL_n] \,
\cong \, \Zeps \otimes_\Z U_\Z\big({\gergl_n}^{\!*}\big)  $$
 \vskip-1pt  
\noindent   
 coinciding, via Theorem 4.8 and Proposition 4.9, with (2.5)
and with  $ \calFr_{M_n}^{\,\Zeps} $  of claim (a), and
described by the same formul{\ae}.    
%
%
 \eject   
   (c) \,  The quantum Frobenius morphism (2.3) for  $ \calfqsl $ 
is defined over  $ \Zeps \, $,  i.e.~it restricts to a Hopf 
$ \, \Zeps $--algebra  epimorphism   
 \vskip6pt  
\noindent   
  $$ \calFr_{{SL}_n}^{\,\Zeps} \, \colon \, \calfesl
\llongtwoheadrightarrow \, \Zeps \otimes_\Z \calF_1[{SL}_n] \,
\cong \, \Zeps \otimes_\Z U_\Z\big({\gersl_n}^{\!*}\big)  $$   
 \vskip-1pt  
\noindent   
coinciding with (2.4) via  Theorem 4.8 and Proposition 4.9,
described by the formul{\ae}  in (a).   
\endproclaim

\demo{Proof} By definition  $ \; \calFr_{M_n}^{\,\Qeps} \, \colon
\, \Qeps \otimes_{\Zeps} \calfem \llongtwoheadrightarrow \, \Qeps
\otimes_\Z U_\Z \big( {\gergl_n}^{\!*} \big) \, = \, \Qeps \otimes_\Q
U \big( {\gergl_n}^{\!*} \big) \; $  is the restriction (via  $ \,
\widehat{\xi} : \calfqm \longhookrightarrow \gerH^g_q \, $  at  $ \,
q = \varepsilon \, $  and  $ \, q = 1 \, $)  of the similar epimorphism 
$ \; \gerFr_{\gergl_n^{\,*}}^{\,\Qeps} \, \colon \, \Qeps \otimes_{\Zeps}
\gerH^g_\varepsilon \, \llongtwoheadrightarrow \, \Qeps \otimes_\Z
\gerH^g_1 \, \cong \, \Qeps \otimes_\Z U_\Z\big({\gergl_n}^{\!*}\big)
\; $  obtained by scalar extension from (2.5), see [Ga1].  Note also
that all the generators  $ \Big( {{t_{i,i} \, ; \, c} \atop k} \Big) $'s 
(for all  $ i $,  $ c $  and  $ k \, $)  are contained in the subalgebra
generated by the  $ \Big( {{t_{j,j} \, ; \, 0} \atop h} \Big) $'s  alone
(for all  $ j $  and  $ k \, $),  thanks to relations (1.2).  Therefore,
the  formul{\ae}  in the claim uniquely determine 
$ \calFr_{M_n}^{\,\Zeps} \, $.   
                                                    \par   
Now, the formul{\ae}  in Lemma 4.2 give, for  $ \, i < j \, $  (noting
that  $ \, {[\ell\,]}_\varepsilon = 0 \, $),   
 \vskip-9pt  
  $$  \displaylines{
   \gerFr_{\gergl_n^{\,*}}^{\,\Qeps} \Big( {\widehat{\xi}\,
\Big|_{q=\varepsilon}} \Big( \t_{i,j}^{\,(k)}\Big|_{q=\varepsilon} \Big)
\Big) \, = \, {\textstyle \sum\limits_{\sum_{h=0}^{n-j} e_h = k}} 
\hskip-9pt  \varepsilon^{e_0 + \sum_{s=1}^{n-j} {e_s \choose 2} +
k(j-i-2) - {k \choose 2}} \, {(-1)}^{k(j-i)} \, {\big( \varepsilon^{-1}
\! - \varepsilon \big)}^{k - e_0} \; \times   \hfill  \cr   
   \times \; {\textstyle \prod\limits_{s=1}^{n-j}}
{[e_s]}_\varepsilon! \, \cdot {\textstyle \prod\limits_{r=0}^{n-j}}
\, \gerFr_{\gergl_n^{\,*}}^{\,\Zeps} \Big(\! E_{i,j+r}^{\,(e_r)}
\Big|_{q=\varepsilon} \Big) \cdot {\textstyle \prod\limits_{r=0}^{n-j}}
\, \gerFr_{\gergl_n^{\,*}}^{\,\Zeps} {\big( \Lambda_{j+r}
\big|_{q=\varepsilon} \big)}^{e_r} \cdot {\textstyle
\prod\limits_{r=1}^{n-j}} \, \gerFr_{\gergl_n^{\,*}}^{\,\Zeps}
\Big(\! F_{j+r,j}^{\,(e_r)}\Big|_{q=\varepsilon} \Big)
\; =  \cr    
   \hfill   = \; \varepsilon^{k(j-i-1) - {k \choose 2}} \,
{(-1)}^{k(j-i)} \cdot \, \gerFr_{\gergl_n^{\,*}}^{\,\Zeps}
\Big(\! E_{i,j}^{\,(k)} \Big|_{q=\varepsilon} \Big)  \cr }  $$   
 \vskip-3pt  
The very last term is equal to  $ \; {(-1)}^{k(j-i)} \,
E_{i,j}^{\,(k/\ell)}{\Big|}_{q=1} = \; {(-1)}^{(j-i) \,
k/\ell} \, \e_{i,j}^{\,(k/\ell)} \, = \; \t_{i,j}^{\,
(k/\ell)}\Big|_{q=1}  \, $  if  $ \, \ell \Big| n \, $, 
\; while it is zero if  $ \, \ell \! \not{\hskip-1pt\Big|}
\, n \, $  (note that in the present case  $ k $  and  $ k/\ell $ 
have the same parity), via  Theorem 4.8{\it (a)}.  This proves the first
formula in the claim.   
                                                       \par   
   Similarly, one proves the third formula.  As for the second formula,
we have   
  $$  \displaylines{
   \gerFr_{\gergl_n^{\,*}}^{\,\Qeps} \!
\left( {\widehat{\xi}\,{\Big|}_{q=\varepsilon}} \!
\left( \! \left(\! {{t_{i,i} \, ; \, 0} \atop k} \!\right)
\bigg\vert_{q=\varepsilon} \right) \!\right) \,  = \;
{\textstyle \sum\limits_{r=0}^k} \, \varepsilon^{-r(k+2)
- {r \choose 2}} \, {\big( \varepsilon - \varepsilon^{-1}
\big)}^r {\textstyle \sum\limits_{\sum_{s=1}^{n-i} e_s = r}}
\hskip3pt {\textstyle \prod\limits_{s=1}^{n-i}} \hskip1pt
\varepsilon^{e_s \choose 2} \, {[e_s]}_\varepsilon! \;
\times   \hfill  \cr   
   \times {\textstyle \prod\limits_{r=1}^{n-i}}
\gerFr_{\gergl_n^{\,*}}^{\,\Zeps} \!\! \left( \!\!
E_{i,i+r}^{\,(e_r)}{\Big|}_{q=\varepsilon} \right)
\, \gerFr_{\gergl_n^{\,*}}^{\,\Zeps} \!\! \left( \!\!
\left\{ {\! {\Lambda_i ; -r} \atop {k \; , \, r} \!\!\!\!\!}
\right\}\!{\bigg|}_{q=\varepsilon} \!\right) \, {\textstyle
\prod\limits_{r=1}^{n-i}} \gerFr_{\gergl_n^{\,*}}^{\,\Zeps}
\! {\big( \Lambda_{i+r}{\big|}_{q=\varepsilon} \big)}^{e_r}
\, {\textstyle \prod\limits_{r=1}^{n-i}}
\gerFr_{\gergl_n^{\,*}}^{\,\Zeps}
\!\! \left( \!\! F_{i+r,i}^{\,(e_r)}
{\Big|}_{q=\varepsilon} \!\right)
\! =  \cr   
   \hfill   = \; \gerFr_{\gergl_n^{\,*}}^{\,\Zeps}
\!\! \left( \! \left\{ {{\Lambda_i \, ; \, 0} \atop {k \, , 0}}
\right\}\!{\bigg|}_{q=\varepsilon} \right) \; = \;
\gerFr_{\gergl_n^{\,*}}^{\,\Zeps} \!\!
\left( \! \left( {{\Lambda_i \, ; \, 0} \atop k} \right)
\!{\bigg|}_{q=\varepsilon} \right) \; = \; \cases
     \hskip-3pt  \Big(\! {{\Lambda_i \, ; \, 0} \atop {k / \ell}}
\!\Big){\Big|}_{q=1}  \hskip-3pt  = \Big(\! {\text{l}_i \atop
{k / \ell}} \!\Big)  &  \text{if} \quad \ell \,\big\vert\, k  \\
     \hskip-3pt  \quad\;  0  &  \text{if} \quad \ell \!\not\big\vert\, k
   \endcases  \cr }  $$
by the very definition of  $ \gerFr_{{\gergl_n}^{\!*}}^{\,\Zeps} $ 
again.  On the other hand, if  $ \, \ell \big\vert\, k \, $  then 
$ \; {\widehat{\xi}\,{\Big|}_{q=1}} \Big( \! \Big(\! {{t_{i,i} \, ;
\, 0} \atop {k / \ell}} \Big)\Big\vert_{q=1} \Big) = \, \Big(\!
{{\Lambda_i \, ; \, 0} \atop {k / \ell}} \Big)\Big\vert_{q=1}
\hskip-1pt = \, \Big(\! {\g_i \atop {k / \ell}} \!\Big) \; $   
%
%
thanks to Theorem 4.8 once more, whence our formula follows.   
                                             \par   
   All this accounts for claim  {\it (a)}.  Claims  {\it (b)\/}  and 
{\it (c)\/}  can be proved with the same arguments, or deduced from 
{\it (a)\/}  in force of Proposition 4.9 and of Corollary 4.5.  
\hskip65pt \hfill \qed\break
\enddemo

%
%
 \vfill   
 \eject    

   \centerline{\sc List of Symbols}   
 \vskip7pt
   \centerline{ $ {(n)}_q \, $,  $ \; {(n)}_q! \, $,  $ \; {\Big(
{n \atop s} \Big)}_{\!q} \, $,  $ \; {\Big( {n \atop {k_1, \dots,
k_r}} \Big)}_{\!q} \, $,  $ \; {[n]}_q \, $  etc.,  $ \; X^{(n)}
\, $,  $ \; \Big( {{X \, ; \, c} \atop n} \Big) \, $,  $ \; \left[
{{X \, ; \, c} \atop n} \right] \, $,  $ \; \left\{ {{X \, ; \, c}
\atop {n \, , \, r}} \right\} \; $:  \;\; see \S 1.3 }   
 \vskip3pt
   \centerline{ $ \uqgl \, $,  $ \; \uqsl \; $:  \;\; see
\S 2.1   \;\;\; --- \;\;\;   $ \geruqg \, $,  $ \caluqg \, $, 
$ \, \gerFr_\gerg^{\,\Zeps} \, $,  $ \, \calFr_\gerg^{\,\Zeps}
\, $  (for  $ \, \gerg \in \big\{ \gergl_n, \gersl_n \big\}\, $): 
\;\; see \S 2.3 }
 \vskip3pt
   \centerline{ $ \fqg \, $,  $ \; \gerfqg \, $,  $ \calfqg \, $, 
$ \, \gerFr_G^{\,\Z} \, $,  $ \, \calFr_G^{\,\Qeps} \, $  (for 
$ \, G \in \big\{ M_n, {GL}_n, {SL}_n \big\}\, $):  \;\; see
\S 2.4 }
 \vskip3pt
   \centerline{ $ \uqgs \, $,  $ \, \geruqgs \, $,  $ \caluqgs \, $, 
$ \H_q^x \, $,  $ \, \gerH_q^x \, $,  $ \calH_q^x \, $, 
$ \, \gerFr_{\gerg^*}^{\,\Zeps} \, $,  $ \, \calFr_{\gerg^*}^{\,\Zeps}
\, $  (for  $ \, x \in \{g,s\} \, $,  $ \, \gerg \in \big\{ \gergl_n,
\gersl_n \big\}\, $):  \;\; see \S 2.5 }
 \vskip3pt
   \centerline{ $ \widetilde{\xi} \, $,  $ \, \xi \, $, 
$ \widehat{\xi} \; $  {\sl (embeddings of quantum function algebras
into dual quantum groups)\/}:  \;\; see \S 3.2 }
  
\vskip1,7truecm

\vskip23pt

\enddocument
\bye
\end